\documentclass{article}%
\usepackage{graphicx}
\usepackage{amsmath}%
\usepackage{amsfonts}%
\usepackage{amssymb}
\newtheorem{theorem}{Theorem}[section]

\newtheorem{corollary}[theorem]{Corollary}

\newtheorem{definition}[theorem]{Definition}

\newtheorem{lemma}[theorem]{Lemma}

\newtheorem{proposition}[theorem]{Proposition}
\newtheorem{remark}[theorem]{Remark}

\newtheorem{theorem and definition}[theorem]{Theorem and Definition}
\newenvironment{proof}[1][Proof]{\textbf{#1.} }{\ \rule{0.5em}{0.5em}}

\begin{document}

\title{Quantum Markov Processes \\(Correspondences and Dilations)}
\author{Paul S. Muhly\thanks{Supported by grants from the U.S. National Science
Foundation and from the U.S.-Israel Binational Science Foundation.}\\Department of Mathematics\\University of Iowa\\Iowa City, IA 52242\\\texttt{muhly@math.uiowa.edu}
\and Baruch Solel\thanks{Supported by the U.S.-Israel Binational Science Foundation
and by the Fund for the Promotion of Research at the Technion}\\Department of Mathematics\\Technion\\32000 Haifa\\Israel\\\texttt{mabaruch@techunix.technion.ac.il}}
\maketitle

\begin{abstract}
We study the structure of quantum Markov Processes from the point of view of
product systems and their representations.

\textbf{2000 Subject Classificiation }Primary: 46L53, 46L55, 46L57, 46L60,
81S25. Secondary: 46L07, 46L08, 47L90.

\end{abstract}

\section{Introduction}

A \emph{quantum Markov process} is a pair, $(\mathcal{M},\{P_{t}\}_{t\geq0})$,
consisting of a von Neumann algebra $\mathcal{M}$ and a semigroup
$\{P_{t}\}_{t\geq0}$ of unital, completely positive, normal linear maps on
$\mathcal{M}$ such that $P_{0}$ is the identity mapping on $\mathcal{M}$ and
such that the map $t\rightarrow P_{t}(a)$ from $[0,\infty)$ to $\mathcal{M}$
is continuous with respect to the $\sigma$-weak topology on $\mathcal{M}$ for
each $a\in\mathcal{M}$. Over the years, there have been numerous studies
wherein the authors ``dilate'' the Markov semigroup $\{P_{t}\}_{t\geq0}$ to an
$E_{0}$-semigroup, in the sense of Arveson \cite{wA89} and Powers \cite{rP88},
of \emph{endomorphisms} $\{\alpha_{t}\}_{t\geq0}$ of a larger von Neumann
algebra $\mathcal{R}$. Depending on context, the process of dilation has taken
different meanings. Here we mean the following: Suppose $\mathcal{M}$ acts on
a Hilbert space $H$, then a quadruple $(K,\mathcal{R},\{\alpha_{t}\}_{t\geq
0},u_{0})$, consisting of a Hilbert space $K$, a von Neumann algebra
$\mathcal{R}$, an $E_{0}$-semigroup $\{\alpha_{t}\}_{t\geq0}$ of $\ast
$-endomorphisms of $\mathcal{R}$ and an isometric embedding $u_{0}%
:H\rightarrow K$, will be called an \emph{E}$_{0}$-\emph{dilation }of the
quantum Markov process $(\mathcal{M},\{P_{t}\}_{t\geq0})$ (or, simply, of
$\{P_{t}\}_{t\geq0}$) in case for all $T\in\mathcal{M}$, all $S\in\mathcal{R}$
and all $t\geq0$ the following equations hold%
\[
P_{t}(T)=u_{0}^{\ast}\alpha_{t}(u_{0}Tu_{0}^{\ast})u_{0}%
\]
and%
\[
P_{t}(u_{0}^{\ast}Su_{0})=u_{0}^{\ast}\alpha_{t}(S)u_{0}\text{.}%
\]
Our objective in this paper is to prove that if the Hilbert space $H$ on which
$\mathcal{M}$ acts is separable, then such a dilation always exists.

What is novel in our approach is that we recognize the space of the
Stinespring dilation of each $P_{t}$ as a correspondence $\mathcal{E}_{t}$
over the \emph{commutant} of $\mathcal{M}$, $\mathcal{M}^{\prime}$. (All
relevant terms will be defined below.) These correspondences are then
assembled and ``dilated'' to a product system $\{E(t)\}_{t\geq0}$ of
correspondences over $\mathcal{M}^{\prime}$, very similar to the product
systems that Arveson defined in \cite{wA89}. Then we describe $\{P_{t}%
\}_{t\geq0}$ explicitly in terms of what we call a ``fully coisometric,
completely contractive covariant representation'' of $\{E(t)\}_{t\geq0}$,
denoted $\{T_{t}\}_{t\geq0}$, in a fashion that derives immediately from our
work in \cite{MS98}. A bit more explicitly, but still incompletely, we find
that $\{P_{t}\}_{t\geq0}$ may be expressed in terms of $\{T_{t}\}_{t\geq0}$
via the formula%
\[
P_{t}(a)=\widetilde{T}_{t}(I_{E(t)}\otimes a)\widetilde{T}_{t}^{\ast}\text{,}%
\]
$a\in\mathcal{M}$ and $t\geq0$, where $\widetilde{T}_{t}$ is the operator from
$E(t)\otimes H$ to $H$ defined by the equation $\widetilde{T}_{t}(\xi\otimes
h)=T_{t}(\xi)h$. Then we dilate $\{T_{t}\}_{t\geq0}$ to what is called an
isometric representation $\{V_{t}\}_{t\geq0}$ of $\{E(t)\}_{t\geq0}$ on a
Hilbert space $K$. If $u_{0}:H\rightarrow K$ is the embedding that goes along
with $\{V_{t}\}_{t\geq0}$, the we find that $T_{t}(\xi)=u_{0}^{\ast}V_{t}%
(\xi)u_{0}$ for all $\xi\in E(t)$ and that the $E_{0}$-semigroup of
endomorphisms $\{\alpha_{t}\}_{t\geq0}$ that we want is given by the formula%
\[
\alpha_{t}(R)=\widetilde{V}_{t}(I_{E(t)}\otimes R)\widetilde{V}_{t}^{\ast
}\text{,}%
\]
where $R$ runs over the von Neumann algebra $\mathcal{R}$ generated by
$\{\alpha_{t}(u_{0}\mathcal{M}u_{0}^{\ast})\}_{t\geq0}$. That is,
$(K,\mathcal{R},\{\alpha_{t}\}_{t\geq0},u_{0})$ is the dilation of
$(\mathcal{M},\{P_{t}\}_{t\geq0})$.

An important part of our analysis was inspired by Bhat's paper \cite{bB96}.
Recently, Bhat and Skeide \cite{BS00} have dilated a quantum Markov process
$(\mathcal{M},\{P_{t}\}_{t\geq0})$ using a product system over the von Neumann
algebra $\mathcal{M}$ (ours is over $\mathcal{M}^{\prime}$). The precise
connection between their work and ours has still to be determined. However,
what we find attractive about our approach is the close explicit connection
between dilations of quantum Markov processes and the classical dilation
theory of contraction operators on Hilbert space pioneered by B. Sz-Nagy (see
\cite{szNF70}).

In the next section we develop the theory of correspondences over von Neumann
algebras sufficiently so that we can link up with theory developed in
\cite{MS98} in which representations and dilations of $C^{\ast}$%
-correspondences are considered. We also show how what we call the Arveson
correspondence $\mathcal{E}_{P}$ associated with the Stinespring dilation of a
completely positive map $P$ can be dilated to a bigger correspondence $E$ in
such a way that a certain completely contractive covariant representation of
$E$ that gives $P$ is dilated to a fully coisometric, isometric representation
of $E$. This representation of $E$ gives a ``power'' dilation of $P$.

Then, in Section 3, we construct a ``discrete'' dilation $(K,\mathcal{R}%
,\{\alpha_{t}\}_{t\geq0},u_{0})$ of the quantum Markov process $(\mathcal{M}%
,\{P_{t}\}_{t\geq0})$. It is here, following ideas developed in Section 2,
that we dilate the family $\{\mathcal{E}_{P_{t}}\}_{t\geq0}$ to a product
system of correspondences $\{E(t)\}_{t\geq0}$ over $\mathcal{M}^{\prime}$. In
Section 4, we show that if the Hilbert space on which $\mathcal{M}$ acts is
separable, then the dilation $(K,\mathcal{R},\{\alpha_{t}\}_{t\geq0},u_{0})$
we construct in Section 2 is, in fact, an $E_{0}$-dilation.

We adopt the standard notation that if $A$ is a subset of a Hilbert space $H$,
then $[A]$ will denote the closed linear span of $A$.

\section{Dilations of Completely Positive Maps\label{section1}}

Throughout, $\mathcal{M}$ will denote a von Neumann algebra. While much of
what we will have to say about von Neumann algebras can be formulated in a
space-free fashion, it will be convenient to view $\mathcal{M}$ as acting on a
fixed Hilbert space $H$. Thus, we will work inside $\mathcal{B}(H)$, the
bounded operators on $H$. Also, throughout, $P$ will denote a fixed completely
positive, unital and normal map from $\mathcal{M}$ to $\mathcal{M}$. We need
to call attention to specific features of the minimal Stinespring dilation of
$P$ \cite{fS55, wA69, wA97}.

Form the algebraic tensor product, $\mathcal{M}\otimes H$ and define the
sesquilinear form $\langle\cdot,\cdot\rangle$ on this space by the formula
\[
\langle T_{1}\otimes h_{1},T_{2}\otimes h_{2}\rangle=\langle h_{1}%
,P(T_{2}^{\ast}T_{1})h_{2}\rangle\text{,}%
\]
$T_{i}\otimes h_{i}\in\mathcal{M}\otimes H$. The complete positivity of $P$
guarantees that this form is positive semidefinite. Therefore, the Hausdorff
completion of $\mathcal{M}\otimes H$ is a Hilbert space, which we shall denote
by $\mathcal{M}\otimes_{P}H$. We shall not distinguish between an element in
$\mathcal{M}\otimes H$ and its image in $\mathcal{M}\otimes_{P}H$. The
formula
\[
\pi_{P}(S)(T\otimes h):=ST\otimes h\text{,}%
\]
$S\in\mathcal{M}$, $T\otimes h\in\mathcal{M}\otimes_{P}H$ defines a
representation of $\mathcal{M}$ on $\mathcal{M}\otimes_{P}H$ that is normal
because $P$ is normal. Also, the formula
\[
W_{P}(h):=I\otimes h,
\]
$h\in H$, defines an isometric imbedding of $H$ in $\mathcal{M}\otimes_{P}H$,
and there results the fundamental equation
\[
P(T)=W_{P}^{\ast}\pi_{P}(T)W_{P}\text{,}%
\]
$T\in\mathcal{M}$. It is an easy matter to check that $\mathcal{M}\otimes
_{P}H$ is minimal in the sense that the smallest subspace of $\mathcal{M}%
\otimes_{P}H$ containing $W_{P}H$ and reducing $\pi_{P}$ is all of
$\mathcal{M}\otimes_{P}H$. Consequently, the triple $(\pi_{P},\mathcal{M}%
\otimes_{P}H,W_{P})$ is the unique minimal triple, $(\pi,K,W)$, up to unitary
equivalence, such that
\[
P(T)=W^{\ast}\pi(T)W\text{,}%
\]
$T\in\mathcal{M}$. We therefore refer to $(\pi_{P},\mathcal{M}\otimes
_{P}H,W_{P})$ as \emph{the} \emph{Stinespring dilation }of $P$.

The adjoint $W_{P}^{\ast}$ of the isometric embedding $W_{P}$ of $H$ in
$\mathcal{M}\otimes_{P}H$ has an explicit form that we will need throughout
our analysis:%
\begin{equation}
W_{P}^{\ast}(X\otimes h)=P(X)h\text{, }X\otimes h\in\mathcal{M}\otimes
H\text{.} \label{Wpstar}%
\end{equation}
This is easy to see because
\[
\langle W_{P}^{\ast}(X\otimes h),k\rangle=\langle X\otimes h,W_{P}%
k\rangle=\langle X\otimes h,I\otimes k\rangle=\langle h,P(X^{\ast}%
)k\rangle=\langle P(X)h,k\rangle.
\]

A space of critical importance for us will be the intertwining space,%

\[
\mathcal{L}_{\mathcal{M}}(H,\mathcal{M}\otimes_{P}H):=\{X:H\rightarrow
\mathcal{M}\otimes_{P}H\mid XT=\pi_{P}(T)X,T\in M\}\text{.}%
\]
That is, $\mathcal{L}_{\mathcal{M}}(H,\mathcal{M}\otimes_{P}H)$ is the space
of operators that intertwine the identity representation of $\mathcal{M}$ on
$H$ and $\pi_{P}$. This space turns out to be a $W^{\ast}$-correspondence over
the \emph{commutant} $\mathcal{M}^{\prime}$ of $\mathcal{M}$. The notion of a
$W^{\ast}$-correspondence is fundamental in this study, and therefore we pause
to develop the terminology and to cite some important facts.

We follow Lance \cite{cL94} for the general theory of Hilbert $C^{\ast}%
$-modules that we shall use. In particular, unless indicated to the contrary,
a Hilbert module $\mathcal{X}$ over a $C^{\ast}$-algebra $A$, will be a
\emph{right }Hilbert $C^{\ast}$-module. We write $\mathcal{L}(\mathcal{X})$
for the space of continuous, adjointable $A$-module maps on $\mathcal{X}$
(which we shall write on the left of $\mathcal{X}$) and we shall write
$\mathcal{K}(\mathcal{X})$ for the space of (generalized) compact operators on
$\mathcal{X}$, i.e., $\mathcal{K}(\mathcal{X})$ is the span of the the rank
one operators $\xi\otimes\eta^{\ast}$, $\xi,\eta\in\mathcal{X}$, where
$\xi\otimes\eta^{\ast}(\zeta)=\xi\langle\eta,\zeta\rangle$.

\begin{definition}
Let $A$ and $B$ be $C^{\ast}$-algebras. A $C^{\ast}$-\emph{correspondence}
\emph{from }$A$\emph{\ to }$B$ is a Hilbert $C^{\ast}$-module $\mathcal{X}$
\emph{over} $B$ endowed with the structure of a left module over $A$ via a
$\ast$-homomorphism $\varphi:A\rightarrow\mathcal{L}(\mathcal{X})$. A
$C^{\ast}$-\emph{correspondence over }$A$ is simply a $C^{\ast}$%
-correspondence from $A$ to $A$.
\end{definition}

When dealing with specific $C^{\ast}$-correspondences, $\mathcal{X}$ from a
$C^{\ast}$-algebra $A$ to a $C^{\ast}$-algebra $B$, it will be convenient to
suppress the $\varphi$ in formulas involving the left action and simply write
$a\xi$ or $a\cdot\xi$ for $\varphi(a)\xi$. \ This should cause no confusion in context.

$C^{\ast}$-correspondences should be viewed as generalized $C^{\ast}%
$-homomorphisms. Indeed, the collection of $C^{\ast}$-algebras together with
(isomorphism classes of) $C^{\ast}$-correspondences is a category that
contains (contravariantly) the category of $C^{\ast}$-algebras and (conjugacy
classes of) $C^{\ast}$-homomorphisms. Of course, for this to make sense, one
has to have a notion of composition of correspondences and a precise notion of
isomorphism. The notion of isomorphism is the obvious one: a bijective,
bimodule map that preserves inner products. Composition is ``tensoring'': If
$\mathcal{X}$ is a $C^{\ast}$-correspondence from $A$ to $B$ and if
$\mathcal{Y}$ is a correspondence from $B$ to $C$, then the balanced tensor
product, $\mathcal{X}\otimes_{B}\mathcal{Y}$ is an $A,C$-bimodule that carries
the inner product defined by the formula
\[
\langle\xi_{1}\otimes\eta_{1},\xi_{2}\otimes\eta_{2}\rangle_{\mathcal{X}%
\otimes_{B}\mathcal{Y}}:=\langle\eta_{1},\varphi(\langle\xi_{1},\xi_{2}%
\rangle_{\mathcal{X}})\eta_{2}\rangle_{\mathcal{Y}}\text{.}%
\]
The Hausdorff completion of this bimodule is again denoted by $\mathcal{X}%
\otimes_{B}\mathcal{Y}$ and is called the \emph{composition} of $\mathcal{X}$
and $\mathcal{Y}$. At the level of correspondences, composition is not
associative. However, if we pass to isomorphism classes, it is. That is, we
only have an isomorphism $(\mathcal{X}\otimes\mathcal{Y})\otimes
\mathcal{Z}\simeq\mathcal{X}\otimes(\mathcal{Y}\otimes\mathcal{Z})$. It is
worthwhile to emphasize here that while it often is safe to ignore the
distinction between correspondences and their isomorphism classes, at times,
as we shall see, the distinction is of critical importance.

If $\mathcal{N}$ is a von Neumann algebra and if $\mathcal{X}$ is a Hilbert
$C^{\ast}$-module over $\mathcal{N}$, then $\mathcal{X}$ is called
\emph{self-dual} in case every continuous $\mathcal{N}$-module map $\Phi$ from
$\mathcal{X}$ to $\mathcal{N}$ is implemented by an element of $\mathcal{X}$,
i.e., in case there is an $\xi_{\Phi}\in\mathcal{X}$ so that $\Phi
(\xi)=\langle\xi_{\Phi},\xi\rangle$, $\xi\in\mathcal{X}$. There is a
topological characterization of self-dual Hilbert $C^{\ast}$-modules over von
Neumann algebras given in \cite{BDH88} that will be useful for us. To state
it, recall that their $\sigma$-topology on a Hilbert $C^{\ast}$-module
$\mathcal{X}$ over a von Neumann algebra $\mathcal{N}$ is the topology defined
by the functionals
\[
f(\cdot):=\sum_{n=1}^{\infty}w_{n}(\langle\eta_{n},\cdot\rangle)
\]
where the $\eta_{n}$ lie in $\mathcal{X}$, the $w_{n}$ lie in $\mathcal{N}%
_{\ast}$, and $\sum\left\|  w_{n}\right\|  \left\|  \eta_{n}\right\|  <\infty
$. Baillet, Denizeau, and Havet proved that a Hilbert $C^{\ast}$-module
$\mathcal{X}$ over a von Neumann algebra $\mathcal{N}$ is self-dual if and
only if the unit ball in $\mathcal{X}$ is compact in the $\sigma$-topology
\cite[Proposition 1.7]{BDH88}. In \cite{wP73}, Paschke proved that if
$\mathcal{X}$ is a self-dual Hilbert $C^{\ast}$-module over a von Neumann
algebra $\mathcal{N}$, then $\mathcal{L}(\mathcal{X})$ is a von Neumann
algebra, i.e., $\mathcal{L}(\mathcal{X})$ is a $C^{\ast}$-algebra which is
also a dual space and which, therefore, may be represented faithfully on
Hilbert space in such a way that the weak-$\ast$ topology on $\mathcal{L}%
(\mathcal{X})$ coincides with the $\sigma$-weak topology on the image.

\begin{definition}
Let $\mathcal{M}$ and $\mathcal{N}$ be von Neumann algebras and let
$\mathcal{X}$ be a Hilbert $C^{\ast}$-module over $\mathcal{N}$. Then
$\mathcal{X}$ is called a \emph{Hilbert }$W^{\ast}$\emph{-module} over
$\mathcal{N}$ in case $\mathcal{X}$ is self-dual. The module $\mathcal{X}$ is
called a $W^{\ast}$\emph{-correspondence from }$\mathcal{M}$\emph{\ to
}$\mathcal{N}$\emph{\ }in case $\mathcal{X}$ is a self-dual $C^{\ast}%
$-correspondence from $\mathcal{M}$ to $\mathcal{N}$ such that the $\ast
$-homomorphism $\varphi:\mathcal{M}\rightarrow\mathcal{L}(\mathcal{X})$ giving
the left module structure on $\mathcal{X}$ is normal.
\end{definition}

It is evident that the composition of $W^{\ast}$-correspondences is again a
$W^{\ast}$-correspondence. The following proposition shows one way to
construct $W^{\ast}$-correspondences.

\begin{proposition}
\label{Lemma 1.2}Let $H$ and $K$ be Hilbert spaces. Let $\mathcal{M}$ be a von
Neumann algebra on $K$, let $\mathcal{N}$ be a von Neumann algebra on $H$, and
let $\mathcal{Y}\subseteq B(H,K)$ be a $\sigma$-weakly closed linear space of
operators. Suppose that $\mathcal{MYN}\subseteq\mathcal{Y}$ and that
$\mathcal{Y}^{\ast}\mathcal{Y}:=\{Y^{\ast}Y\mid Y\in\mathcal{Y}\}$ is
contained in $\mathcal{N}$. Then $\mathcal{Y}$ is a self-dual Hilbert
$W^{\ast}$-module over $\mathcal{N}$ that has the structure of a $W^{\ast}%
$-correspondence from $\mathcal{M}$ to $\mathcal{N}$. Further, the $\sigma
$-topology on $\mathcal{Y}$ coincides with the $\sigma$-weak topology on
$\mathcal{Y}$ as a subspace of $B(H,K)$.
\end{proposition}

\begin{proof}
It is evident that $\mathcal{Y}$ has the structure of a $C^{\ast}%
$-correspondence from $\mathcal{M}$ to $\mathcal{N}$. The main point of the
proposition is the assertion about self-duality and the topologies. The
functionals $f$ defining the $\sigma$-topology are of the form $f(\cdot
):=\sum_{n=1}^{\infty}w_{n}(\langle\eta_{n},\cdot\rangle)$ where the $\eta
_{n}$ lie in $\mathcal{Y}$, the $w_{n}$ lie in $\mathcal{N}_{\ast}$, and
$\sum\left\|  w_{n}\right\|  \left\|  \eta_{n}\right\|  <\infty$. Evidently,
each of these is $\sigma$-weakly continuous. Conversely, given a functional on
$\mathcal{Y}$ of the form $g(Y)=\langle Yh,k\rangle$, we may assume that $k$
is in the closed span of $\{Yh\mid Y\in\mathcal{Y}$, $h\in H\}$ and
approximate $g$ in norm by functionals of the form
\[
\tilde{g}(Y)=\sum_{m=1}^{r}\langle Yh,Y_{m}h_{m}^{\prime}\rangle=\sum
_{m=1}^{r}\langle Y^{\ast}Yh_{m},h_{m}^{\prime}\rangle.
\]
Each of these functionals is continuous in the $\sigma$-topology. Since the
space of $\sigma$-continuous functionals is a Banach space \cite[1.2]{BDH88},
the functional $Y\rightarrow\langle Yh,k\rangle$ is in this space, and so is
every $\sigma$-weakly continuous functional on $\mathcal{Y}$. It follows that
the two topologies coincide on $\mathcal{Y}$. Since the closed unit ball in
$\mathcal{Y}$ is $\sigma$-weakly compact, it must be compact in the $\sigma
$-topology. By \cite[Proposition 1.7]{BDH88}, $\mathcal{Y}$ is self-dual.
\end{proof}

\begin{remark}
(i) The theory developed in \cite{dB97} can be used to prove a converse to
this result: Given a $W^{\ast}$-correspondence $\mathcal{Y}$ from
$\mathcal{M}$ to $\mathcal{N}$, then there are faithful normal representations
$\pi:\mathcal{M}\rightarrow B(K)$ and $\rho:\mathcal{N}\rightarrow B(H)$ and
there is a linear map $\Phi:\mathcal{Y}\rightarrow B(H,K)$ such that
$\Phi(\varphi(T)YS)=\pi(T)\Phi(Y)\rho(S)$ and $\rho(\langle X,Y\rangle
_{\mathcal{Y}})=\Phi(X)^{\ast}\Phi(Y)$ for all $X,Y\in\mathcal{Y}$,
$T\in\mathcal{M}$, and $S\in\mathcal{N}$, and such that $\Phi$ is a
homeomorphism with respect to the $\sigma$-topology on $\mathcal{Y}$ and the
$\sigma$-weak topology on $\Phi(\mathcal{Y})$. Thus, in a sense, the
construction in Proposition \ref{Lemma 1.2} is universal.

(ii) Suppose $\mathcal{X}$ is a self-dual Hilbert $W^{\ast}$-module over a von
Neumann algebra $\mathcal{N}$ and that $\pi:\mathcal{M}\rightarrow
\mathcal{L}(\mathcal{X})$ is a $C^{\ast}$-homomorphism on the von Neumann
algebra $\mathcal{M}$. Then $\pi$ is normal if for every bounded net
$\{A_{\alpha}\}\subseteq\mathcal{N}$, with $A_{\alpha}\rightarrow A$ weakly,
every $g\in\mathcal{N}_{\ast}$, and every $X,Y\in\mathcal{X}$, we have
$g(\langle\pi(A_{\alpha})X,Y\rangle)\rightarrow g(\langle\pi(A)X,Y\rangle)$.
This follows from the fact that $\mathcal{L}(\mathcal{X})$ is the dual space
of the tensor product $\mathcal{X}\otimes\mathcal{X}^{\ast}\otimes
\mathcal{N}_{\ast}$ equipped with the greatest cross norm \cite[Proposition
3.10]{wP73}.
\end{remark}

\begin{proposition}
\label{Lemma1.3}Let $(\pi_{P},\mathcal{M}\otimes_{P}H,W_{P})$ be the
Stinespring dilation of a completely positive map $P$ on the von Neumann
algebra $\mathcal{M}$. Then $\mathcal{L}_{\mathcal{M}}(H,\mathcal{M}%
\otimes_{P}H)$ is a $\sigma$-weakly closed subspace of $B(H,\mathcal{M}%
\otimes_{P}H)$ that is closed under left multiplication by $I\otimes
\mathcal{M}^{\prime}$ and right multiplication by $\mathcal{M}^{\prime}$ and
has the property that $\mathcal{L}_{\mathcal{M}}(H,\mathcal{M}\otimes
_{P}H)^{\ast}\mathcal{L}_{\mathcal{M}}(H,\mathcal{M}\otimes_{P}H)\subseteq
\mathcal{M}^{\prime}$. Thus $\mathcal{L}_{\mathcal{M}}(H,\mathcal{M}%
\otimes_{P}H)$ has the structure of a $W^{\ast}$-correspondence over
$\mathcal{M}^{\prime}$.
\end{proposition}

\begin{proof}
Evidently, if $X\in\mathcal{L}_{\mathcal{M}}(H,\mathcal{M}\otimes_{P}H)$ and
$T\in\mathcal{M}^{\prime}$, then $XT\in\mathcal{L}_{\mathcal{M}}%
(H,\mathcal{M}\otimes_{P}H)$. Indeed, if $S\in\mathcal{M}$, then
$XTS=XST=\pi_{P}(S)XT$, showing that $XT\in\mathcal{L}_{\mathcal{M}%
}(H,\mathcal{M}\otimes_{P}H)$. Also, if $X,Y\in\mathcal{L}_{\mathcal{M}%
}(H,\mathcal{M}\otimes_{P}H)$, then $X^{\ast}Y\in\mathcal{M}^{\prime}$ because
for all $T\in\mathcal{M}$, $X^{\ast}YT=X^{\ast}\pi_{P}(T)Y=TX^{\ast}Y$. Thus,
by Proposition \ref{Lemma 1.2}, it remains to give the left action of
$\mathcal{M}^{\prime}$. On the face of it, this is evidently given by the
formula $\varphi(T)=I\otimes T$, $T\in\mathcal{M}^{\prime}$. However, the
meaning of $I\otimes T$, $T\in\mathcal{M}^{\prime}$, and the expression
$I\otimes\mathcal{M}^{\prime}$, need a little development. For $T\in
\mathcal{M}^{\prime}$, we prove that the algebraic tensor product $I\otimes T$
is bounded as follows: Observe that for $\sum_{i=1}^{n}S_{i}\otimes h_{i}\in
M\otimes H$, we have
\begin{align*}
\left\|  (I\otimes T)(\sum_{i=1}^{n}S_{i}\otimes h_{i})\right\|  ^{2}  &
=\langle I\otimes T(\sum_{i=1}^{n}S_{i}\otimes h_{i}),I\otimes T(\sum
_{i=1}^{n}S_{i}\otimes h_{i})\rangle\\
&  =\langle(\sum_{i=1}^{n}S_{i}\otimes Th_{i}),(\sum_{i=1}^{n}S_{i}\otimes
Th_{i})\rangle\\
&  =\sum_{i,j=1}^{n}\langle Th_{i},P(S_{i}^{\ast}S_{j})Th_{j}\rangle\\
&  =\sum_{i,j=1}^{n}\langle h_{i},T^{\ast}P(S_{i}^{\ast}S_{j})Th_{j}%
\rangle\text{.}%
\end{align*}
However, since $P$ is completely positive the operator matrix $\left(
P(S_{i}^{\ast}S_{j})\right)  $ is a positive element in $M_{n}(\mathcal{M}%
^{\prime}\mathbb{)}$ and so can be written as $C^{\ast}C$, for an element
$C\in M_{n}(\mathcal{M}^{\prime}\mathbb{)}$. Therefore, $\left(  T^{\ast
}P(S_{i}^{\ast}S_{j})T\right)  =\hat{T}^{\ast}C^{\ast}C\hat{T}=C^{\ast}\hat
{T}^{\ast}\hat{T}C\leq\left\|  T\right\|  ^{2}C^{\ast}C$, where $\hat{T}$ is
the $n$-fold inflation of $T$. Consequently, the last term in the displayed
equation is dominated by
\[
\left\|  T\right\|  ^{2}\sum_{i,j=1}^{n}\langle h_{i},T^{\ast}P(S_{i}^{\ast
}S_{j})Th_{j}\rangle=\left\|  T\right\|  ^{2}\left\|  (\sum_{i=1}^{n}%
S_{i}\otimes h_{i})\right\|  ^{2}\text{.}%
\]
Thus $I\otimes T$ extends to an element in $\pi_{P}(\mathcal{M})^{\prime}$,
which we continue to denote by $I\otimes T$. The collection of all these
operators on $M\otimes_{P}H$ is denoted by $I\otimes\mathcal{M}^{\prime}$.
Evidently, the map $T\rightarrow I\otimes T$ is a (not-necessarily-injective)
normal $\ast$-homomorphism of $\mathcal{M}^{\prime}$ onto its range.
Nevertheless, we denote the range by $I\otimes\mathcal{M}^{\prime}$ and note
that $\mathcal{L}_{\mathcal{M}}(H,\mathcal{M}\otimes_{P}H)$ is a left
$\mathcal{M}^{\prime}$ module through this homomorphism. Since $\mathcal{L}%
_{\mathcal{M}}(H,\mathcal{M}\otimes_{P}H)$ is manifestly $\sigma$-weakly
closed, the proof is completed by appeal to Proposition \ref{Lemma 1.2}.
\end{proof}

For our purposes, a drawback of $\mathcal{L}_{\mathcal{M}}(H,\mathcal{M}%
\otimes_{P}H)$ is that it is a space of operators acting between two
\emph{different} Hilbert spaces. We want to ``pull $\mathcal{L}_{\mathcal{M}%
}(H,\mathcal{M}\otimes_{P}H)$ back'' to $H$ using $W_{P}$ and the following
device that is due to Arveson \cite{wA97}. Given $Y\in B(H)$, and $S\otimes h$
in the algebraic tensor product $\mathcal{M}\otimes H$, we set
\[
\Phi_{Y}(S\otimes h):=SY^{\ast}h\text{,}%
\]
and extend $\Phi_{Y}$ by linearity to a densely defined linear map from
$\mathcal{M}\otimes_{P}H$ to $H$ with domain $\mathcal{M}\otimes H$. We write
$\mathcal{E}_{P}$ for the space of all operators $Y\in B(H)$ such that
$\Phi_{Y}$ is continuous. (In this case, of course, we continue to write
$\Phi_{Y}$ for the unique continuous extension to all of $\mathcal{M}%
\otimes_{P}H$.)

\begin{proposition}
\label{Lemma1.4}The space $\mathcal{E}_{P}$ is a linear space that is stable
under left and right multiplication by elements from $\mathcal{M}^{\prime}$,
and the pairing $\langle Y,Z\rangle:=\Phi_{Y}\Phi_{Z}^{\ast}$ converts
$\mathcal{E}_{P}$ into a $W^{\ast}$-correspondence that is isomorphic to
$\mathcal{L}_{\mathcal{M}}(H,\mathcal{M}\otimes_{P}H)$ under the map
$Y\mapsto\Phi_{Y}^{\ast}$.
\end{proposition}

The advantage of $\mathcal{E}_{P}$ over $\mathcal{L}_{\mathcal{M}%
}(H,\mathcal{M}\otimes_{P}H)$ is not only that $\mathcal{E}_{P}\subseteq
B(H)$, but also, as we shall see shortly, given two completely positive maps
$P$ and $Q$ on $\mathcal{M}$, the relations among $\mathcal{E}_{P}$,
$\mathcal{E}_{Q}$, and $\mathcal{E}_{PQ}$ are easier to work with than those
among $\mathcal{L}_{\mathcal{M}}(H,\mathcal{M}\otimes_{P}H)$, $\mathcal{L}%
_{\mathcal{M}}(H,\mathcal{M}\otimes_{Q}H)$, and $\mathcal{L}_{\mathcal{M}%
}(H,\mathcal{M}\otimes_{PQ}H)$.\medskip

\begin{proof}
Evidently, $\mathcal{E}_{P}$ is a linear space. For $R\in\mathcal{M}^{\prime}$
and $Y\in\mathcal{E}_{P}$, $\Phi_{YR}=R^{\ast}\Phi_{Y}$, and so $\mathcal{E}%
_{P}\mathcal{M}^{\prime}\subseteq\mathcal{E}_{P}$. For the other side, fix
$Y\in\mathcal{E}_{P}$ and $T\in\mathcal{M}^{\prime}$. Then for $(\sum
_{i=1}^{n}S_{i}\otimes h_{i})\in\mathcal{M}\otimes H$, $\Phi_{TY}(\sum
_{i=1}^{n}S_{i}\otimes h_{i})=\sum_{i=1}^{n}S_{i}Y^{\ast}T^{\ast}h_{i}%
=\Phi_{Y}(\sum_{i=1}^{n}S_{i}\otimes T^{\ast}h_{i})$. Therefore $\Phi
_{TY}=\Phi_{Y}(I\otimes T^{\ast})$ on $\mathcal{M}\otimes H$, showing that
$\Phi_{TY}$ is bounded; i.e., $\mathcal{M}^{\prime}\mathcal{E}_{P}%
\subseteq\mathcal{E}_{P}$.

Next note that for $Y\in\mathcal{E}_{P}$, $S,T\in\mathcal{M}$, and $h\in H$,
\[
\Phi_{Y}(\pi_{P}(T)(S\otimes h))=\Phi_{Y}(TS\otimes h)=TSY^{\ast}h=T\Phi
_{Y}(S\otimes h)\text{; }%
\]
i.e., $\Phi_{Y}\pi_{P}(T)=T\Phi_{Y}$. Taking adjoints, we conclude that
$\Phi_{Y}^{\ast}\in\mathcal{L}_{\mathcal{M}}(H,\mathcal{M}\otimes_{P}H)$.
Thus, from the properties of $\mathcal{L}_{\mathcal{M}}(H,\mathcal{M}%
\otimes_{P}H)$, we know that the formula $\langle Z,Y\rangle:=\Phi_{Z}\Phi
_{Y}^{\ast}$, defines an $\mathcal{M}^{\prime}$-valued function on
$\mathcal{E}_{P}$. In fact, it is clearly an $\mathcal{M}^{\prime}$-valued
sesquilinear form, since $\Phi_{YR}=R^{\ast}\Phi_{Y}$ and $\Phi_{TY}=\Phi
_{Y}(I\otimes T^{\ast})$, $R,T\in\mathcal{M}^{\prime}$, and it is clearly
positive semidefinite. It is definite because if $\langle Y,Y\rangle=\Phi
_{Y}\Phi_{Y}^{\ast}=0$, then $\Phi_{Y}=0,$ and so $TY^{\ast}h=0$ for all
$T\in\mathcal{M}$ and $h\in H$. Taking $T=I$ we conclude that $Y=0$.

The map $Y\rightarrow\Phi_{Y}^{\ast}$ preserves inner products by definition.
Further, it is a bimodule map since $\left(  \Phi_{YR}\right)  ^{\ast
}=(R^{\ast}\Phi_{Y})^{\ast}=\Phi_{Y}^{\ast}R$ and $\left(  \Phi_{RY}\right)
^{\ast}=\left(  \Phi_{Y}(I\otimes R^{\ast})\right)  ^{\ast}=(I\otimes
R)\Phi_{Y}^{\ast}$, for all $R\in\mathcal{M}^{\prime}$. Thus, to show that
$\mathcal{E}_{P}$ is isomorphic to $\mathcal{L}_{\mathcal{M}}(H,\mathcal{M}%
\otimes_{P}H)$ under this map, we need only show that it is onto. However, we
assert that for all $X\in\mathcal{L}_{\mathcal{M}}(H,\mathcal{M}\otimes_{P}%
H)$, $X=\Phi_{(W_{P}^{\ast}X)}^{\ast}$. Indeed, for $h,h^{\prime}\in H$, and
$S\in\mathcal{M}$, the fact that $(I\otimes S)X=XS$ implies that
\begin{align*}
\Phi_{(W_{P}^{\ast}X)}(S\otimes h)=SX^{\ast}W_{P}h=SX^{\ast}(I\otimes h)\\
=X^{\ast}(S\otimes I)(I\otimes h)=X^{\ast}(S\otimes h)\text{.}%
\end{align*}
This shows that $(W_{P}^{\ast}X)$ is in $\mathcal{E}_{P}$ and that
$X=\Phi_{(W_{P}^{\ast}X)}^{\ast}$. The facts that $\mathcal{E}_{P}$ and
$\mathcal{L}_{\mathcal{M}}(H,\mathcal{M}\otimes_{P}H)$ are isomorphic and that
$\mathcal{L}_{\mathcal{M}}(H,\mathcal{M}\otimes_{P}H)$ is self-dual imply that
$\mathcal{E}_{P}$ is self-dual.
\end{proof}

\begin{definition}
\label{arvesoncorresp}The $W^{\ast}$-correspondence $\mathcal{E}_{P}$ over
$\mathcal{M}^{\prime}$ associated with a normal, unital completely positive
map $P$ on a von Neumann algebra $\mathcal{M}$ will be called the
\emph{Arveson correspondence }associated with $P$.
\end{definition}

The following corollary is immediate from the self-duality of the spaces
involved. We call attention to it because it will be used several times in the sequel.

\begin{corollary}
\label{corollary1.5}If a subspace $\mathcal{Y}$ of $\mathcal{E}_{P}$ or of
$\mathcal{L}_{\mathcal{M}}(H,\mathcal{M}\otimes_{P}H)$ has zero annihilator,
i.e., if $\mathcal{Y}^{\perp}=0$, then $\mathcal{Y}$ is dense.
\end{corollary}

A special case of our analysis so far needs to be singled out. \ Suppose that
$P=\alpha$ is a unital, normal, $\ast$-\emph{endomorphism} of $\mathcal{M}$.
(All endomorphisms will be unital, normal, and preserve adjoints.) Then
$\mathcal{M}\otimes_{\alpha}H$ is isomorphic to $H$ under the map $T\otimes
h\rightarrow\alpha(T)h$, which is $W_{\alpha}^{-1}=W_{\alpha}^{\ast}$.
Further, $\pi_{\alpha}$ is unitarily equivalent to $\alpha$. Thus, we may
identify $\mathcal{E}_{P}=\mathcal{E}_{\alpha}$ directly with $\mathcal{L}%
_{\mathcal{M}}(H,\mathcal{M}\otimes_{P}H)$ as in the following corollary of
Proposition \ref{Lemma1.4}.

\begin{corollary}
\label{corollary1.6}If $\alpha$ is a unital, normal, $\ast$-endomorphism of
$\mathcal{M}$, then $\mathcal{E}_{\alpha}=\{X\in B(H)\mid XT=\alpha
(T)X,\;T\in\mathcal{M}\}$ with the inner product $\langle X_{1},X_{2}%
\rangle=X_{1}^{\ast}X_{2}$.
\end{corollary}

The next lemma may seem like a technicality, but among other things, it
establishes that the modules $\mathcal{L}_{\mathcal{M}}(H,\mathcal{M}%
\otimes_{P}H)$ and $\mathcal{E}_{P}$ are nonzero. It plays other useful roles
in the sequel.

\begin{lemma}
\label{lemma1.7}In the setting of a normal, unital completely positive map $P$
on $\mathcal{M}$ that we have been studying,
\[
\mathcal{M}\otimes_{P}H=\bigvee\{\Phi_{Y}^{\ast}(H)\mid Y\in\mathcal{E}%
_{P}\}=\bigvee\{X(H)\mid X\in\mathcal{L}_{\mathcal{M}}(H,\mathcal{M}%
\otimes_{P}H)\}\text{.}%
\]
\end{lemma}

\begin{proof}
Since $\pi_{P}$ is a normal $\ast$-representation, its kernel is of the form
$\mathcal{M}q$ for a central projection $q$. Write $\pi^{\prime}$ for the
representation of $\mathcal{M}$ that is reduction by $I-q$, i.e., $\pi
^{\prime}(S)=S(I-q)$. Then for $S\in\mathcal{M}$, $\left\|  \pi^{\prime
}(S)\right\|  =\left\|  \pi_{P}(S)\right\|  $, so that $\pi^{\prime}$ and
$\pi_{P}$ are quasiequivalent. If $Q$ is the projection of $\mathcal{M}%
\otimes_{P}H$ onto $\bigvee\{X(H)\mid X\in\mathcal{L}_{\mathcal{M}%
}(H,\mathcal{M}\otimes_{P}H)\}$, then for every $L\in\mathcal{L}_{\mathcal{M}%
}(H,\mathcal{M}\otimes_{P}H)$, $L(H)$ is $\pi_{P}(\mathcal{M})$-invariant.
Hence $Q\in\pi_{P}(\mathcal{M})^{\prime}$. If $\pi_{0}$ is the reduction of
$\pi_{P}$ to the range of $I-Q$, then, on the one hand, $\pi_{0}\leq\pi_{P}$
and on the other, $\pi_{0}$ is disjoint from $\pi^{\prime}$. Since
$\pi^{\prime}$ is quasiequivalent to $\pi_{P}$ we conclude that $\pi_{0}=0$,
i.e., that $Q=I$.
\end{proof}

We next want to illuminate the relation between the composition of two
completely positive maps on $\mathcal{M}$ and the composition of their Arveson
correspondences. This was worked out in the case when $\mathcal{M}=B(H)$ by
Arveson in \cite[Theorem 1.12]{wA97}. Given two normal, unital completely
positive maps $P,Q:\mathcal{M}\rightarrow\mathcal{M}$, we shall write $m$ for
the multiplication map from $\mathcal{E}_{P}\otimes_{\mathcal{M}^{\prime}%
}\mathcal{E}_{Q}$ to $B(H)$. That is, $m(Y\otimes Z)=YZ$.

\begin{lemma}
\label{lemma1.8}The range of $m$ is contained in $\mathcal{E}_{PQ}$.
\end{lemma}

\begin{proof}
First observe that if $Y\in\mathcal{E}_{P}$, if $A=(a_{ij})$ is a positive
semidefinite element in $M_{n}(\mathcal{M}\mathbb{)}$, and if $\mathbf{h}%
=(h_{1},h_{2},\ldots,h_{n})$ is an $n$-tuple of elements from $H$, then
\begin{equation}
\left\|  A^{1/2}(Y^{\ast}\otimes I)\mathbf{h}\right\|  ^{2}\leq\left\|
\Phi_{Y}\right\|  ^{2}\langle h_{i},\sum P(a_{ij})h_{j}\rangle\text{.}
\label{inflation}%
\end{equation}
To see this, note that if $A$ is a diad, i.e., if $A$ has the form
\[
A=(S_{1},S_{2},\ldots,S_{n})^{\ast}(S_{1},S_{2},\ldots,S_{n}),\;\;S_{i}%
\in\mathcal{M}\text{,}%
\]
then the left hand side of the inequality is simply $\left\|  \sum
S_{i}Y^{\ast}h_{i}\right\|  ^{2}$, while the right hand side is $\left\|
\Phi_{Y}\right\|  ^{2}\left\|  \sum S_{i}\otimes h_{i}\right\|  ^{2}$. So, the
inequality is valid by definition. However, every non-negative $A\in
M_{n}(\mathcal{M}\mathbb{)}$ is a sum of diads. (See \cite[Lemma 3.11]{vP86}.)
So the inequality is valid as claimed.

Now fix $Y\in\mathcal{E}_{P}$, $Z\in\mathcal{E}_{Q}$, $\sum S_{i}\otimes
h_{i}\in\mathcal{M}\otimes H$. Then
\begin{align*}
\left\|  \sum S_{i}(YZ)^{\ast}h_{i}\right\|  ^{2}  &  =\left\|  \sum
S_{i}Z^{\ast}(Y^{\ast}h_{i})\right\|  ^{2}\leq\left\|  \Phi_{Z}\right\|
^{2}\left\|  \sum S_{i}\otimes_{Q}Y^{\ast}h_{i}\right\|  ^{2}\\
&  =\left\|  \Phi_{Z}\right\|  ^{2}\langle\sum Y^{\ast}h_{i},Q(S_{i}^{\ast
}S_{j})Y^{\ast}h_{j}\rangle\\
&  =\left\|  \Phi_{Z}\right\|  ^{2}\left\|  (Q(S_{i}^{\ast}S_{j}%
))^{1/2}(Y^{\ast}\otimes I)\mathbf{h}\right\|  ^{2}\\
&  \leq\left\|  \Phi_{Z}\right\|  ^{2}\left\|  \Phi_{Y}\right\|  ^{2}%
\langle\sum h_{i},P(Q(S_{i}^{\ast}S_{j}))Y^{\ast}h_{j}\rangle\\
&  =\left\|  \Phi_{Z}\right\|  ^{2}\left\|  \Phi_{Y}\right\|  ^{2}\left\|
S\otimes_{PQ}h\right\|  ^{2}\text{.}%
\end{align*}
This shows that $\Phi_{YZ}$ is bounded, and that $\left\|  \Phi_{YZ}\right\|
\leq\left\|  \Phi_{Z}\right\|  \left\|  \Phi_{Y}\right\|  $. Thus
$YZ\in\mathcal{E}_{PQ}$.
\end{proof}

We also want to express $m$ in terms of the space $\mathcal{L}_{\mathcal{M}%
}(H,\mathcal{M}\otimes_{P}H)$. For this purpose, fix two normal, unital,
completely positive maps $P,Q:\mathcal{M}\rightarrow\mathcal{M}$. We define a
map $\Psi:\mathcal{L}_{\mathcal{M}}(H,\mathcal{M}\otimes_{P}H)\otimes
_{\mathcal{M}^{\prime}}\mathcal{L}_{\mathcal{M}}(H,\mathcal{M}\otimes
_{Q}H)\rightarrow\mathcal{L}_{\mathcal{M}}(H,\mathcal{M}\otimes_{Q}%
\mathcal{M}\otimes_{P}H)$ by the formula $\Psi(X\otimes Y)=(I\otimes X)Y$,
where $I\otimes X$ is the map from $\mathcal{M}\otimes_{Q}H$ to $\mathcal{M}%
\otimes_{Q}\mathcal{M}\otimes_{P}H$ given by the equation $I\otimes X(S\otimes
h)=S\otimes Xh$. We also define a map $V_{0}:\mathcal{M}\otimes_{PQ}%
H\rightarrow\mathcal{M}\otimes_{Q}\mathcal{M}\otimes_{P}H$ via the equation
$V_{0}(S\otimes h)=S\otimes I\otimes h$, and we define $V:\mathcal{L}%
_{\mathcal{M}}(H,\mathcal{M}\otimes_{P}H)\rightarrow\mathcal{L}_{\mathcal{M}%
}(H,\mathcal{M}\otimes_{Q}\mathcal{M}\otimes_{P}H)$ by the formula
$V(X)=V_{0}X$.

\begin{proposition}
\label{Lemma1.9}In the notation just established, $\Psi$ is an isomorphism of
correspondences and $V$ is an isometry whose range is $\{X\in\mathcal{L}%
_{\mathcal{M}}(H,\mathcal{M}\otimes_{Q}\mathcal{M}\otimes_{P}H)\mid
X(H)\subseteq\mathcal{M}\otimes_{Q}I\otimes_{P}H\}$. Further, if we write
$U_{P}:\mathcal{E}_{P}\rightarrow\mathcal{L}_{\mathcal{M}}(H,\mathcal{M}%
\otimes_{P}H)$ for the isomorphism defined above, and similarly write $U_{Q}$
and $U_{PQ}$, then
\begin{equation}
U_{PQ}m(U_{P}^{-1}\otimes U_{Q}^{-1})=V^{\ast}\Psi\text{,} \label{comult}%
\end{equation}
showing that $m$ is coisometric and $m^{\ast}$ is isometric.
\end{proposition}

\begin{proof}
On the one hand, $\Psi(X\otimes Y)^{\ast}\Psi(X^{\prime}\otimes Y^{\prime
})=Y^{\ast}(I\otimes X^{\ast})(I\otimes X^{\prime})Y^{\prime}=Y^{\ast
}(I\otimes X^{\ast}X^{\prime})Y^{\prime}$ - the inner product in
$\mathcal{L}_{\mathcal{M}}(H,\mathcal{M}\otimes_{Q}\mathcal{M}\otimes_{P}H)$.
On the other hand, recall that the left action of $Z\in\mathcal{M}^{\prime}$
on $\mathcal{L}_{\mathcal{M}}(H,\mathcal{M}\otimes_{Q}H)$ is given by the
equation $(I\otimes Z)Y$, $Y\in\mathcal{L}_{\mathcal{M}}(H,\mathcal{M}%
\otimes_{Q}H)$. Consequently, in $\mathcal{L}_{\mathcal{M}}(H,\mathcal{M}%
\otimes_{P}H)\otimes_{\mathcal{M}^{\prime}}\mathcal{L}_{\mathcal{M}%
}(H,\mathcal{M}\otimes_{Q}H)$,
\begin{align*}
\langle X\otimes Y,X^{\prime}\otimes Y^{\prime}\rangle &  =\langle Y,\langle
X,X^{\prime}\rangle_{\mathcal{L}_{\mathcal{M}}(H,\mathcal{M}\otimes_{P}%
H)}Y^{\prime}\rangle_{\mathcal{L}_{\mathcal{M}}(H,\mathcal{M}\otimes_{Q}H)}\\
&  =\langle Y,(I\otimes X^{\ast}X^{\prime})Y^{\prime}\rangle_{\mathcal{L}%
_{\mathcal{M}}(H,\mathcal{M}\otimes_{Q}H)}=Y^{\ast}(I\otimes X^{\ast}%
X^{\prime})Y^{\prime}\text{.}%
\end{align*}
Thus $\Psi$ preserves the inner products.

To see that $\Psi$ is a bimodule map, let $S\in\mathcal{M}^{\prime}$,
$X\in\mathcal{L}_{\mathcal{M}}(H,\mathcal{M}\otimes_{P}H)$, and $Y\in
\mathcal{L}_{\mathcal{M}}(H,\mathcal{M}\otimes_{Q}H)$. Then
\begin{align*}
\Psi(S(X\otimes Y))  &  =\Psi((I\otimes S)X\otimes Y)\\
&  =(I\otimes(I\otimes S)X)Y\\
&  =(I\otimes I\otimes S)(I\otimes X)Y\\
&  =S\Psi(X\otimes Y)\text{,}%
\end{align*}
while
\[
\Psi(X\otimes YS)=(I\otimes X)(YS)=((I\otimes X)Y)S=\Psi(X\otimes Y)S\text{.}%
\]

To see that $\Psi$ is surjective, we use the fact that $\mathcal{L}%
_{\mathcal{M}}(H,\mathcal{M\otimes}_{Q}\mathcal{M}\otimes_{P}H)$ is self-dual
(see Proposition \ref{Lemma 1.2}) and show that $(\operatorname{Im}%
\Psi)^{\perp}=\{0\}$. Corollary \ref{corollary1.5}, then, will yield the
result. If $Z$ annihilates $\operatorname{Im}\Psi$, then for every
$X\in\mathcal{L}_{\mathcal{M}}(H,\mathcal{M}\otimes_{P}H)$ and for every
$Y\in\mathcal{L}_{\mathcal{M}}(H,\mathcal{M}\otimes_{Q}H)$, $Y^{\ast}(I\otimes
X^{\ast})Z=0$. Observe that $(I\otimes X^{\ast})Z$ is a map from $H$ to
$\mathcal{M}\otimes_{Q}H$. By Lemma \ref{lemma1.7}, $\mathcal{M}\otimes_{Q}H$
is the span of $Y(H)$, $Y\in\mathcal{L}_{\mathcal{M}}(H,\mathcal{M}\otimes
_{Q}H)$. Consequently, $\cap\{\ker Y^{\ast}\mid Y\in\mathcal{L}_{\mathcal{M}%
}(H,\mathcal{M}\otimes_{Q}H)\}=\{0\}$. Since $Y^{\ast}(I\otimes X^{\ast})Z=0$
for all such $Y$, we conclude that $(I\otimes X^{\ast})Z=0$ for all
$X\in\mathcal{L}_{\mathcal{M}}(H,\mathcal{M}\otimes_{P}H)$. By Lemma
\ref{lemma1.7} again, $\cap\{\ker X^{\ast}\mid X\in\mathcal{L}_{\mathcal{M}%
}(H,\mathcal{M}\otimes_{P}H)\}=\{0\}$, and so $Z=0$.

Turning now to $V$, note that it is an easy matter to check that $V_{0}$ is an
isometry. Consequently, for $X\in\mathcal{L}_{\mathcal{M}}(H,\mathcal{M}%
\otimes_{PQ}H)$, $V(X)^{\ast}V(X)=X^{\ast}V_{0}^{\ast}V_{0}X=X^{\ast}X$. So
$V$ is isometric. Also, it is evident from the definitions that $V$ is a
bimodule map and that its image is $\{X\in\mathcal{L}_{\mathcal{M}%
}(H,\mathcal{M}\otimes_{Q}\mathcal{M}\otimes_{P}H)\mid X(H)\subseteq
\mathcal{M}\otimes_{Q}I\otimes_{P}H\}$. Thus, we are left to prove equation
(\ref{comult}). To this end, let $Y\in\mathcal{E}_{P}$ and let $Z\in
\mathcal{E}_{Q}$. As an equation between maps from $\mathcal{M}\otimes_{PQ}H$
to $H$, the relation $\Phi_{Z}(I\otimes\Phi_{Y})V_{0}=\Phi_{YZ}$ is immediate.
Therefore, $\Phi_{YZ}^{\ast}=V_{0}^{\ast}(I\otimes\Phi_{Y}^{\ast})\Phi
_{Z}^{\ast}$. Since $U_{P}(Y)=\Phi_{Y}^{\ast},\;U_{Q}(Z)=\Phi_{Z}^{\ast}\;$and
$m(Y\otimes Z)=YZ$, we find that
\begin{align*}
U_{PQ}m(U_{P}^{-1}\otimes U_{Q}^{-1})(\Phi_{Y}^{\ast}\otimes\Phi_{Z}^{\ast})
&  =U_{PQ}(YZ)=\Phi_{YZ}^{\ast}\\
&  =V_{0}^{\ast}\Psi(\Phi_{Y}^{\ast}\otimes\Phi_{Z}^{\ast})=V^{\ast}(\Psi
(\Phi_{Y}^{\ast}\otimes\Phi_{Z}^{\ast}))\text{.}%
\end{align*}
\end{proof}

\begin{corollary}
\label{Cor 1.10}Under the hypotheses of \ Proposition \ref{Lemma1.9},
\[
\bigvee\{T(H)\mid T\in\mathcal{L}_{\mathcal{M}}(H,\mathcal{M}\otimes
_{Q}\mathcal{M}\otimes_{P}H)\}=\mathcal{M}\otimes_{Q}\mathcal{M}\otimes_{P}H.
\]
\end{corollary}

\begin{proof}
The span $\bigvee\{T(H)\mid T\in\mathcal{L}_{\mathcal{M}}(H,\mathcal{M}%
\otimes_{Q}\mathcal{M}\otimes_{P}H)\}$ contains the span $\bigvee\{(I\otimes
X)Y(H)\mid X\in\mathcal{L}_{\mathcal{M}}(H,\mathcal{M}\otimes_{P}H)$,
$Y\in\mathcal{L}_{\mathcal{M}}(H,\mathcal{M}\otimes_{Q}H)\}$. But using Lemma
\ref{lemma1.7} twice, we see that this space is $\bigvee\{(I\otimes
X)(M\otimes_{Q}H)\mid X\in\mathcal{L}_{\mathcal{M}}(H,\mathcal{M}\otimes
_{P}H)\}=\mathcal{M}\otimes_{Q}\mathcal{M}\otimes_{P}H$.
\end{proof}

Although the product of two correspondences associated with completely
positive maps does not coincide with the correspondence of the product of the
maps, there are important special situations when they do. This, and more, is
spelled out in the following proposition. (For a related result, see
\cite[Theorem 2.12]{BDH88}.)

\begin{proposition}
\label{Proposition1.11}(1) If $\alpha\in Aut(\mathcal{M})$, then
$\mathcal{E}_{\alpha^{-1}}=\mathcal{E}_{\alpha}^{\ast}$, and $\mathcal{E}%
_{\alpha}\mathcal{E}_{\alpha}^{\ast}$ and $\mathcal{E}_{\alpha}^{\ast
}\mathcal{E}_{\alpha}$ are $\sigma$-weakly dense in $\mathcal{M}^{\prime}$. In
particular, $\mathcal{E}_{\alpha}$ is an $\mathcal{M}^{\prime}-\mathcal{M}%
^{\prime}$ equivalence bimodule.

(2) If $\alpha$ is an \emph{endomorphism }of $\mathcal{M}$ and if $Q$ is a
normal, unital, completely positive map of $\mathcal{M}$, then the
multiplication map $m:\mathcal{E}_{\alpha}\otimes_{\mathcal{M}^{\prime}%
}\mathcal{E}_{Q}\rightarrow\mathcal{E}_{\alpha\circ Q}$, defined above, is an isomorphism.

(3) If $\alpha\in Aut(\mathcal{M})$ and if $P$ is a normal, unital completely
positive map, then the map $m:\mathcal{E}_{P}\otimes_{\mathcal{M}^{\prime}%
}\mathcal{E}_{\alpha}\rightarrow\mathcal{E}_{P\circ\alpha}$ is an isomorphism.

(4) If $P$ and $Q$ are conjugate normal, unital completely positive maps
(i.e., if there is an automorphism $\alpha$ of $\mathcal{M}$ such that
$P=\alpha\circ Q\circ\alpha^{-1}$, then
\begin{equation}
\mathcal{E}_{P}\simeq\mathcal{E}_{\alpha}\otimes_{\mathcal{M}^{\prime}%
}\mathcal{E}_{Q}\otimes_{\mathcal{M}^{\prime}}\mathcal{E}_{\alpha}^{\ast
}\text{.}\label{abstracttensor}%
\end{equation}
In particular, $(\mathcal{M}^{\prime},\mathcal{E}_{P})$ and $(\mathcal{M}%
^{\prime},\mathcal{E}_{Q})$ are strongly Morita equivalent in the sense of
\cite{MS00}.
\end{proposition}

We note for the sake of emphasis, that the tensor products described in
equation (\ref{abstracttensor}) are realized through operator multiplication
and adjunction. That is, if $R$ and $T\in\mathcal{E}_{\alpha}$ and
$S\in\mathcal{E}_{Q}$, then $T^{\ast}\in\mathcal{E}_{\alpha}^{\ast}$ and
$R\otimes S\otimes T^{\ast}=RST^{\ast}$.\medskip

\begin{proof}
(1) By Corollary \ref{corollary1.6}, $\mathcal{E}_{\alpha}=\{X\in B(H)\mid
XT=\alpha(T)X,\;T\in\mathcal{M}\}$. The inner product is given by the formula
$\langle X_{1},X_{2}\rangle=X_{1}^{\ast}X_{2}$. It follows easily that
$\mathcal{E}_{\alpha^{-1}}=\mathcal{E}_{\alpha}^{\ast}$ and that
$\mathcal{E}_{\alpha}^{\ast}\mathcal{E}_{\alpha}\subseteq\mathcal{M}^{\prime}%
$. In fact, the $\sigma$-weak closure of $\mathcal{E}_{\alpha}^{\ast
}\mathcal{E}_{\alpha}$ is a $2$-sided ideal in $\mathcal{M}^{\prime}$, and
therefore is of the form $q\mathcal{M}^{\prime}$, for some central projection
in $\mathcal{M}^{\prime}$. However, given $X\in\mathcal{E}_{\alpha^{-1}}$,
with polar decomposition $X=V|X|$, we see that $V\in\mathcal{E}_{\alpha^{-1}}$
and $VV^{\ast}$ is the projection onto the closure of the range of $X$. Since
$VV^{\ast}\in\mathcal{E}_{\alpha^{-1}}\mathcal{E}_{\alpha^{-1}}^{\ast
}=\mathcal{E}_{\alpha}^{\ast}\mathcal{E}_{\alpha}\subseteq q\mathcal{M}%
^{\prime}$, the range of $X$ is contained in the range of $q$. \ However,
Lemma \ref{lemma1.7} implies, now, that $q=I$; i.e., that $\mathcal{E}%
_{\alpha}^{\ast}\mathcal{E}_{\alpha}$ is $\sigma$-weakly dense in
$\mathcal{M}^{\prime}$. Hence, $\mathcal{E}_{\alpha}$ is a normal equivalence bimodule.

(2) From Proposition \ref{Lemma1.9}, we know in general that $m$ is an
isomorphism if $\mathcal{M}\otimes_{Q}I\otimes_{P}H=\mathcal{M}\otimes
_{Q}\mathcal{M}\otimes_{P}H$. If $P=\alpha$ is an endomorphism of
$\mathcal{M}$, then for every $T,S\in\mathcal{M}$ and $h,k\in H$, we have
\begin{align*}
\langle T\otimes h-I\otimes\alpha(T)h,S\otimes k\rangle_{\mathcal{M}%
\otimes_{\alpha}H}  &  =\langle T\otimes h,S\otimes k\rangle-\langle
I\otimes\alpha(T)h,S\otimes k\rangle\\
&  =\langle h,\alpha(T^{\ast}S)k\rangle-\langle\alpha(T)h,\alpha(S)k\rangle=0.
\end{align*}
Hence $T\otimes h=I\otimes\alpha(T)h$ and so $\mathcal{M}\otimes_{\alpha
}H=I\otimes H=H$. Thus, $\mathcal{M}\otimes_{Q}\mathcal{M}\otimes_{\alpha
}H=\mathcal{M}\otimes_{Q}I\otimes_{\alpha}H$, as required.

(3) As in (2), we need to show that $\mathcal{M}\otimes_{\alpha}I\otimes_{Q}H$
coincides with $\mathcal{M}\otimes_{\alpha}\mathcal{M}\otimes_{Q}H$. This is
obvious, since for all $S,T\in\mathcal{M}$, and $h\in H$, $S\otimes T\otimes
h=S\alpha^{-1}(T)\otimes I\otimes h$.

(4) From (2) and (3), we know that $\mathcal{E}_{P}\simeq\mathcal{E}_{\alpha
}\otimes\mathcal{E}_{Q}\otimes\mathcal{E}_{\alpha}^{\ast}$, and (1) implies
then that $\mathcal{E}_{P}\otimes\mathcal{E}_{\alpha}\simeq\mathcal{E}%
_{\alpha}\otimes\mathcal{E}_{Q}$. Part (1) also asserts that $\mathcal{E}_{Q}$
is an $\mathcal{M}^{\prime}$-$\mathcal{M}^{\prime}$-equivalence bimodule and
so $\mathcal{E}_{P}$ are strongly Morita equivalent in the sense of
\cite{MS00}.
\end{proof}

Although $\mathcal{E}_{P}\otimes\mathcal{E}_{Q}\ncong\mathcal{E}_{PQ}$, and
especially, $\mathcal{E}_{P}^{\otimes n}\ncong\mathcal{E}_{P^{n}}$ in general,
we shall soon see that it is possible to ``dilate'' $\mathcal{E}_{P}$ to a
correspondence $\mathcal{F}_{\alpha}$ where $\alpha$ is an endomorphism of the
commutant an isomorphic copy of $\mathcal{M}^{\prime}$. We then have
$\mathcal{F}_{\alpha}^{\otimes n}\simeq\mathcal{F}_{\alpha^{n}}$, by part (2)
of Proposition \ref{Proposition1.11}. This $\alpha$, then, will turn out to be
a ``dilation'' of $P$. To effect this program, we require some of the
technology from \cite{MS98}. We generally adopt the terminology and notation
of \cite{MS98}, but with some minor modifications because we are working in
the category of von Neumann algebras and \emph{normal }maps - representations,
and completely positive maps.

\begin{definition}
\label{Definition1.12}Let $\mathcal{E}$ be a $W^{\ast}$-correspondence over a
von Neumann algebra $\mathcal{N}$ and let $H_{0}$ be a Hilbert space.

\begin{enumerate}
\item A \emph{completely contractive covariant representation }of
$\mathcal{E}$ in $B(H_{0})$ is a pair $(T,\sigma)$, where

\begin{enumerate}
\item $\sigma$ is a normal $\ast$-representation of $\mathcal{N}$ in
$B(H_{0})$.

\item $T$ is a linear, completely contractive map from $\mathcal{E}$ to
$B(H_{0})$ that is continuous in the $\sigma$-topology of \cite{BDH88} on
$\mathcal{E}$ and the $\sigma$-weak topology on $B(H_{0}).$

\item $T$ is a bimodule map in the sense that $T(S\xi R)=\sigma(S)T(\xi
)\sigma(R)$, $\xi\in\mathcal{E}$, and $S,R\in\mathcal{N}$.
\end{enumerate}

\item A completely contractive covariant representation $(T,\sigma)$ of
$\mathcal{E}$ in $B(H_{0})$ is called \emph{isometric }in case
\begin{equation}
T(\xi)^{\ast}T(\eta)=\sigma(\langle\xi,\eta\rangle)\text{,} \label{isometric}%
\end{equation}
for all $\xi,\eta\in\mathcal{E}$.
\end{enumerate}
\end{definition}

To lighten the terminology, we shall refer to an isometric, completely
contractive, covariant representation simply as an isometric covariant
representation. There is no problem doing this because it is easy to see that
if one has a pair $(T,\sigma)$ satisfying all the conditions of part 1 of
Definition \ref{Definition1.12}, except possibly the complete contractivity
assumption, but which is isometric in the sense of equation (\ref{isometric}),
then necessarily $T$ is completely contractive.

The theory developed in \cite{MS98} applies here to prove that if a completely
contractive covariant representation, $(T,\sigma)$, of $\mathcal{E}$ in $B(H)$
is given, then it determines a contraction $\tilde{T}:\mathcal{E}%
\otimes_{\sigma}H\rightarrow H$ defined by the formula $\tilde{T}(\eta\otimes
h):=T(\eta)h$, $\eta\otimes h\in\mathcal{E}\otimes_{\sigma}H$. Here,
$\mathcal{E}\otimes_{\sigma}H$ denotes the Hausdorff completion of the
algebraic tensor product $\mathcal{E}\otimes H$ in the pre-inner product given
by the formula $\langle\xi\otimes h,\eta\otimes k\rangle:=\langle
h,\sigma(\langle\xi,\eta\rangle)k\rangle$. (See \cite[Lemma 3.5]{MS98}.) Also,
there is an \emph{induced }representation $\sigma^{\mathcal{E}}:\mathcal{L}%
(\mathcal{E})\rightarrow B(\mathcal{E}\otimes_{\sigma}H)$ defined by the
formula $\sigma^{\mathcal{E}}(S):=S\otimes I$ \cite[Lemma 3.4]{MS98}.
Recalling that $\mathcal{L}(\mathcal{E})$ is a von Neumann algebra, it is not
hard to see that $\sigma^{\mathcal{E}}$ is a normal representation. The
operator $\tilde{T}$ and $\sigma^{\mathcal{E}}$ are related by the equation%
\begin{equation}
\tilde{T}\sigma^{\mathcal{E}}\circ\varphi=\sigma\tilde{T}. \label{covariance}%
\end{equation}
In fact we have the following lemma that is immediate from \cite{MS98} and
\cite{MS99}. See, in particular, \cite[Lemmas 3.4-3.6]{MS98} and \cite[Lemma
2.1]{MS99}.

\begin{lemma}
\label{CovRep}The map $(T,\sigma)\rightarrow\tilde{T}$ is a bijection between
all completely contractive covariant representations $(T,\sigma)$ of
$\mathcal{E}$ on the Hilbert space $H$ and contractive operators $\tilde
{T}:\mathcal{E}\otimes_{\sigma}H\rightarrow H$ that satisfy equation
(\ref{covariance}). Given such a $\tilde{T}$ satisfying this equation, $T$,
defined by the formula $T(\xi)h:=\tilde{T}(\xi\otimes h)$, together with
$\sigma$ is a completely contractive covariant representation of $\mathcal{E}$
on $H$. Further, $(T,\sigma)$ is isometric if and only if $\tilde{T}$ is an isometry.
\end{lemma}

We note in passing that this lemma shows that the $\sigma$-weak continuity of
$T$ really depends only on the fact that $\sigma$ is normal.

The map $\Psi:\mathcal{L}(\mathcal{E})\rightarrow B(H)$ defined, then, by the
formula
\[
\Psi(S):=\tilde{T}\sigma^{\mathcal{E}}(S)\tilde{T}^{\ast}\text{,}%
\]
$S\in\mathcal{L}(\mathcal{E})$, evidently is completely positive, normal, and contractive.

\begin{definition}
\label{Definition1.12bis}Given a completely contractive covariant
representation $(T,\sigma)$ of $\mathcal{E}$ in $B(H),$ the map $\Psi$ is
called the \emph{completely positive extension} of $(T,\sigma)$, and the
representation $(T,\sigma)$ is called \emph{fully coisometric} in case
$\Psi(I_{\mathcal{E}})=I_{H}$.
\end{definition}

The terminology is reminiscent of the theory of a single contraction. A
completely contractive covariant representation $(T,\sigma)$ is isometric
precisely when $\tilde{T}$ is an isometry. Likewise, it is fully coisometric
precisely when $\tilde{T}$ is a coisometry. The map $\Psi$ is a normal $\ast
$-representation precisely when $(T,\sigma)$ is isometric and it is a unital
$\ast$-representation precisely when $(T,\sigma)$ is both isometric and fully
coisometric. (We have, however, resisted the temptation to call $(T,\sigma)$
unitary in this case.) Our next result, which is a variant of \cite[Corollary
5.21]{MS98}, shows that a completely contractive covariant representation
$(T,\sigma)$ can be dilated to an isometric covariant representation in the
following sense.

\begin{theorem and definition}
\label{Theorem 1.13}Let $\mathcal{E}$ be a $W^{\ast}$-correspondence over a
von Neumann algebra $\mathcal{N}$ and let $(T,\sigma)$ be a completely
contractive covariant representation of $\mathcal{E}$ on the Hilbert space
$H$. Then there is a Hilbert space $K$ containing $H$ and an isometric
covariant representation $(V,\rho)$ of $\mathcal{E}$ on $K$ such that if $P$
is the projection of $K$ onto $H$, then

\begin{enumerate}
\item $P$ commutes with $\rho(\mathcal{N})$ and $\rho(A)P=\sigma(A)P,$
$A\in\mathcal{N}$; and

\item for all $\eta\in\mathcal{E}$, $V(\eta)^{\ast}$ leaves $H$ invariant and
$PV(\eta)P=T(\eta)P$.
\end{enumerate}

The representation $(V,\rho)$ may be chosen so that the smallest subspace $K$
containing $H$ that reduces $(V,\rho)$ is $K$. When this is done, $(V,\rho)$
is unique up to unitary equivalence and is called \emph{the minimal isometric
dilation of }$(T,\sigma)$.

Further, if $(T,\sigma)$ is fully coisometric, the (unique minimal) isometric
dilation $(V,\rho)$ is fully coisometric, too.
\end{theorem and definition}

\begin{proof}
One can construct a proof following the steps leading to Theorem 3.3 and
Corollary 5.21 in \cite{MS98}. However, continuity issues must be dealt with
along the way and one needs to observe that the ideal $J$ discussed there
plays no role here. Rather than doing this, it is easier and it may be more
revealing to appeal to Lemma \ref{CovRep} and simply write down the operator
$\tilde{V}$ and representation $\rho$ that lead to the dilation $(V,\rho)$ of
$(T,\sigma)$. The remaining details will be very easy to verify.

To this end, let $\Delta=(I-\tilde{T}^{\ast}\tilde{T})^{1/2}$ and let
$\mathcal{D}$ be its range. Then $\Delta$ is an operator on $\mathcal{E}%
\otimes_{\sigma}H$ and commutes with the representation $\sigma^{\mathcal{E}%
}\circ\varphi$ of $\mathcal{N}$, by equation (\ref{covariance}). Write
$\sigma_{1}$ for the restriction of $\sigma^{\mathcal{E}}\circ\varphi$ to
$\mathcal{D}$. Let $\sigma_{2}=\sigma_{1}^{\mathcal{E}}\circ\varphi$ on
$\mathcal{E}\otimes_{\sigma_{1}}\mathcal{D}$, and let $\sigma_{3}=\sigma
_{2}^{\mathcal{E}}\circ\varphi$ on $\mathcal{E}\otimes_{\sigma_{2}%
}(\mathcal{E}\otimes_{\sigma_{1}}\mathcal{D})$. It is easy to see that
$\sigma_{3}$ is naturally unitarily equivalent to $\sigma_{1}^{\mathcal{E}%
^{\otimes2}}\circ\varphi_{2}$ on $\mathcal{E}^{\otimes2}\otimes_{\sigma_{1}%
}\mathcal{D}$, where $\varphi_{2}$ is the representation of $\mathcal{N}$ in
$\mathcal{L}(\mathcal{E}^{\otimes2})$ defined by the formula $\varphi
_{2}(a)(\xi\otimes\eta)=(\varphi(a)\xi)\otimes\eta$. We shall identify them
henceforth and in general, we write $\sigma_{n+1}$for $\sigma_{1}%
^{\mathcal{E}^{\otimes n}}\circ\varphi_{n}$ on $\mathcal{E}^{\otimes n}%
\otimes_{\sigma_{1}}\mathcal{D}$, where $\varphi_{n}$ has its obvious meaning.
It is evident that all the $\sigma_{n}$ are normal. We let
\[
K=H\oplus\mathcal{D}\oplus\sum_{n=1}^{\infty}\oplus\mathcal{E}^{\otimes
n}\otimes_{\sigma_{1}}\mathcal{D}%
\]
and we let $\rho=\sigma\oplus\sigma_{1}\oplus\bigoplus_{n=1}^{\infty}%
\sigma_{n+1}$, i.e., thinking matricially, $\rho=diag(\sigma,\sigma_{1}%
,\sigma_{2},\ldots)$. Then a moment's reflection reveals that $\rho$ is a
normal representation of $\mathcal{N}$ on $K$ whose restriction to $H$ is
$\sigma$, of course. Form $\mathcal{E}\otimes_{\rho}K$ and define $\tilde
{V}:\mathcal{E}\otimes_{\rho}K\rightarrow K$ matricially as
\[
\left[
\begin{array}
[c]{cccccc}%
\tilde{T} & 0 & 0 & \cdots &  & \\
\Delta & 0 & 0 &  & \ddots & \\
0 & I & 0 & \ddots &  & \\
0 & 0 & I & 0 & \ddots & \\
\vdots & 0 & 0 & I & \ddots & \\
&  &  &  & \ddots & \ddots
\end{array}
\right]  .
\]
Of course the identity operators in this matrix really must be interpreted as
the operators that identify $\mathcal{E}\otimes_{\sigma_{n+1}}(\mathcal{E}%
^{\otimes n}\otimes_{\sigma_{1}}\mathcal{D})$ with $\mathcal{E}^{\otimes
(n+1)}\otimes_{\sigma_{1}}\mathcal{D}$.

It is easily checked that $\tilde{V}$ is an isometry and that the associated
covariant representation $(V,\rho)$ is an isometric dilation $(T,\sigma)$.
Moreover, it is easily checked that $(V,\rho)$ is minimal, i.e., that the
smallest subspace of $K$ containing $H$ and reducing $(V,\rho)$ is $K$.
Further, if $(T,\sigma)$ is fully coisometric, so that $\tilde{T}$ is a
coisometry, then so is $\tilde{V}$ a coisometry and $(V,\rho)$ is\ fully coisometric.

The proof of the uniqueness of $(V,\rho)$ is the same as in the $C^{\ast}%
$-setting and is given in \cite[Proposition 3.2]{MS98}.

Finally, to see that $V$ is fully coisometric if $T$ is, observe that if $T$
is fully coisometric, then $\widetilde{T}$ is a coisometry as we noted
earlier. Thus $\widetilde{T}\Delta^{2}=0$. This implies that $\widetilde
{T}\Delta=0$. Therefore, from the form of $\widetilde{V}$, we see that
$\widetilde{V}\widetilde{V}^{\ast}=I$, which proves that $V$ is fully coisometric.
\end{proof}

We shall use Theorem \ref{Theorem 1.13} only for the module $\mathcal{L}%
_{\mathcal{M}}(H,\mathcal{M}\otimes_{P}H)$ associated to a completely positive
map $P$ on a von Neumann algebra $\mathcal{M}$ and only for the special
covariant representation $(T,\sigma)$ which identifies $\mathcal{L}%
_{\mathcal{M}}(H,\mathcal{M}\otimes_{P}H)$ with $\mathcal{E}_{P}$. However, we
shall employ a picture of the dilation $(V,\rho)$ that is different from the
one constructed in Theorem \ref{Theorem 1.13}. It will play a critical role in
our analysis of semigroups of completely positive maps.

The definition of $(T,\sigma)$ is simple: $T$ maps $\mathcal{L}_{\mathcal{M}%
}(H,\mathcal{M}\otimes_{P}H)$ to $B(H)$ via the formula:
\begin{equation}
T(X):=W_{P}^{\ast}X\text{, \ \ \ }X\in\mathcal{L}_{\mathcal{M}}(H,\mathcal{M}%
\otimes_{P}H)\text{,} \label{covT}%
\end{equation}
and $\sigma$ is the identity representation,%
\[
\sigma(S)=S\text{,\ \ \ \ }S\in\mathcal{M}^{\prime}\text{.}%
\]
Of course, $\sigma$ is $\sigma$-weakly continuous. Also, a straightforward
calculation shows that $T$ is a bimodule map. To see that $T$ is completely
contractive, we appeal to \cite[Lemma 3.5]{MS98} and show that the linear
transformation $\tilde{T}:\mathcal{L}_{\mathcal{M}}(H,\mathcal{M}\otimes
_{P}H)\otimes_{\sigma}H\rightarrow H$ defined by the formula
\[
\tilde{T}(\sum X_{j}\otimes h_{j})=\sum W_{P}^{\ast}X_{j}h_{j}%
\]
is contractive. However, this is immediate:
\begin{align*}
\left\|  \sum W_{P}^{\ast}X_{j}h_{j}\right\|  ^{2}  &  \leq\left\|  \sum
X_{j}h_{j}\right\|  ^{2}\\
&  =\sum\langle h_{k},X_{k}^{\ast}X_{j}h_{j}\rangle=\left\|  \sum X_{j}\otimes
h_{j}\right\|  ^{2}\text{.}%
\end{align*}
As we remarked after Lemma \ref{CovRep}, $T$ is continuous with respect to the
$\sigma$-topology on $\mathcal{L}_{\mathcal{M}}(H,\mathcal{M}\otimes_{P}H)$
and the $\sigma$-weak topology on $B(H),$ and so $(T,\sigma)$ is a completely
contractive representation of $\mathcal{L}_{\mathcal{M}}(H,\mathcal{M}%
\otimes_{P}H)$ on $H$.

Evidently, $T$ is really the inverse of the map $X\rightarrow\Phi_{X}^{\ast}$
that we used to identify $\mathcal{L}_{\mathcal{M}}(H,\mathcal{M}\otimes
_{P}H)$ with $\mathcal{E}_{P}$ in Proposition \ref{Lemma1.4}. Indeed, using
the notation of that proposition, we see that for $Y\in\mathcal{E}_{P}$,
$T(\Phi_{Y}^{\ast})=W_{P}^{\ast}\Phi_{Y}^{\ast}=(\Phi_{Y}W_{P})^{\ast}=Y$. Now
all this may look trivial. It appears that after identifying $\mathcal{L}%
_{\mathcal{M}}(H,\mathcal{M}\otimes_{P}H)$ with $\mathcal{E}_{P}$ we are
simply studying the identity covariant representation of $\mathcal{E}_{P}$.
However, we need to emphasize that the heart of the matter lies in the fact
that the inner product on $\mathcal{E}_{P}$ is \emph{not} the one coming from
operator multiplication in $B(H)$ (unless $P$ is an endomorphism - see
Corollary \ref{corollary1.6}). Rather, it is defined through the map
$X\rightarrow\Phi_{X}^{\ast}$ (or through its inverse $T$) which identifies
$\mathcal{E}_{P}$ with $\mathcal{L}_{\mathcal{M}}(H,\mathcal{M}\otimes_{P}H)$.

\begin{definition}
\label{identityrep}The completely contractive covariant representation
$(T,\sigma)$ of $\mathcal{L}_{\mathcal{M}}(H,\mathcal{M}\otimes_{P}H)$, where
$T$ is defined by (\ref{covT}) and where $\sigma$ is the identity
representation, will be called the \emph{identity covariant representation }of
$\mathcal{L}_{\mathcal{M}}(H,\mathcal{M}\otimes_{P}H)$.
\end{definition}

As we noted above, and as we shall use to good effect, $(T,\sigma)$ really
identifies $\mathcal{L}_{\mathcal{M}}(H,\mathcal{M}\otimes_{P}H)$ with
$\mathcal{E}_{P}$ and when this identification is made, the maps $T$ and
$\sigma$ are both the identity maps.

To present the model for the minimal isometric dilation $(V,\rho)$ of
$(T,\sigma)$ with which we will work, we define, for $0\leq k<\infty$, maps
\[
\iota_{k}:\underset{k\;\text{times}}{\underbrace{\mathcal{M}\otimes
_{P}\mathcal{M}\otimes_{P}\cdots\otimes_{P}\mathcal{M}}}\otimes_{P}%
H\rightarrow\underset{k+1\;\text{times}}{\underbrace{\mathcal{M}\otimes
_{P}\mathcal{M}\otimes_{P}\cdots\otimes_{P}\mathcal{M}}}\otimes_{P}H
\]
by the formula $\iota_{k}(T_{1}\otimes T_{2}\otimes\cdots T_{k}\otimes
h)=I\otimes T_{1}\otimes T_{2}\otimes\cdots T_{k}\otimes h$. Of course,
$\iota_{0}=W_{P}$. Since $P$ is unital, this map is a well defined isometry of
$H_{k}:=\underset{k\;\text{times}}{\underbrace{\mathcal{M}\otimes
_{P}\mathcal{M}\otimes_{P}\cdots\otimes_{P}\mathcal{M}}}\otimes_{P}H$ into
$H_{k+1}:=\underset{k+1\;\text{times}}{\underbrace{\mathcal{M}\otimes
_{P}\mathcal{M}\otimes_{P}\cdots\otimes_{P}\mathcal{M}}}\otimes_{P}H$. We
write $H_{\infty}$ for the Hilbert space inductive limit, $\underrightarrow
{\lim}(H_{k},\iota_{k})$, and we write $W_{k}$ for the canonical (isometric)
embeddings of $H_{k}$ into $H_{\infty}$. Given $X\in\mathcal{L}_{\mathcal{M}%
}(H,\mathcal{M}\otimes_{P}H)$, we define $X_{k}:H_{k}\rightarrow H_{k+1}$ by
the formula $X_{k}(T_{1}\otimes T_{2}\otimes\cdots T_{k}\otimes h)=T_{1}%
\otimes T_{2}\otimes\cdots T_{k}\otimes Xh$. A straight forward calculation
using the fact that $X$ intertwines the actions of $\mathcal{M}$ on $H$ and on
$\mathcal{M}\otimes_{P}H$ shows that $X_{k}$ is bounded with $\left\|
X_{k}\right\|  \leq\left\|  X\right\|  $. Further, the diagram
\begin{equation}%
\begin{array}
[c]{ccccccccccccc}%
H & \overset{\iota_{0}}{\rightarrow} & H_{1} & \overset{\iota_{1}}%
{\rightarrow} & \cdots & \overset{\iota_{k-1}}{\rightarrow} & H_{k} &
\overset{\iota_{k}}{\rightarrow} & H_{k+1} & \rightarrow & \cdots &
\rightarrow & H_{\infty}\\
& \underset{X}{\searrow} &  & \underset{X}{\searrow}_{1} &  & \underset
{X}{\searrow}_{k-1} &  & \underset{X}{\searrow}_{k} &  & \underset{X}%
{\searrow}_{k+1} &  &  & \\
H & \overset{\iota_{0}}{\rightarrow} & H_{1} & \overset{\iota_{1}}%
{\rightarrow} & \cdots & \overset{\iota_{k-1}}{\rightarrow} & H_{k} &
\overset{\iota_{k}}{\rightarrow} & H_{k+1} & \rightarrow & \cdots &
\rightarrow & H_{\infty}%
\end{array}
\label{Vinfinity}%
\end{equation}
commutes and so defines an operator $X_{\infty}\in B(H_{\infty})$. We shall
see in a moment that the map $X\rightarrow X_{\infty}$, which we shall call
$V$, is part of an isometric covariant representation of $\mathcal{L}%
_{\mathcal{M}}(H,\mathcal{M}\otimes_{P}H)$, $(V,\rho)$.

To this end, we must first define $\rho$ through the following diagram, where
$S\in\mathcal{M}^{\prime}$. The diagram
\begin{equation}%
\begin{array}
[c]{ccccccccc}%
H & \overset{\iota_{0}}{\rightarrow} & H_{1} & \overset{\iota_{1}}%
{\rightarrow} & \cdots & \overset{\iota_{k-1}}{\rightarrow} & H_{k} &
\overset{\iota_{k}}{\rightarrow} & \cdots\\
\downarrow &
\begin{array}
[c]{rr}%
S & \;
\end{array}
& \downarrow &
\begin{array}
[c]{rr}%
I\otimes S & \;\;\;\;
\end{array}
&  &  & \downarrow &
\begin{array}
[c]{rr}%
\underset{}{I\otimes\cdots\otimes I}\otimes S & \;\;\;
\end{array}
& \\
H & \overset{\iota_{0}}{\rightarrow} & H_{1} & \overset{\iota_{1}}%
{\rightarrow} & \cdots & \overset{\iota_{k-1}}{\rightarrow} & H_{k} &
\overset{\iota_{k}}{\rightarrow} & \cdots
\end{array}
\label{rho}%
\end{equation}
commutes and, therefore, defines an operator $\rho(S)$ on $H_{\infty}$. Note
that
\[
W_{k}^{\ast}\rho(S)W_{k}=\underset{k\;\text{times}}{\underbrace{I\otimes
\cdots\otimes I}}\otimes S,
\]
where, recall, $W_{k}$ is the canonical embedding of $H_{k}$ in $H_{\infty}$.
>From this it is obvious that $\rho$ is a normal representation of
$\mathcal{M}^{\prime}$ on $H_{\infty}$ that is reduced by each of the spaces
$W_{k}H_{k}$. In particular, note that $W_{0}^{\ast}\rho(\cdot)W_{0}=\sigma$.

If the diagrams that define $V$ and $\rho$, (\ref{Vinfinity}) and (\ref{rho}),
resp., are combined in the obvious way, it becomes clear that for
$X\in\mathcal{L}_{\mathcal{M}}(H,\mathcal{M}\otimes_{P}H)$ and $S\in
\mathcal{M}^{\prime}$,
\[
V(XS)=(XS)_{\infty}=X_{\infty}\rho(S)=V(X)\rho(S)
\]
while
\[
V(S\cdot X)=V((I\otimes S)\circ X)=\rho(S)V(X)
\]
so that $(V,\rho)$ is covariant.

Next we show that $(V,\rho)$ is isometric. To this end, fix $X$ and $Y$ in
$\mathcal{L}_{\mathcal{M}}(H,\mathcal{M}\otimes_{P}H)$ and recall that
$X_{k}:H_{k}\rightarrow H_{k+1}$ is defined by the formula $X_{k}(T_{1}\otimes
T_{2}\otimes\cdots T_{k}\otimes h)=T_{1}\otimes T_{2}\otimes\cdots
T_{k}\otimes Xh$ and similarly for $Y_{k}$. Consequently, we find that
\[
X_{k}^{\ast}(T_{1}\otimes T_{2}\otimes\cdots T_{k+1}\otimes h)=T_{1}\otimes
T_{2}\otimes\cdots\otimes X^{\ast}(T_{k+1}\otimes h)
\]
because
\begin{align*}
&  \langle T_{1}\otimes T_{2}\otimes\cdots\otimes X^{\ast}(T_{k+1}\otimes
h),S_{1}\otimes S_{2}\otimes\cdots\otimes S_{k}\otimes k\rangle\\
&  =\langle X^{\ast}(T_{k+1}\otimes h),P(T_{k}^{\ast}P(\cdots)S_{k})k\rangle\\
&  =\langle T_{k+1}\otimes h,XP(T_{k}^{\ast}P(\cdots)S_{k})k\rangle\\
&  =\langle T_{k+1}\otimes h,(P(T_{k}^{\ast}P(\cdots)S_{k})\otimes
I)Xk\rangle\;\;\text{(because }X\in\mathcal{L}_{\mathcal{M}}(H,\mathcal{M}%
\otimes_{P}H)\text{)}\\
&  =\langle T_{1}\otimes T_{2}\otimes\cdots\otimes T_{k+1}\otimes
h,S_{1}\otimes S_{2}\otimes\cdots\otimes S_{k}\otimes Xk\rangle\\
&  =\langle T_{1}\otimes T_{2}\otimes\cdots\otimes T_{k+1}\otimes
h,X_{k}(S_{1}\otimes S_{2}\otimes\cdots\otimes S_{k}\otimes k)\rangle\text{.}%
\end{align*}
Therefore,
\begin{align*}
X_{k}^{\ast}Y_{k}(T_{1}\otimes T_{2}\otimes\cdots\otimes T_{k}\otimes h)  &
=X_{k}^{\ast}(T_{1}\otimes T_{2}\otimes\cdots\otimes T_{k}\otimes Yh)\\
&  =T_{1}\otimes T_{2}\otimes\cdots\otimes T_{k}\otimes X^{\ast}Yh\\
&  =W_{k}^{\ast}\rho(X^{\ast}Y)W_{k}(T_{1}\otimes T_{2}\otimes\cdots\otimes
T_{k+1}\otimes h).
\end{align*}
Thus $W_{k}^{\ast}\rho(X^{\ast}Y)W_{k}=X_{k}^{\ast}Y_{k}=W_{k}^{\ast
}V(X)^{\ast}W_{k+1}W_{k+1}^{\ast}V(Y)W_{k}$ for all $k$ , from which it
follows that $V(X)^{\ast}V(Y)=\rho(\langle X,Y\rangle)$, i.e., that $(V,\rho)$
is isometric.

We now show that $(V,\rho)$ dilates $(T,\sigma)$ in the sense described in
Theorem \ref{Theorem 1.13}. Of course to do this, we must, strictly speaking,
identify $H$ with the subspace $W_{0}H$ of $H_{\infty}$. When this is done,
the projection $P$ of $H_{\infty}$ on $H$ is $W_{0}W_{0}^{\ast}$. We already
have seen that $H=W_{0}H$ reduces $\rho$ and that $\rho|H=\sigma$ as is
required in part 1. of Theorem \ref{Theorem 1.13}. Also note that for
$X\in\mathcal{L}_{\mathcal{M}}(H,\mathcal{M}\otimes_{P}H)$,
\[
W_{0}^{\ast}X_{\infty}W_{0}=\iota_{0}^{\ast}X=T(X)\text{,}%
\]
which is an evident consequence of the properties of inductive limits: For
$h\in H$, $X_{\infty}W_{0}h=W_{1}Xh$, and $W_{0}^{\ast}W_{1}=\iota_{0}^{\ast}%
$. This of course means that $T(X)=W_{0}^{\ast}V(X)W_{0}$, so that after
identifying $H$ with $W_{0}H$, through $W_{0}$, we see that $T(X)=PV(X)|H$,
$X\in\mathcal{L}_{\mathcal{M}}(H,\mathcal{M}\otimes_{P}H)$, as required in
part 2. of Theorem \ref{Theorem 1.13}. But we also need to check that
$V(X)^{\ast}$ maps $H$ into itself. Equivalently, we need to show that $V(X)$
maps $H_{\infty}\ominus H$ into itself.

For this purpose, it suffices to show that for each $k\geq1$, $V(X)$ maps
$W_{k}H_{k}\ominus H$ into $H_{\infty}\ominus H$. To show how this is done,
but to keep the matters simple, we show that $V(X)$ maps $W_{2}H_{2}\ominus H$
into $H_{\infty}\ominus H$. So let $W_{2}\sum T_{i}\otimes S_{i}\otimes h_{i}$
be an element of $W_{2}H_{2}.$ To say this is orthogonal to $H$ means that for
all $h\in H$,%
\begin{align*}
0  &  =\langle W_{2}\sum T_{i}\otimes S_{i}\otimes h_{i},W_{0}h\rangle\\
&  =\langle\sum T_{i}\otimes S_{i}\otimes h_{i},I\otimes I\otimes h\rangle\\
&  =\langle\sum P(P(T_{i})S_{i})h_{i},h\rangle\text{;}%
\end{align*}
i.e., $W_{2}\sum T_{i}\otimes S_{i}\otimes h_{i}\in W_{2}H_{2}\ominus H$ if
and only if $\sum P(P(T_{i})S_{i})h_{i}=0$. Now assume that $W_{2}\sum
T_{i}\otimes S_{i}\otimes h_{i}\in W_{2}H_{2}\ominus H$, let $h$ be an element
in $H$, and compute:%
\begin{align*}
&  \langle V(X)W_{2}\sum T_{i}\otimes S_{i}\otimes h_{i},W_{0}h\rangle\\
&  =\langle\sum T_{i}\otimes S_{i}\otimes Xh_{i},I\otimes I\otimes I\otimes
h\rangle_{H_{3}}\\
&  =\sum\langle Xh_{i},P(S_{i}^{\ast}P(T_{i}^{\ast}))\otimes h\rangle_{H_{1}%
}\\
&  =\sum\langle(P(P(T_{i})S_{i})\otimes I)Xh_{i},I\otimes h\rangle_{H_{1}}\\
&  =\sum\langle X(P(P(T_{i})S_{i}))h_{i},I\otimes h\rangle_{H_{1}}\\
&  =\langle X(\sum P(P(T_{i})S_{i}))h_{i}),I\otimes h\rangle_{H_{1}}\\
&  =0\text{,}%
\end{align*}
since $\sum P(P(T_{i})S_{i}))h_{i}=0$. Thus $(V,\rho)$ satisfies condition 2.
in Theorem \ref{Theorem 1.13}.

To show that this $(V,\rho)$ is unitarily equivalent to the dilation of
$(T,\sigma)$ that is provided by Theorem \ref{Theorem 1.13}, we appeal to
Proposition 3.2 of \cite{MS98} (which we stated as part of Theorem
\ref{Theorem 1.13}) and show that $(V,\rho)$ is minimal; i.e., that there are
no closed subspaces $K$ properly contained between $H$ and $H_{\infty}$ that
are invariant under the images of $V$ and $\rho$. So, suppose $K$ is such a
subspace, then for every $X\in\mathcal{L}_{\mathcal{M}}(H,\mathcal{M}%
\otimes_{P}H)$, $X_{\infty}(H)=V(X)(H)$ is contained in $K$. Hence, in
particular, $X_{\infty}(H)=X_{\infty}(W_{0}H)=X(H)\subseteq K$. However, the
span of $\{X(H)\mid X\in\mathcal{L}_{\mathcal{M}}(H,\mathcal{M}\otimes
_{P}H)\}$ is $M\otimes_{P}H$, by Lemma \ref{lemma1.7}, and so we conclude that
$W_{1}H_{1}=W_{1}(\mathcal{M}\otimes_{P}H)\subseteq K.$ This, in turn, implies
that $X_{\infty}(\mathcal{M}\otimes_{P}H)\subseteq K$; i.e., that $(I\otimes
X)(\mathcal{M}\otimes_{P}H)\subseteq K$. Applying Lemma \ref{lemma1.7} again,
we see that $W_{2}H_{2}=W_{2}(\mathcal{M}\otimes_{P}\mathcal{M}\otimes_{P}H)$
is contained in $K$. Continuing in this manner, we find that $W_{k}H_{k}$ is
contained in $K$ for every $k$. Hence $K=H_{\infty}$.

Since our special $(V,\rho)$ is unitarily equivalent to the one provided by
Theorem \ref{Theorem 1.13}, we may infer that $V$ is continuous with respect
to the $\sigma$-topology on $\mathcal{L}_{\mathcal{M}}(H,\mathcal{M}%
\otimes_{P}H)$ and the $\sigma$-weak topology on $B(H)$.

We summarize our discussion of the identity representation of $\mathcal{L}%
_{\mathcal{M}}(H,\mathcal{M}\otimes_{P}H)$ in the following theorem.

\begin{theorem}
\label{idDilation}The maps $V$ and $\rho$ defined by the diagrams,
(\ref{Vinfinity}) and (\ref{rho}) together form an isometric covariant
representation of $\mathcal{L}_{\mathcal{M}}(H,\mathcal{M}\otimes_{P}H)$ that
dilates the identity representation $(T,\sigma)$ of $\mathcal{L}_{\mathcal{M}%
}(H,\mathcal{M}\otimes_{P}H)$. Moreover, $(T,\sigma)$ and $(V,\rho)$ are fully coisometric.
\end{theorem}

\begin{proof}
The only thing that remains to be proved is the last statement about
$(T,\sigma)$ and $(V,\rho)$ being fully coisometric. However, for this
purpose, it suffices to show that $(T,\sigma)$ is fully coisometric, by
Theorem \ref{Theorem 1.13}. Recall that $\tilde{T}$ maps $\mathcal{L}%
_{\mathcal{M}}(H,\mathcal{M}\otimes_{P}H)\otimes_{\sigma}H$ to $H$ by the
formula $\tilde{T}(X\otimes h)=W_{P}^{\ast}Xh$. To calculate $\tilde{T}^{\ast
}$, simply observe that for $X\in\mathcal{L}_{\mathcal{M}}(H,\mathcal{M}%
\otimes_{P}H)$ and $h\in H$, $\langle\tilde{T}^{\ast}k,X\otimes h\rangle
=\langle k,\tilde{T}(X\otimes h)\rangle=\langle k,W_{P}^{\ast}Xh\rangle
=\langle W_{P}k,Xh\rangle.$ However, by Lemma \ref{lemma1.7}, $\{Xh\mid
X\in\mathcal{L}_{\mathcal{M}}(H,\mathcal{M}\otimes_{P}H),\;h\in H\}$ spans
$\mathcal{M}\otimes_{P}H$. So, if we let $u:\mathcal{L}_{\mathcal{M}%
}(H,\mathcal{M}\otimes_{P}H)\otimes H\rightarrow\mathcal{M}\otimes_{P}H$ be
defined by the formula $u(X\otimes h)=Xh$, the $u$ is a Hilbert space
isomorphism such that $\langle\widetilde{T}^{\ast}k,X\otimes h\rangle=\langle
W_{P}k,Xh\rangle=\langle W_{P}k,u(X\otimes h)\rangle=\langle u^{\ast}%
W_{P}k,X\otimes h\rangle$ for all $k$ and all $X\otimes h$. Thus $\tilde
{T}^{\ast}$ is the \emph{isometry }$u^{\ast}W_{P}$, proving that $\tilde{T}$
is a coisometry and, therefore, that $(T,\sigma)$ is fully coisometric.
\end{proof}

If $\mathcal{N}$ is a von Neumann algebra and if $\mathcal{E}$ is a $W^{\ast}%
$-correspondence over $\mathcal{N}$, then we have seen how a completely
contractive covariant representation $(T,\sigma)$ of $\mathcal{E}$ on a
Hilbert space $H$ gives rise to a completely positive map $\Psi=\Psi_{T}$ of
$\mathcal{L}(\mathcal{E})$ on $H$. (See Definition \ref{Definition1.12bis}.)
However, equally important for our purposes is the related completely positive
map $\Theta=\Theta_{T}$ on the \emph{commutant }of $\sigma(\mathcal{N})$,
$\sigma(\mathcal{N})^{\prime}$, that is described in the next proposition.

\begin{proposition}
\label{Proposition 1.15}Let $\mathcal{N}$ be a von Neumann algebra, let
$\mathcal{E}$ be a $W^{\ast}$-correspondence over $\mathcal{N}$, and let
$(T,\sigma)$ be a completely contractive covariant representation of
$\mathcal{E}$ on a Hilbert space $H$. For $S\in\sigma(\mathcal{N})^{\prime}$,
set%
\begin{equation}
\Theta(S)=\Theta_{T}(S):=\tilde{T}(1_{\mathcal{E}}\otimes S)\tilde{T}^{\ast
}\text{.} \label{thetaT}%
\end{equation}
Then $\Theta$ is normal completely positive map from $\sigma(\mathcal{N}%
)^{\prime}$ into itself that is unital if and only if $(T,\sigma)$ is fully
coisometric. Further, if $(T,\sigma)$ is isometric, then $\Theta$ is
multiplicative, i.e., $\Theta$ is an endomorphism of $\sigma(\mathcal{N}%
)^{\prime}$, and, conversely, if $\Theta$ is multiplicative, then the
correspondence $\mathcal{E}$ decomposes as the direct sum of two
subcorrespondences, $\mathcal{E}=\mathcal{E}_{1}\oplus\mathcal{E}_{2}$, so
that $(T|\mathcal{E}_{1},\sigma)$ is isometric and $T|\mathcal{E}_{2}=0$.
\end{proposition}

\begin{proof}
Much of the proof may be dug out of \cite{MS99}. See Lemma 2.3 there, in
particular. Here are the particulars. First, recall the induced representation
$\sigma^{\mathcal{E}}:\mathcal{L}(\mathcal{E})\rightarrow B(\mathcal{E}%
\otimes_{\sigma}H)$, $\sigma^{\mathcal{E}}(X)=X\otimes I_{H}$. As Rieffel
shows in Theorem 6.23 of \cite{mR74}, the commutant of $\sigma^{\mathcal{E}%
}(\mathcal{L}(\mathcal{E}))$ is $\mathbb{C}1_{\mathcal{E}}\otimes
\sigma(\mathcal{N})^{\prime}$, and of course the map $S\rightarrow
1_{\mathcal{E}}\otimes S$ is a normal representation of $\sigma(\mathcal{N}%
)^{\prime}$ onto $\mathbb{C}1_{\mathcal{E}}\otimes\sigma(\mathcal{N})^{\prime
}$. Thus $\Theta$ is a normal completely positive map from $\sigma
(\mathcal{N})^{\prime}$ into $B(H)$. The problem is to locate its range. This,
however, is easy on the basis of equation (\ref{covariance}): Given
$R\in\mathcal{N}$ and $S\in\sigma(\mathcal{N})^{\prime}$, that equation
implies that%
\begin{align*}
\sigma(R)\Theta(S)=\sigma(R)\tilde{T}(1_{\mathcal{E}}\otimes S)\tilde{T}%
^{\ast}=\tilde{T}\sigma^{\mathcal{E}}\circ\varphi(R)(1_{\mathcal{E}}\otimes
S)\tilde{T}^{\ast}\\
=\tilde{T}(1_{\mathcal{E}}\otimes S)\sigma^{\mathcal{E}}\circ\varphi
(R)\tilde{T}^{\ast}=\tilde{T}(1_{\mathcal{E}}\otimes S)\tilde{T}^{\ast}%
\sigma(R)\\
=\Theta(S)\sigma(R),
\end{align*}
so $\Theta(S)\in\sigma(\mathcal{N})^{\prime}$.

Of course $\Theta$ is unital if and only if $(T,\sigma)$ is fully coisometric.

As for the last assertion, the direct statement is proved as Lemma 2.3 of
\cite{MS99}. For the converse, suppose that $\Theta$ is multiplicative. Then
$\tilde{T}\tilde{T}^{\ast}=\Theta(I)$ is a projection. Therefore, $\tilde
{T}^{\ast}\tilde{T}$ is a projection on $\mathcal{E}\otimes_{\sigma}H$, call
it $q$. Since $\Theta$ is multiplicative, we infer that $q(1_{\mathcal{E}%
}\otimes S_{1})q(1_{\mathcal{E}}\otimes S_{2})q=q(1_{\mathcal{E}}\otimes
S_{1}S_{2})q$ for all $S_{1},S_{2}\in\sigma(\mathcal{N})^{\prime}$. This
implies that $q\in(\mathbb{C}1_{\mathcal{E}}\otimes\sigma(\mathcal{N}%
)^{\prime})^{\prime}=\sigma^{\mathcal{E}}(\mathcal{L}(\mathcal{E}))$, by
Rieffel's theorem \cite[Theorem 6.23]{mR74} and the fact that $\sigma
^{\mathcal{E}}$ is a normal representation of the von Neumann algebra
$\mathcal{L}(\mathcal{E})$. Thus, $q=\sigma^{\mathcal{E}}(Q)$ for a projection
$Q\in\mathcal{L}(\mathcal{E})$. If $\mathcal{E}_{1}:=Q\mathcal{E}$ and
$\mathcal{E}_{2}:=(1_{\mathcal{E}}-Q)\mathcal{E}$, then it is easy to see that
$(T|\mathcal{E}_{1},\sigma)$ is isometric, while $T|\mathcal{E}_{2}=0$. We
omit the details.
\end{proof}

\begin{definition}
\label{inducedCPmap}Let $\mathcal{E}$ be a $W^{\ast}$-correspondence over a
von Neumann algebra $\mathcal{N}$ and let $(T,\sigma)$ be a completely
contractive covariant representation of $\mathcal{E}$ on the Hilbert space
$H$, the normal, completely positive map
\[
\Theta_{T}:\sigma(\mathcal{N})^{\prime}\rightarrow\sigma(\mathcal{N})^{\prime}%
\]
defined by equation (\ref{thetaT}) will be called the \emph{induced
(completely positive) map }on $\sigma(\mathcal{N})^{\prime}$. If $T$ is
isometric, then $\Theta_{T}$ will be called the \emph{induced endomorphism} of
$\sigma(\mathcal{N})^{\prime}$.
\end{definition}

If we apply Proposition \ref{Proposition 1.15} to the identity representation
$(T,\sigma)$ of $\mathcal{L}_{\mathcal{M}}(H,\mathcal{M}\otimes_{P}H)$ or of
$\mathcal{E}_{P}$, for a completely positive map $P$ on a von Neumann algebra
$\mathcal{M}$, we recapture $P$. Specifically, we have

\begin{corollary}
\label{lemma 1.16}Let $P$ be a normal, unital, completely positive map on the
von Neumann algebra $\mathcal{M}$, and let $(T,\sigma)$ be the identity
representation of the Arveson correspondence $\mathcal{E}_{P}\simeq
\mathcal{L}_{\mathcal{M}}(H,\mathcal{M}\otimes_{P}H)$ on $H$. Then $\Theta
_{T}=P$.
\end{corollary}

\begin{proof}
We will apply Proposition \ref{Proposition 1.15}, with the von Neumann algebra
$\mathcal{N}$ identified with $\mathcal{M}^{\prime}$. So, for $S\in
\mathcal{M}$ and $h\in H$, we have from the computations in the proof of
Proposition \ref{Proposition 1.15} (the fact that $\tilde{T}^{\ast}=u^{\ast
}W_{P}$, so $\tilde{T}=W_{P}^{\ast}u$
\begin{align*}
\Theta(S)h &  =\tilde{T}(1_{\mathcal{E}}\otimes S)\tilde{T}^{\ast}h\\
&  =\tilde{T}(1_{\mathcal{E}}\otimes S)W_{P}h\\
&  =W_{P}^{\ast}u(I_{\mathcal{E}}\otimes S)u^{\ast}W_{P}h=P(S)h\text{.}%
\end{align*}
\end{proof}

We conclude this section with our principal dilation result for single
completely positive maps. It is the key to our analysis of semigroups.

\begin{theorem}
\label{Theorem1.17}Let $\mathcal{M}$ be a von Neumann algebra acting on a
Hilbert space $H$ and let $P:\mathcal{M}\rightarrow\mathcal{M}$ be a normal,
unital, completely positive map of $\mathcal{M}$. Let $(T,\sigma)$ be the
identity representation on $H$ of the Arveson correspondence $\mathcal{E}_{P}%
$, let $(V,\rho)$ be the minimal isometric dilation of $(T,\sigma)$ on the
Hilbert space $K$, and let $W:H\rightarrow K$ be the associated imbedding. If
$\mathcal{R}:=\rho(\mathcal{M}^{\prime})^{\prime}$, then

\begin{enumerate}
\item $W^{\ast}\mathcal{R}W=\mathcal{M}$, so that $\mathcal{M}$ is a
\emph{corner} of $\mathcal{R}$ (and $\mathcal{R}^{\prime}$ is a normal
homomorphic image of $\mathcal{M}^{\prime}$),

\item $\Theta_{V}$ is a unital, normal $\ast$-endomorphism of $\mathcal{R}$, and

\item for every non-negative integer $n$,
\[
P^{n}(T)=W^{\ast}\Theta_{V}^{n}(WTW^{\ast})W
\]
and%
\[
P^{n}(W^{\ast}SW)=W^{\ast}\Theta_{V}^{n}(S)W
\]
for all $S\in\mathcal{R}$, and $T\in\mathcal{M}$.
\end{enumerate}
\end{theorem}

Thus, the induced endomorphism $\Theta_{V}$ of $\mathcal{R}$ is a \emph{power
dilation} of $P$.

\begin{proof}
>From Corollary \ref{lemma 1.16}, we know that $P$ is the induced completely
positive map $\Theta_{T}$ on $\mathcal{M}$. Also, since $(V,\rho)$ is the
minimal isometric dilation of $(T,\sigma)$ and $W$ is the embedding map, we
know that $WH$ is invariant under $V(Y)^{\ast}$ for all $Y\in\mathcal{E}$ and
$W^{\ast}V(Y)W=T(Y)$. Since $W^{\ast}\rho(S)W=\sigma(S)$ for all
$S\in\mathcal{M}^{\prime}$ by definition of $(V,\rho)$, and since
$\sigma(S)=S$, $S\in\mathcal{M}^{^{\prime}}$, by definition of the identity
representation, we see that%
\begin{equation}
W^{\ast}\mathcal{R}W=W^{\ast}\rho(\mathcal{M}^{\prime})^{\prime}W=(W^{\ast
}\rho(\mathcal{M}^{\prime})W)^{\prime}=(\mathcal{M}^{\prime})^{\prime
}=\mathcal{M}\text{.}\label{compress}%
\end{equation}
By Theorem \ref{idDilation}, $(T,\sigma)$ and $(V,\rho)$ are fully
coisometric, and so, by Proposition \ref{Proposition 1.15}, $\Theta_{V}$ is a
normal, unital, $\ast$-endomorphism of $\mathcal{R}=\rho(\mathcal{M}^{\prime
})^{\prime}$. Since $WH$ is invariant under $V(Y)^{\ast}$, $Y\in\mathcal{E}$,
we see that for $Y\in\mathcal{E}$ and $k\in K$,
\[
WW^{\ast}\tilde{V}(Y\otimes(I-WW^{\ast})k)=WW^{\ast}V(Y)(I-WW^{\ast})k=0
\]
so that $WW^{\ast}\tilde{V}(I\otimes(I-WW^{\ast}))=0$. Therefore, $WW^{\ast
}\Theta_{V}((I-WW^{\ast}))=WW^{\ast}\tilde{V}(I\otimes(I-WW^{\ast}))\tilde
{V}^{\ast}=0$; i.e. $WW^{\ast}\Theta_{V}(WW^{\ast})=WW^{\ast}$. Multiplying
this equation on the left by $W^{\ast}$, we see that
\begin{equation}
W^{\ast}\Theta_{V}(WW^{\ast})=W^{\ast}.\label{key1}%
\end{equation}

Since $T(\cdot)=W^{\ast}V(\cdot)W$, it follows that $\tilde{T}=W^{\ast}%
\tilde{V}(I\otimes W)$. Consequently, for $L\in\mathcal{M}$,
\begin{align}
P(L)  &  =\tilde{T}(I\otimes L)\tilde{T}^{\ast}\label{key2}\\
&  =W^{\ast}\tilde{V}(I\otimes W)(I\otimes L)(I\otimes W^{\ast})\tilde
{V}^{\ast}W=W^{\ast}\Theta_{V}(WLW^{\ast})W\text{.}%
\end{align}
On the other hand, for $S\in\mathcal{R}$, we find from this equation and the
fact that $W^{\ast}SW\in\mathcal{M}$ (by (\ref{compress})) that
\begin{align*}
P(W^{\ast}SW)  &  =W^{\ast}\Theta_{V}(WW^{\ast}SWW^{\ast})W\\
&  =W^{\ast}\Theta_{V}(WW^{\ast}SWW^{\ast})W=W^{\ast}\Theta_{V}(WW^{\ast
})\Theta_{V}(S)\Theta_{V}(WW^{\ast})W\\
&  =W^{\ast}\Theta_{V}(S)W\text{,}%
\end{align*}
using equation (\ref{key1}).

To relate $P^{2}$ to $\Theta_{V}^{2}$, let $L\in\mathcal{M}$. Then, using
equation (\ref{key1}) again, we find that%
\begin{align*}
P^{2}(L)  &  =P(P(L))=P(W^{\ast}\Theta_{V}(WLW^{\ast})W)\\
&  =W^{\ast}\Theta_{V}(WW^{\ast}\Theta_{V}(WLW^{\ast})WW^{\ast})W\\
&  =W^{\ast}\Theta_{V}(WW^{\ast})\Theta_{V}^{2}(WLW^{\ast})\Theta_{V}%
(WW^{\ast})W\\
&  =W^{\ast}\Theta_{V}^{2}(WLW^{\ast})W\text{.}%
\end{align*}

Continuing in this manner, we find that $P^{n}(L)=W^{\ast}\Theta_{V}%
^{n}(WLW^{\ast})W$ for all $L\in\mathcal{M}$.

To show that $P^{n}(W^{\ast}SW)=W^{\ast}\Theta_{V}^{n}(S)W$ for all
$S\in\mathcal{R}$, and all $n$, we need to generalize equation (\ref{key1}) to
$W^{\ast}\Theta_{V}^{n}(WW^{\ast})=W^{\ast}$, for all $n$. However, this is an
easy induction, the general step of which is:
\begin{align*}
WW^{\ast}\Theta_{V}^{n+1}(WW^{\ast})  &  =WW^{\ast}\Theta_{V}(\Theta_{V}%
^{n}(WW^{\ast}))\\
&  =WW^{\ast}\Theta_{V}(WW^{\ast}\Theta_{V}^{n}(WW^{\ast}))\\
&  +WW^{\ast}\Theta_{V}((I-WW^{\ast})\Theta_{V}^{n}(WW^{\ast}))\\
&  =WW^{\ast}\Theta_{V}^{n}(WW^{\ast})+WW^{\ast}\Theta_{V}(I-WW^{\ast}%
)\Theta_{V}^{n}(WW^{\ast})\\
&  =WW^{\ast}\text{.}%
\end{align*}
Thus $WW^{\ast}\Theta_{V}^{n}(WW^{\ast})=WW^{\ast}$ for all $n$. Multiplying
through on the left by $W^{\ast}$ gives the desired formula.

Using this, we see that since $W^{\ast}SW\in\mathcal{M}$ for all
$S\in\mathcal{R}$, our earlier calculation gives
\begin{align*}
P^{n}(W^{\ast}SW)  &  =W^{\ast}\Theta_{V}^{n}(WW^{\ast}SWW^{\ast})W\\
&  =W^{\ast}\Theta_{V}^{n}(WW^{\ast})\Theta_{V}^{n}(S)\Theta_{V}^{n}(WW^{\ast
})W\\
&  =W^{\ast}\Theta_{V}^{n}(S)W.
\end{align*}
\end{proof}

\section{Semigroups of Completely Positive Maps}

In this section, we focus on semigroups $\{P_{t}\}_{t\geq0}$ of unital,
normal, completely positive maps on our basic von Neumann algebra
$\mathcal{M}$ acting on a Hilbert space $H$. That is, we assume that
$P_{t+s}=P_{t}P_{s}$, $s,t\geq0$, and $P_{0}$ is the identity map on
$\mathcal{M}$. We call $\{P_{t}\}_{t\geq0}$ a \emph{completely positive
semigroup} on $\mathcal{M}$, or simply a \emph{cp semigroup}, for short. We
make no continuity assumptions on $\{P_{t}\}_{t\geq0}$ in this section and, in
fact, everything we say is true if the additive semigroup of non-negative real
numbers is replaced by any totally ordered semigroup. Our goal is to dilate
$\{P_{t}\}_{t\geq0}$ to a semigroup of endomorphisms in much the same fashion
that we did for a single completely positive map in Section \ref{section1}.
However, there is a complication that must be addressed.

Let $\mathcal{E}_{t}$ be the Arveson correspondence over $\mathcal{M}^{\prime
}$ associated with $P_{t}$, $t\geq0$. As in Section \ref{section1}, we shall
view $\mathcal{E}_{t}$ as either a space of operators on $H$ or as the space
$\mathcal{L}_{\mathcal{M}}(H,\mathcal{M}\otimes_{P_{t}}H)$. As we noted, the
spaces $\mathcal{E}_{t}$ need not ``multiply'', i.e., $\mathcal{E}_{t}%
\otimes\mathcal{E}_{s}$ need not be isomorphic to $\mathcal{E}_{t+s}$. So, we
will have to ``dilate'' these to a family $\{E(t)\}_{t\geq0}$ of
$\mathcal{M}^{\prime}$ correspondences such that $E(t)\otimes E(s)\simeq
E(t+s)$. That is, we need to dilate these to a (discrete) product system over
$\mathcal{M}^{\prime}$ - a notion that is inspired by Arveson's product
systems in \cite{wA89}. This we do following in outline arguments of Bhat in
\cite{bB96}. There are similarities also between our arguments and arguments
in \cite{BS00}, but our correspondences are over $\mathcal{M}^{\prime}$ as
opposed to being over $\mathcal{M}$ and we cannot tap directly into their
arguments. Once $\{E(t)\}_{t\geq0}$ is constructed, we promote the identity
representations of the $\mathcal{E}_{t}$'s to completely contractive
representations of the $E(t)$'s and then dilate these to isometric
representations of the $E(t)$'s. These last representations will implement a
semigroup of endomorphisms of a bigger von Neumann algebra in which
$\mathcal{M}$ sits as a corner. The semigroup of endomorphisms will be the
desired dilation of $\{P_{t}\}_{t\geq0}$.

Let $\mathfrak{P}(t)$ denote the collection of partitions of the closed
interval $[0,t]$ and order these by refinement. For a $\mathfrak{p}%
\in\mathfrak{P}(t)$, we shall write $\mathfrak{p}=\{0=t_{0}<t_{1}<t_{2}%
<\cdots<t_{n-1}<t_{n}=t\}$. For such a $\mathfrak{p}$, we shall write
\[
H_{\mathfrak{p},t}:=\mathcal{M}\otimes_{P_{t_{1}}}\mathcal{M}\otimes
_{P_{t_{2}-t_{1}}}\otimes\cdots\mathcal{M}\otimes_{P_{t-t_{n-1}}}H\text{.}%
\]
Then it is easy to see that $H_{\mathfrak{p},t}$ is a left $\mathcal{M}%
$-module via the formula $S\cdot(T_{1}\otimes T_{2}\otimes\cdots T_{n}\otimes
h):=(ST_{1})\otimes T_{2}\otimes\cdots T_{n}\otimes h$, $S\in\mathcal{M}$,
$(T_{1}\otimes T_{2}\otimes\cdots T_{n}\otimes h)\in H_{\mathfrak{p},t}$.
Also, it is easy to see that $\mathcal{L}_{\mathcal{M}}(H,H_{\mathfrak{p},t})$
becomes an $\mathcal{M}^{\prime}$-correspondence via the actions%
\[
(XR)h:=X(Rh)
\]
and%
\[
(RX)h:=(I\otimes R)Xh\text{,}%
\]
$R\in\mathcal{M}^{\prime}$, $X\in\mathcal{L}_{\mathcal{M}}(H,H_{\mathfrak{p}%
,t})$ and $h\in H$, where $I$ is the identity operator on $H_{\mathfrak{p},t}%
$. The inner product is given by the formula $\langle X_{1},X_{2}%
\rangle:=X_{1}^{\ast}X_{2}$. Note that the map $R\mapsto\langle X_{1}%
,RX_{2}\rangle=X_{1}^{\ast}(I\otimes R)X_{2}$ is $\sigma$-weakly continuous,
so that $\mathcal{L}_{\mathcal{M}}(H,H_{\mathfrak{p},t})$ is, indeed, an
$\mathcal{M}^{\prime}$-correspondence.

We shall write $\mathcal{L}_{t}$ for the $\mathcal{M}^{\prime}$-correspondence
$\mathcal{L}_{\mathcal{M}}(H,\mathcal{M}\otimes_{P_{t}}H)$. Then Proposition
\ref{Lemma1.9} shows that $\mathcal{L}_{\mathcal{M}}(H,H_{\mathfrak{p},t})$ is
isomorphic to $\mathcal{L}_{t-t_{n-1}}\otimes_{\mathcal{M}^{\prime}%
}\mathcal{L}_{t_{n-1}-t_{n-2}}\otimes\cdots\otimes_{\mathcal{M}^{\prime}%
}\mathcal{L}_{t_{1}}$ as $\mathcal{M}^{\prime}$-correspondences.

We next want to show that the Hilbert spaces $H_{\mathfrak{p},t}$ and
$\mathcal{M}^{\prime}$-correspondences $\mathcal{L}_{\mathcal{M}%
}(H,H_{\mathfrak{p},t})$ form inductive systems so that we can take their
direct limits. For this purpose, consider first the case when $\mathfrak{p}%
^{\prime}:=\{0=t_{0}<t_{1}<t_{2}<\cdots<t_{k}<\tau<t_{k+1}<\cdots
t_{n-1}<t_{n}=t\}$, a one point refinement of $\mathfrak{p}$. Then we obtain a
Hilbert space isometry $v_{0}:H_{\mathfrak{p},t}\longrightarrow
H_{\mathfrak{p}^{\prime},t}$ defined by the formula $v_{0}(T_{1}\otimes
T_{2}\otimes\cdots T_{n}\otimes h)=T_{1}\otimes\cdots\otimes T_{k}\otimes
I\otimes T_{k+1}\otimes\cdots T_{n}\otimes h$ and an $\mathcal{M}^{\prime}%
$-correspondence isometry $v:\mathcal{L}_{\mathcal{M}}(H,H_{\mathfrak{p}%
,t})\longrightarrow\mathcal{L}_{\mathcal{M}}(H,H_{\mathfrak{p}^{\prime},t})$
defined by the formula $v(X):=v_{0}\circ X.$ The proof of these facts is a
minor modification of the proof of Proposition \ref{Lemma1.9} and so will be omitted.

Since every refinement of a partition can be obtained by a sequence of
one-point refinements, it is clear that for every pair of partitions
$(\mathfrak{p},\mathfrak{p}^{\prime})$, with $\mathfrak{p}^{\prime}$ refining
$\mathfrak{p}$, we have Hilbert space isometries $v_{0,\mathfrak{p,p}^{\prime
}}:H_{\mathfrak{p},t}\rightarrow H_{\mathfrak{p}^{\prime},t}$ and
$\mathcal{M}^{\prime}$-correspondence isometries $v_{\mathfrak{p,p}^{\prime}%
}:\mathcal{L}_{\mathcal{M}}(H,H_{\mathfrak{p},t})\longrightarrow
\mathcal{L}_{\mathcal{M}}(H,H_{\mathfrak{p}^{\prime},t})$ so that
$v_{0,\mathfrak{p}^{\prime}\mathfrak{,p}^{\prime\prime}}\circ
v_{0,\mathfrak{p,p}^{\prime}}=v_{0,\mathfrak{p,p}^{\prime\prime}}$ and
$v_{\mathfrak{p}^{\prime}\mathfrak{,p}^{\prime\prime}}\circ v_{\mathfrak{p,p}%
^{\prime}}=v_{\mathfrak{p,p}^{\prime\prime}}$, when $\mathfrak{p}%
^{\prime\prime}$ refines $\mathfrak{p}^{\prime}$ and $\mathfrak{p}^{\prime}$
refines $\mathfrak{p}$. The Hilbert space isometry $v_{0,\mathfrak{p,p}%
^{\prime}}$ simply sends a decomposable tensor $T_{1}\otimes T_{2}%
\otimes\cdots T_{n}\otimes h\in H_{\mathfrak{p},t}$ to the decomposable tensor
in $H_{\mathfrak{p}^{\prime},t}$ obtained from $T_{1}\otimes T_{2}%
\otimes\cdots T_{n}\otimes h$ by inserting identity operators in those
positions where new indices have been added to $\mathfrak{p}$ to obtain
$\mathfrak{p}^{\prime}$. The $\mathcal{M}^{\prime}$-correspondence isometry
$v_{\mathfrak{p,p}^{\prime}}$ is defined by the formula $v_{\mathfrak{p,p}%
^{\prime}}(X):=v_{0,\mathfrak{p,p}^{\prime}}\circ X$, $X\in\mathcal{L}%
_{\mathcal{M}}(H,H_{\mathfrak{p},t})$.

We may thus form the direct limits%
\[
H_{t}:=\underrightarrow{\lim}(H_{\mathfrak{p},t},v_{0,\mathfrak{p,p}^{\prime}%
})
\]
and
\[
E(t):=\underrightarrow{\lim}(\mathcal{L}_{\mathcal{M}}(H,H_{\mathfrak{p}%
,t}),v_{\mathfrak{p,p}^{\prime}})\text{.}%
\]
Note that $H_{t}$ is a left $\mathcal{M}$ module since each $H_{\mathfrak{p}%
,t}$ is and the maps $v_{0,\mathfrak{p,p}^{\prime}}$ respect the action of
$\mathcal{M}$. It is also a left $\mathcal{M}^{\prime}$-module, since
$\mathcal{M}^{\prime}$ acts on each $H_{\mathfrak{p},t}$ via the formula
$R(T_{1}\otimes T_{2}\otimes\cdots T_{n}\otimes h)=(I\otimes R)(T_{1}\otimes
T_{2}\otimes\cdots T_{n}\otimes h)=T_{1}\otimes T_{2}\otimes\cdots
T_{n}\otimes Rh$, $T_{1}\otimes T_{2}\otimes\cdots T_{n}\otimes h\in
H_{\mathfrak{p},t}$, $R\in\mathcal{M}^{\prime}$ and the maps
$v_{0,\mathfrak{p,p}^{\prime}}$ respect this action. It is now easy to see
that $\mathcal{L}_{M}(H,H_{t})$ has the structure of an $\mathcal{M}^{\prime}%
$-correspondence. Indeed, the bimodule structure has just been indicated. One
passes to the limit when writing $H_{t}=\underrightarrow{\lim}(H_{\mathfrak{p}%
,t},v_{0,\mathfrak{p,p}^{\prime}})$ in $\mathcal{L}_{M}(H,H_{t})$. Since each
$\mathcal{L}_{\mathcal{M}}(H,H_{\mathfrak{p},t})$ is an $\mathcal{M}^{\prime}%
$-correspondence in an obvious way, so is $\mathcal{L}_{\mathcal{M}}(H,H_{t})$
via the limit of the inner products on the $\mathcal{L}_{\mathcal{M}%
}(H,H_{\mathfrak{p},t})$.

\begin{lemma}
\label{Lemma2.1}Each $E(t)$ is isomorphic, as an $\mathcal{M}^{\prime}%
$-correspondence, to $\mathcal{L}_{\mathcal{M}}(H,H_{t})$.
\end{lemma}

\begin{proof}
For each $\mathfrak{p\in P}(t)$, we write $v_{0,\mathfrak{p,\infty}}$ for the
canonical isometric embedding of $H_{\mathfrak{p},t}$ in $H_{t}$. Since the
$v_{0,\mathfrak{p,p}^{\prime}}$ are $\mathcal{M}$-module maps, so is
$v_{0,\mathfrak{p,\infty}}$. Hence we obtain $\mathcal{M}^{\prime}%
$-correspondence isometries $v_{\mathfrak{p,\infty}}:\mathcal{L}_{\mathcal{M}%
}(H,H_{\mathfrak{p},t})\longrightarrow\mathcal{L}_{\mathcal{M}}(H,H_{t})$ by
setting $v_{\mathfrak{p,\infty}}(X)=v_{0,\mathfrak{p,\infty}}\circ X.$
However, for $\mathfrak{p}^{\prime}$ finer than $\mathfrak{p}$, we have
$v_{0,\mathfrak{p}^{\prime}\mathfrak{,\infty}}\circ v_{0,\mathfrak{p,p}%
^{\prime}}=v_{0,\mathfrak{p,\infty}}$. Hence $v_{\mathfrak{p}^{\prime
}\mathfrak{,\infty}}\circ v_{\mathfrak{p,p}^{\prime}}=v_{\mathfrak{p,\infty}}%
$. Thus, by the universal properties of inductive limits, we obtain an
$\mathcal{M}^{\prime}$-correspondence isometry $v:E(t)\longrightarrow
\mathcal{L}_{\mathcal{M}}(H,H_{t})$. We need to show that $v$ is surjective.
To this end, observe that if $P$ and $Q$ are two normal, unital, completely
positive maps on $\mathcal{M}$, then using Proposition \ref{Lemma1.9} (and
applying Lemma \ref{lemma1.7}), we find that $\bigvee\{X(H)\mid X\in
\mathcal{L}_{\mathcal{M}}(H,\mathcal{M}\otimes_{Q}\mathcal{M}\otimes
_{\mathcal{P}}H)\}\supseteq\bigvee\{(I\otimes X)Y(H)\mid X\in\mathcal{L}%
_{\mathcal{M}}(H,\mathcal{M}\otimes_{P}H)$, $Y\in\mathcal{L}_{\mathcal{M}%
}(H,\mathcal{M}\otimes_{Q}H)\}=\bigvee\{(I\otimes X)(\mathcal{M}\otimes
_{Q}H)\mid X\in\mathcal{L}_{\mathcal{M}}(H,\mathcal{M}\otimes_{P}%
H)\}=\mathcal{M}\otimes_{P}\mathcal{M}\otimes_{Q}H$. The same argument,
applied to more than two maps shows that for any partition $\mathfrak{p}$ in
$\mathfrak{P}(t)$, $\bigvee\{X(H)\mid X\in\mathcal{L}_{\mathcal{M}%
}(H,H_{\mathfrak{p},t})\}=H_{\mathfrak{p},t}$ and $\bigvee
\{v_{\mathfrak{p,\infty}}(X)(H)\mid X\in\mathcal{L}_{\mathcal{M}%
}(H,H_{\mathfrak{p},t})\}=v_{0,\mathfrak{p,\infty}}(H_{\mathfrak{p}%
,t})\subseteq H_{t}$. Hence, $\bigvee\{Y(H)\mid Y\in v(E(t))\}=H_{t}$.
Consequently, given any $X\in\mathcal{L}_{\mathcal{M}}(H,H_{t})$ satisfying
$X^{\ast}Y=0$ for all $Y\in v(E(t)),$ we have $X=0$. That is, the orthogonal
complement of $v(E(t))$ in $\mathcal{L}_{\mathcal{M}}(H,H_{t})$ is zero. Since
$\mathcal{L}_{\mathcal{M}}(H,H_{t})$ is self-dual, by Proposition \ref{Lemma
1.2}, we conclude that $v(E(t))=\mathcal{L}_{\mathcal{M}}(H,H_{t})$.
\end{proof}

If $\mathfrak{p}_{1}\in\mathfrak{P}(t)$ and $\mathfrak{p}_{2}\in
\mathfrak{P}(s)$, then we shall write $\mathfrak{p}_{2}\vee\mathfrak{p}_{1}+s$
for the following partition in $\mathfrak{P}(t+s)$:%
\begin{multline*}
\{0=s_{0}<s_{1}<\cdots<s_{m-1}<s_{m}(=s=t_{0}+s)<t_{1}+s<\\
t_{2}+s\cdots<t_{n-1}+s<t_{n}+s=t+s\}\text{,}%
\end{multline*}
where $\mathfrak{p}_{1}=\{0=t_{0}<t_{1}<t_{2}<\cdots<t_{n-1}<t_{n}=t\}$ and
$\mathfrak{p}_{2}=\{0=s_{0}<s_{1}<s_{2}<\cdots<s_{n-1}<s_{m}=s\}$. Note the
order in the definition of $\mathfrak{p}_{2}\vee\mathfrak{p}_{1}+s$. The
``concatination'' of partitions is \emph{not} commutative. It is designed to
support the isomorphism of $E(s)\otimes E(t)$ with $E(t+s)$ that we are about
to describe.

\begin{lemma}
\label{Lemma2.2}Let $\mathfrak{p}_{1}\in\mathfrak{P}(t)$ and $\mathfrak{p}%
_{2}\in\mathfrak{P}(s)$ and write $\mathfrak{p}$ for $\mathfrak{p}_{2}%
\vee\mathfrak{p}_{1}+s$. Then the map that sends $X\otimes Y\in\mathcal{L}%
_{\mathcal{M}}(H,H_{\mathfrak{p}_{1},t})\otimes\mathcal{L}_{\mathcal{M}%
}(H,H_{\mathfrak{p}_{2},s})$ to $(I_{s}\otimes X)Y$ in $\mathcal{L}%
_{\mathcal{M}}(H,H_{\mathfrak{p},t+s})$ extends to an isomorphism of
$\mathcal{M}^{\prime}$-correspondences, where $I_{s}$ denotes the identity map
on $\mathcal{M}\otimes_{P_{s_{1}}}\mathcal{M}\otimes_{P_{s_{2}-s_{1}}%
}\mathcal{M}\otimes\cdots\otimes_{P_{s-s_{m-1}}}\mathcal{M}$ and where
$\mathfrak{p}_{2}=\{0=s_{0}<s_{1}<s_{2}<\cdots<s_{n-1}<s_{m}=s\}$. Further,
this isomorphism induces a natural isomorphism of $\mathcal{M}^{\prime}%
$-correspondences from $E(t)\otimes E(s)$ onto $E(t+s)$.
\end{lemma}

\begin{proof}
That the map $X\otimes Y\rightarrow(I_{s}\otimes X)Y$ induces an isomorphism
of $\mathcal{M}^{\prime}$-correspondences from $\mathcal{L}_{\mathcal{M}%
}(H,H_{\mathfrak{p}_{1},t})\otimes\mathcal{L}_{\mathcal{M}}(H,H_{\mathfrak{p}%
_{2},s})$ into $\mathcal{L}_{\mathcal{M}}(H,H_{\mathfrak{p},t+s})$ is
essentially proved in Proposition \ref{Lemma1.9}. To see that the isomorphism
is surjective, simply apply Corollary \ref{Cor 1.10} (several times).

To get the isomorphism from $E(t)\otimes E(s)$ onto $E(t+s)$, we appeal to
universal properties of inductive limits. Let $\mathfrak{p}_{1}$ and
$\mathfrak{p}_{1}^{\prime}$ be partitions in $\mathfrak{P}(t)$, with
$\mathfrak{p}_{1}^{\prime}$ finer than $\mathfrak{p}_{1}$, and let
$\mathfrak{p}_{2}$ be a partition in $\mathfrak{P}(s)$. Write $\mathfrak{p}%
=\mathfrak{p}_{2}\vee\mathfrak{p}_{1}+s$ and $\mathfrak{p}^{\prime
}=\mathfrak{p}_{2}\vee\mathfrak{p}_{1}^{\prime}+s$. Also let $\alpha
_{\mathfrak{p}_{1},\mathfrak{p}_{2}}$ be the isomorphism from $\mathcal{L}%
_{\mathcal{M}}(H,H_{\mathfrak{p}_{1},t})\otimes\mathcal{L}_{\mathcal{M}%
}(H,H_{\mathfrak{p}_{2},s})$ onto $\mathcal{L}_{\mathcal{M}}(H,H_{\mathfrak{p}%
,t+s})$ that sends $X\otimes Y$ to $(I_{s}\otimes X)Y$, and let $\alpha
_{\mathfrak{p}_{1}^{\prime},\mathfrak{p}_{2}}$ be the similarly defined
isomorphism from $\mathcal{L}_{\mathcal{M}}(H,H_{\mathfrak{p}_{1}^{\prime}%
,t})\otimes\mathcal{L}_{\mathcal{M}}(H,H_{\mathfrak{p}_{2},s})$ onto
$\mathcal{L}_{\mathcal{M}}(H,H_{\mathfrak{p}^{\prime},t+s})$. Then we have the
following diagram, which is easily seen to be commutative:%
\[%
\begin{array}
[c]{ccc}%
\mathcal{L}_{\mathcal{M}}(H,H_{\mathfrak{p}_{1},t})\otimes\mathcal{L}%
_{\mathcal{M}}(H,H_{\mathfrak{p}_{2},s}) & \overset{\alpha_{\mathfrak{p}%
_{1},\mathfrak{p}_{2}}}{\longrightarrow} & \mathcal{L}_{\mathcal{M}%
}(H,H_{\mathfrak{p},t+s})\\%
\begin{array}
[c]{ccc}%
v_{\mathfrak{p}_{1}\mathfrak{,p}_{1}^{\prime}}\otimes I & \downarrow &
\;\;\;\;\;\;\;\;\;\;\;\;\;\;\;
\end{array}
&  &
\begin{array}
[c]{ccc}%
\;\;\;\;\;\;\; & \downarrow & v_{\mathfrak{p,p}^{\prime}}%
\end{array}
\\
\mathcal{L}_{\mathcal{M}}(H,H_{\mathfrak{p}_{1}^{\prime},t})\otimes
\mathcal{L}_{\mathcal{M}}(H,H_{\mathfrak{p}_{2},s}) & \overset{\alpha
_{\mathfrak{p}_{1}^{\prime},\mathfrak{p}_{2}}}{\longrightarrow} &
\mathcal{L}_{\mathcal{M}}(H,H_{\mathfrak{p}^{\prime},t+s})
\end{array}
\]
In the limit, we obtain an isometry from $E(t)\otimes\mathcal{L}_{\mathcal{M}%
}(H,H_{\mathfrak{p}_{2},s})$ into $E(t+s)$. A similar argument yields an
isometry from $E(t)\otimes E(s)$ into $E(t+s)$. It is clear from the
definition of this map that its image contains all the $\mathcal{L}%
_{\mathcal{M}}(H,H_{\mathfrak{p},t+s})$, where $\mathfrak{p}$ is constructed
as $\mathfrak{p}_{2}\vee\mathfrak{p}_{1}+s$. (We shall view these spaces as
contained in $E(t+s)$ without reference to the isomorphic embeddings.) For a
given partition $\mathfrak{p\in P}(s+t)$, we can refine it by adding $s$ to
get $\mathfrak{p}^{\prime}$, say. Then $\mathcal{L}_{\mathcal{M}%
}(H,H_{\mathfrak{p}^{\prime},t+s})$ is contained in $E(t+s)$ and contains (a
copy of) $\mathcal{L}_{\mathcal{M}}(H,H_{\mathfrak{p},t+s})$. Hence the image
contains all the $\mathcal{L}_{\mathcal{M}}(H,H_{\mathfrak{p},t+s})$ and so
must be all of $E(t+s)$.
\end{proof}

\begin{remark}
\label{Associativity} Given $t,\;s,\;r\in(0,\infty)$ and partitions
$\mathfrak{p}_{1}\in\mathfrak{P}(t)$, $\mathfrak{p}_{2}\in\mathfrak{P}(s)$,
and $\mathfrak{p}_{3}\in\mathfrak{P}(r)$, one can define an isomorphism of
$\mathcal{M}^{\prime}$-correspondences between $\mathcal{L}_{\mathcal{M}%
}(H,H_{\mathfrak{p}_{1},t})\otimes\mathcal{L}_{\mathcal{M}}(H,H_{\mathfrak{p}%
_{2},s})\otimes\mathcal{L}_{\mathcal{M}}(H,H_{\mathfrak{p}_{3},r})$ and
$\mathcal{L}_{\mathcal{M}}(H,H_{\mathfrak{p},t+s+r})$ in two different, but
natural, ways, where $\mathfrak{p}=\mathfrak{p}_{3}\vee(\mathfrak{p}%
_{2}+r)\vee(\mathfrak{p}_{1}+s+r)$: In the first, we map the left hand side,
$\mathcal{L}_{\mathcal{M}}(H,H_{\mathfrak{p}_{1},t})\otimes\mathcal{L}%
_{\mathcal{M}}(H,H_{\mathfrak{p}_{2},s})\otimes\mathcal{L}_{\mathcal{M}%
}(H,H_{\mathfrak{p}_{3},r})$, to $\mathcal{L}_{\mathcal{M}}(H,H_{\mathfrak{p}%
_{1},t})\otimes\mathcal{L}_{\mathcal{M}}(H,H_{\mathfrak{p}^{\prime},t})$,
where $\mathfrak{p}^{\prime}=\mathfrak{p}_{3}\vee(\mathfrak{p}_{2}+r)$, and
then to $\mathcal{L}_{\mathcal{M}}(H,H_{\mathfrak{p},t+s+r})$, while in the
second, we map $\mathcal{L}_{\mathcal{M}}(H,H_{\mathfrak{p}_{1},t}%
)\otimes\mathcal{L}_{\mathcal{M}}(H,H_{\mathfrak{p}_{2},s})\otimes
\mathcal{L}_{\mathcal{M}}(H,H_{\mathfrak{p}_{3},r})$ to $\mathcal{L}%
_{\mathcal{M}}(H,H_{\mathfrak{p}^{\prime\prime},t+s})\otimes\mathcal{L}%
_{\mathcal{M}}(H,H_{\mathfrak{p}_{3},r})$, where $\mathfrak{p}^{\prime\prime
}=\mathfrak{p}_{2}\vee(\mathfrak{p}_{1}+s)$ and then to $\mathcal{L}%
_{\mathcal{M}}(H,H_{\mathfrak{p},t+s+r})$. These two ways of identifying
$\mathcal{L}_{\mathcal{M}}(H,H_{\mathfrak{p}_{1},t})\otimes\mathcal{L}%
_{\mathcal{M}}(H,H_{\mathfrak{p}_{2},s})\otimes\mathcal{L}_{\mathcal{M}%
}(H,H_{\mathfrak{p}_{3},r})$ and $\mathcal{L}_{\mathcal{M}}(H,H_{\mathfrak{p}%
,t+s+r})$ amount to nothing more than identifying \thinspace$X_{1}\otimes
X_{2}\otimes X_{3}$ with $(I_{s+r}\otimes X_{1})\circ(I_{r}\otimes X_{2})\circ
X_{3}$, as we may. Passing to the limit yields the natural isomorphisms%
\[
(E(t)\otimes E(s))\otimes E(r)\simeq E(t)\otimes(E(s)\otimes E(r))\simeq
E(t+s+r)\text{.}%
\]
\end{remark}

Our analysis to this point shows that if we set $E(0)=\mathcal{M}^{\prime}$,
then $\{E(t)\}_{t\geq0}$ is a discrete product system in the sense of

\begin{definition}
\label{Definition2.3a}Let $\mathcal{N}$ be a von Neumann algebra. A
\emph{discrete product system} over $\mathcal{N}$ is simply a family
$\{E(t)\}_{t\geq0}$ of $W^{\ast}$-correspondences over $\mathcal{N}$ such that
$E(0)=\mathcal{N}$ and such that $E(t+s)\simeq E(t)\otimes E(s)$ for all
$t,s\in\lbrack0,\infty)$.

The particular product system that we associated with the semigroup
$\{P_{t}\}_{t\geq0}$ in the preceding paragraphs will be called \emph{the}
\emph{(discrete) product system of }$M^{\prime}$\emph{-correspondences}
associated with $\{P_{t}\}_{t\geq0}$.

A (completely contractive) \emph{covariant representation} of a discrete
product system $\{E(t)\}_{t\geq0}$ on a Hilbert space $H$ is simply a family
$\{T_{t}\}_{t\geq0}$ of completely contractive linear maps, where $T_{t}$ maps
from $E(t)$ to $\mathcal{B}(H)$ such that each $T_{t}$ is continuous with
respect to the $\sigma$-topology on $E(t)$ and the $\sigma$-weak topology on
$\mathcal{B}(H)$, $T_{0}$ is a $\ast$-representation of $E(0)=\mathcal{N}$ on
$H$, and such that $T_{t}\otimes T_{s}=T_{t+s}$ (after identifying
$E(t+s)\simeq E(t)\otimes E(s)$.)
\end{definition}

\begin{remark}
It is useful to think of product systems as semigroups and then to view
covariant representations as representations of such a semigroup. However,
when working with any particular product system and representation, it
frequently becomes necessary to make explicit the isomorphisms between
$E(t)\otimes E(s)$ and $E(t+s)$ and then, of course, the formulas involving
$\{T_{t}\}_{t\geq0}$ become correspondingly more complicated.

Note, too, that the definition of a covariant representation implies that
$T_{t}(a\xi b)\allowbreak=T_{0}(a)T_{t}(\xi)T_{0}(b)$ for all $t\geq0$,
$\xi\in E(t)$, $a,b\in\mathcal{N}$. Thus if $\{T_{t}\}_{t\geq0}$ is a
covariant representation of the product system then for each $t$,
$(T_{t},T_{0})$ is a completely contractive covariant representation of $E(t)$
in the sense of Definition \ref{Definition1.12}.
\end{remark}

\begin{definition}
\label{Definition2.3b} A covariant representation $\{T_{t}\}_{t\geq0}$ of a
product system $\{E(t)\}_{t\geq0}$ is called \emph{isometric} in case for each
$t$, $(T_{t},T_{0})$ is isometric in the sense of Definition
\ref{Definition1.12}. It is called \emph{fully coisometric} in case for each
$t$, $(T_{t},T_{0})\;$is fully coisometric in the sense of Definition
\ref{Definition1.12bis}.
\end{definition}

Our next objective is to show how a fully coisometric covariant representation
of a product system $\{E(t)\}_{t\geq0}$ can be dilated to a fully coisometric
and isometric representation of $\{E(t)\}_{t\geq0}$.

\begin{theorem and definition}
\label{Theorem2.6}Let $\{E(t)\}_{t\geq0}$ be a discrete product system over a
von Neumann algebra $\mathcal{N}$ and let $\{T_{t}\}_{t\geq0}$ be a
\emph{fully coisometric }covariant representation of $\{E(t)\}_{t\geq0}$ on a
Hilbert space $H$. Then there is another Hilbert space $K$, an isometry
$u_{0}$ mapping $H$ into $K$, and fully coisometric, isometric covariant
representation $\{V_{t}\}_{t\geq0}$ of $E$ on $K$ so that

\begin{enumerate}
\item $u_{0}^{\ast}V_{t}(\xi)u_{0}=T_{t}(\xi)$ for all $\xi\in E(t)$,
$t\geq0;$ and

\item For $\xi\in E(t)$, $t\geq0$, $V_{t}(\xi)^{\ast}$ leaves $u_{0}(H)$ invariant.
\end{enumerate}

The smallest subspace of $K$ containing $u_{0}(H)$ and reducing $V_{t}(\xi)$
for every $\xi\in E(t)$, $t\geq0$, is all of $K$. If $(\{V_{t}^{\prime
}\}_{t\geq0},u_{0}^{\prime},K^{\prime})$ is another triple with same
properties as $(\{V_{t}\}_{t\geq0},u_{0},K)$, then there is a Hilbert space
isomorphism $W$ from $K$ to $K^{\prime}$ such that $WV_{t}(\xi)W^{-1}%
=V_{t}^{^{\prime}}(\xi)$ for all $\xi\in E(t)$ and $t\geq0$, and $W\circ
u_{0}=u_{0}^{\prime}$. We therefore call the triple, $(\{V_{t}\}_{t\geq
0},u_{0},K)$, \emph{the} \emph{minimal isometric dilation of }$\{T_{t}%
\}_{t\geq0}$.
\end{theorem and definition}

\begin{proof}
For $0\leq t<s$, we write $U_{t,s}$ for the isomorphism from $E(t)\otimes
E(s-t)$ to $E(s).$ Then the associativity of tensor products implemented
through these isomorphisms coupled with the identification of $E(t)\otimes
E(s)\otimes E(r)$ with $E(t+s+r)$ imply that $U_{s,r}(U_{t,s}\otimes
I_{r-s})=U_{t,r}$.\ Further, for any $t$, we write $\widetilde{T}_{t}$ for the
operator from $E(t)\otimes_{T_{0}}H$ to $H$ defined by the formula
$\widetilde{T}_{t}(\xi\otimes h)=T_{t}(\xi)h$. (See Lemma \ref{CovRep} and the
discussion surrounding it.) For $0\leq t<s$, we define $u_{t,s}$ from
$E(t)\otimes_{T_{0}}H$ to $E(s)\otimes_{T_{0}}H$ by the formula%
\[
u_{t,s}:=(U_{t,s}\otimes I_{H})(I_{E(t)}\otimes\widetilde{T}_{s-t}^{\ast
})\text{.}%
\]
Observe that each space $E(t)\otimes_{T_{0}}H$ is a left $\mathcal{N}$-module
and that the $u_{t,s}$ are $\mathcal{N}$-module maps.

We claim that each $u_{t,s}$ is an isometry. Indeed, since $U_{t,s}$ is a
Hilbert module isomorphism, $U_{t,s}\otimes I_{H}$ is a Hilbert space
isomorphism, i.e., a unitary, and so%
\begin{align*}
u_{t,s}^{\ast}u_{t,s}  &  =(I_{E(t)}\otimes\widetilde{T}_{s-t})(U_{t-s}\otimes
I_{H})^{\ast}(U_{t-s}\otimes I_{H})(I_{E(t)}\otimes\widetilde{T}_{s-t}^{\ast
})\\
&  =(I_{E(t)}\otimes\widetilde{T}_{s-t})(I_{E(t)}\otimes\widetilde{T}%
_{s-t}^{\ast})\\
&  =I_{E(t)}\otimes\widetilde{T}_{s-t}\widetilde{T}_{s-t}^{\ast}\text{.}%
\end{align*}
However, this last term is the identity on $E(t)\otimes_{T_{0}}H$ because
$\{T_{t}\}_{t\geq0}$ is assumed to be fully coisometric.

Further observe that the composition properties of the $U_{t,s}$ coupled with
the fact that $\{T_{t}\}_{t\geq0}$ is a covariant representation imply that
for $0\leq t<s<r$, $u_{t,s}u_{s,r}=u_{t,r}$. Hence, if we agree to set
$u_{t,t}$ equal to the identity on $E(t)\otimes_{T_{0}}H$ for each $t$, then
$\{\{E(t)\otimes_{T_{0}}H\}_{t\geq0},\{u_{t,s}\}_{0\leq t\leq s}\}$ is an
inductive system of Hilbert spaces - in fact, it is an inductive system of
$\mathcal{N}$-modules and module maps. We set $K=\underrightarrow{\lim
}(E(t)\otimes_{T_{0}}H,u_{t,s})$ and we write $u_{t}:E(t)\otimes_{T_{0}%
}H\longrightarrow K$ for the canonical embeddings. Note that the $u_{t}$'s are
isometries and $\mathcal{N}$-module maps.

To construct the dilation $\{V_{t}\}_{t\geq0}$, we begin by defining $V_{t}$
on the range of each $u_{s}$ by the formula%
\[
V_{t}(\xi)\cdot u_{s}(\eta\otimes h):=u_{t+s}(U_{t,t+s}(\xi\otimes\eta)\otimes
h)\text{,}%
\]
where $t,s\geq0$, $\xi\in E(t)$, $\eta\in E(s)$, and $h\in H$. To see that
$V_{t}(\xi)$ is well-defined on the union of the ranges of the $u_{s}$, simply
note that for $s_{1}>s_{2}\geq0$, $t\geq0$, $\eta\in E(s_{2})$, $\xi\in E(t)$,
and $h\in H$, we have
\[
V_{t}(\xi)u_{s_{1}}u_{s_{2},s_{1}}(\eta\otimes h)=V_{t}(\xi)u_{s_{2}}%
(\eta\otimes h)\text{.}%
\]
The $\mathcal{N}$-module structure on $K$ is just that afforded by $V_{0}$.
That is, for $a\in\mathcal{N}$, $V_{0}(a)\cdot u_{s}(\eta\otimes
h)=u_{s}(a\eta\otimes h),\;\eta\otimes h\in E(s)\otimes_{T_{0}}H$. So, to show
that every other $V_{t}$ extends to all of $K$, and yields an isometric
representation of the $\mathcal{N}$-correspondence $E(t)$, we first simply
compute to see that for $\xi_{1},\xi_{2}\in E(t)$, and $\eta\otimes
h,\zeta\otimes k\in E(s)\otimes_{T_{0}}H$, we have%
\begin{multline*}
\langle V_{t}^{\ast}(\xi_{1})V_{t}(\xi_{2})u_{s}(\eta\otimes h),u_{s}%
(\zeta\otimes k)\rangle=\langle u_{t+s}(\xi_{2}\otimes\eta\otimes
h),u_{t+s}(\xi_{1}\otimes\eta\otimes k)\rangle\\
=\langle\xi_{2}\otimes\eta\otimes h,\xi_{1}\otimes\eta\otimes k\rangle
=\langle\eta\otimes h,\langle\xi_{2},\xi_{1}\rangle\eta\otimes k\rangle
=\langle u_{s}(\eta\otimes h),V_{0}(\langle\xi_{2},\xi_{1}\rangle)u_{s}%
(\eta\otimes k)\rangle\text{.}%
\end{multline*}
This shows that on the range of each $u_{s}$ $(V_{t},V_{0})$ is an isometric
covariant representation of $E(t)$. Thus on the range of each $u_{s}$,
$V_{t}(\xi)$ is a bounded operator with norm bounded by $\left\|  \xi\right\|
$. Hence, $V_{t}(\xi)$ extends to all of $K$ as a bounded operator. Further,
this equation shows that if we denote the projection of $K$ onto the range of
$u_{s}$ by $Q_{s}$, i.e., if we let $Q_{s}=u_{s}u_{s}^{\ast}$, then
\[
Q_{s}(V_{t}(\xi_{2})^{\ast}V_{t}(\xi_{1}))Q_{s}=Q_{s}V_{0}(\langle\xi_{2}%
,\xi_{1}\rangle)Q_{s}\text{.}%
\]
Since this so for all $s$, it follows that $V_{t}(\xi_{2})^{\ast}V_{t}(\xi
_{1})=V_{0}(\langle\xi_{2},\xi_{1}\rangle)$ on all of $K$. Thus, for each $t$,
$(V_{t},V_{0})$ is an isometric covariant representation of $E(t)$.

To show that $\{V_{t}\}_{t\geq0}$ satisfies the semi-group property, let
$t=t_{1}+t_{2}$, let $\xi_{1}\in E(t_{1})$, $\xi_{2}\in E(t_{2})$, and let
$\eta\otimes h\in E(s)\otimes_{T_{0}}H$. Then on the one hand we have%
\[
V_{t}(U_{t_{1},t}(\xi_{1}\otimes\xi_{2}))u_{s}(\eta\otimes h)=u_{t+s}%
(U_{t,t+s}(U_{t_{1},t}(\xi_{1}\otimes\xi_{2})\otimes\eta)\otimes h)\text{,}%
\]
while on the other we have%
\begin{align*}
V_{t_{1}}(\xi_{1})V_{t_{2}}(\xi_{2})u_{s}(\eta\otimes h)  &  =V_{t_{1}}%
(\xi_{1})u_{t_{2}+s}U_{t_{2},t_{2}+s}(\xi_{2}\otimes\eta)\otimes h)\\
&  =u_{t+s}(U_{t_{1},t+s}(\xi_{1}\otimes U_{t_{2},t_{2}+s}(\xi_{2}\otimes
\eta))\otimes h)\text{.}%
\end{align*}
By Remark \ref{Associativity}, we conclude that $V_{t}(U_{t_{1},t}(\xi
_{1}\otimes\xi_{2}))=V_{t_{1}}(\xi_{1})V_{t_{2}}(\xi_{2})$. Ignoring the
$U_{t,t+s}$ when identifying $E(t)\otimes E(s)$ with $E(t+s)$, we obtain the
desired result: $V_{t}(\xi_{1}\otimes\xi_{2})=V_{t_{1}}(\xi_{1})V_{t_{2}}%
(\xi_{2})$.

Next, we show that each $V_{t}$ is continuous with respect to the $\sigma
$-topology on $E(t)$ and the $\sigma$-weak topology on $\mathcal{B}(K)$. For
this, observe that for $\xi,\xi_{1}\in E(t)$, $\eta,\xi_{2}\in E(s)$, and
$h,k\in H$, we have $\langle V_{t}(\xi)u_{s}(\eta\otimes h),u_{t+s}(\xi
_{1}\otimes\xi_{2}\otimes k)\rangle=\langle u_{t+s}(\xi\otimes\eta\otimes
h),u_{t+s}(\xi_{1}\otimes\xi_{2}\otimes k)\rangle=\langle h,T_{0}(\langle
\xi\otimes\eta,\xi_{1}\otimes\xi_{2}\rangle)k\rangle=\langle h,T_{0}%
(\langle\eta,\langle\xi,\xi_{1}\rangle\xi_{2}\rangle)k\rangle$. Thus, for each
$s\geq0$, the map $\xi\mapsto V_{t}(\xi)|u_{s}(E(s)\otimes_{T_{0}}H)$ has the
desired continuity properties. Since the union of the ranges of the $u_{s}$ is
dense and since $\left\|  V_{t}(\xi)\right\|  \leq\left\|  \xi\right\|  $, we
conclude that $V_{t}$ is continuous with respect to the $\sigma$-topology on
$E(t)$ and the $\sigma$-weak topology on $\mathcal{B}(K)$.

To see that $u_{0}^{\ast}V_{t}(\xi)u_{0}=T_{t}(\xi)$, i.e., to see that
$\{V_{t}\}_{t\geq0}$ dilates $\{T_{t}\}_{t\geq0}$, simply note that for
$h,h^{\prime}\in H,\;t>0$ and $\xi\in E(t)$, we have $\langle u_{0}^{\ast
}V_{t}(\xi)u_{0}(h),h^{\prime}\rangle=\langle u_{t}(\xi\otimes h),u_{0}%
(h^{\prime})\rangle=\langle u_{t}(\xi\otimes h),u_{t}(u_{0,t}(h^{\prime
}))\rangle=\langle\xi\otimes h,\widetilde{T}_{t}^{\ast}h^{\prime}%
\rangle_{E(t)\otimes H}=\langle\widetilde{T}_{t}(\xi\otimes h),h^{\prime
}\rangle=\langle T_{t}(\xi)h,h^{\prime}\rangle$.

To check that $V_{t}(\xi)^{\ast}$, $\xi\in E(t)$, leaves $u_{0}(H)$ invariant,
first note that the computation just completed shows that for $\zeta\in E(r)$,
$r\geq0$, and $h\in H$, $u_{0}^{\ast}u_{r}(\zeta\otimes h)=T_{r}(\zeta)h$.
Hence, for $\xi\in E(t)$, $\eta\in E(s)$, and $h\in H$, $u_{0}^{\ast}V_{t}%
(\xi)u_{s}(\eta\otimes h)=u_{0}^{\ast}u_{t+s}(\xi\otimes\eta\otimes
h)=T_{t+s}(\xi\otimes\eta)h=T_{t}(\xi)T_{t}(\eta)h=u_{0}^{\ast}u_{t}%
(\xi\otimes T_{s}(\eta)h)=u_{0}^{\ast}V_{t}(\xi)u_{0}u_{0}^{\ast}u_{s}%
(\eta\otimes h).$ Since this holds for all $s\geq0,$ we see that $u_{0}^{\ast
}V_{t}(\xi)=u_{0}^{\ast}V_{t}(\xi)u_{0}u_{0}^{\ast}.$ Taking adjoints and
multiplying the resulting equation on the right by $u_{0}$, we conclude that
$V_{t}(\xi)^{\ast}u_{0}u_{0}^{\ast}=u_{0}u_{0}^{\ast}V_{t}(\xi)^{\ast}%
u_{0}u_{0}^{\ast}$ for all $\xi\in E(t)$, and $t\geq0$, which shows that
$V_{t}(\xi)^{\ast}$, $\xi\in E(t)$, leaves $u_{0}(H)$ invariant.

To see that $\{V_{t}\}_{t\geq0}$ is fully coisometric because $\{T_{t}%
\}_{t\geq0}$ is, we need to show that $\widetilde{V}_{t}$ is a coisometry for
each $t$. Since $\{V_{t}\}_{t\geq0}$ is isometric, each $\widetilde{V}_{t}$ is
an isometry. Hence, all we need to do is to show that the range of each
$\widetilde{V}_{t}$ is dense. For this, it suffices to show that for every
$s\geq0$ the span of $\{V_{t}(\xi)u_{s}(\eta\otimes h)\mid\eta\in E(s)$,
$\xi\in E(t)$, $h\in H\}$ equals $u_{t+s}(E(t+s)\otimes H)$. However,
$V_{t}(\xi)u_{s}(\eta\otimes h)=u_{t+s}(\xi\otimes\eta\otimes h)$ and, since
$E(t)\otimes E(s)$ is isomorphic to $E(t+s)$, we see that the range of
$\widetilde{V}_{t}$ is, indeed, dense.

>From what we have shown so far, it is clear that the smallest subspace of $K$
that contains $u_{0}(H)$ and reduces every $V_{t}(\xi)$ is all of $K$. The
uniqueness of $(\{V_{t}\}_{t\geq0},u_{0},K)$ up to unitary equivalence is
proved just as in Proposition 3.2 of \cite{MS98}, and so will be omitted here.
\end{proof}

\begin{remark}
\label{alternate}It is worthwhile pointing out that the relation $u_{0}^{\ast
}V_{t}(\xi)u_{0}=T_{t}(\xi)$ in the preceding theorem is equivalent to the
relation%
\[
u_{0}^{\ast}\widetilde{V}_{t}(I\otimes u_{0})=\widetilde{T}_{t}\text{.}%
\]
Further, the invariance of $u_{0}(H)$ under $V_{t}(\xi)^{\ast}$ is equivalent
to the equation $u_{0}u_{0}^{\ast}\widetilde{V}_{t}(I\otimes u_{0}u_{0}^{\ast
})=u_{0}u_{0}^{\ast}\widetilde{V}_{t}$. These assertions are immediate from
the proof.
\end{remark}

We return to our semigroup, $\{P_{t}\}_{t\geq0}$, of completely positive maps
on the von Neumann algebra $\mathcal{M}$ and to the associated product system
$\mathcal{M}^{\prime}$-correspondences $\{E(t)\}_{t\geq0}$ that we constructed
at the outset of this section. Our next objective, Theorem \ref{Lemma2.8}, is
to show that there is a fully coisometric, completely contractive covariant
representation $\{T_{t}\}_{t\geq0}$ of $\{E(t)\}_{t\geq0}$ on $H$ (the Hilbert
space of $\mathcal{M}$) so that $\{P_{t}\}_{t\geq0}$ can be represented by the
formula%
\[
P_{t}(S)=\widetilde{T}_{t}(I_{E(t)}\otimes S)\widetilde{T}_{t}^{\ast}\text{,}%
\]
$S\in\mathcal{M}$, $t\geq0$.

For this purpose, recall that for a partition $\mathfrak{p}\in\mathfrak{P}%
(t)$, the Hilbert space $H_{\mathfrak{p},t}$ is $H_{\mathfrak{p}%
,t}:=\mathcal{M}\otimes_{P_{t_{1}}}\mathcal{M}\otimes_{P_{t_{2}-t_{1}}}%
\otimes\cdots\mathcal{M}\otimes_{P_{t-t_{n-1}}}H$ where $\mathfrak{p}%
=\{0=t_{0}<t_{1}<t_{2}<\cdots<t_{n-1}<t_{n}=t\}$. The map $\iota
_{\mathfrak{p}}:H\longrightarrow H_{\mathfrak{p},t}$, defined by the formula
$\iota_{\mathfrak{p}}(h)=I\otimes I\otimes\cdots\otimes I\otimes h$ is easily
seen to be an isometry, with adjoint $\iota_{\mathfrak{p}}^{\ast}$ given by
the formula%
\[
\iota_{\mathfrak{p}}^{\ast}(X_{1}\otimes X_{2}\otimes\cdots\otimes
X_{n}\otimes h)=P_{t-t_{n-1}}(P_{t_{n-1}-t_{n-2}}(\cdots(P_{t_{1}}(X_{1}%
)X_{2})\cdots)X_{n-1})X_{n})h\text{.}%
\]
Indeed, $\iota_{\mathfrak{p}}$ is just a generalization of the Stinespring
embedding\ $W_{P}$ for a single completely positive map and the the formula
for $\iota_{\mathfrak{p}}^{\ast}$ is an obvious extension of formula
(\ref{Wpstar}). Further, it is easy to check that if $\mathfrak{p}^{\prime}$
is a refinement of $\mathfrak{p}$ in $\mathfrak{P}(t)$, then $\iota
_{\mathfrak{p}}^{\ast}=\iota_{\mathfrak{p}^{\prime}}^{\ast}\circ
v_{0,\mathfrak{p},\mathfrak{p}^{\prime}}$. Hence, by the universal properties
of inductive limits, there is a (unique) map $\iota_{t}^{\ast}:H_{t}%
\;(=\underrightarrow{\lim}(H_{\mathfrak{p},t},v_{0,\mathfrak{p,p}^{\prime}%
}))\longrightarrow H$ so that $\iota_{t}^{\ast}v_{0,\mathfrak{p,\infty}}%
=\iota_{\mathfrak{p}}^{\ast}$, where, recall, $v_{0,\mathfrak{p,\infty}%
}:H_{\mathfrak{p},t}\longrightarrow H_{t}$ is the canonical isometric
embedding associated with the directed system $(H_{\mathfrak{p},t}%
,v_{0,\mathfrak{p,p}^{\prime}})$ and its limit, $H_{t}$. It is easy to check
that $\iota_{t}^{\ast}$ is a coisometry.

To define the covariant representation $\{T_{t}\}_{t\geq0}$ of
$\{E(t)\}_{t\geq0}$ that we want, we recall that $E(t)$ is isomorphic to
$\mathcal{L}_{\mathcal{M}}(H,H_{t})$ and we set%
\begin{equation}
T_{t}(X)=\iota_{t}^{\ast}\circ X\text{,} \label{idcovrep}%
\end{equation}
for $X\in\mathcal{L}_{\mathcal{M}}(H,H_{t})$.

\begin{theorem and definition}
\label{Lemma2.8}Let $\{E(t)\}_{t\geq0}$ be the discrete product system of
$\mathcal{M}^{\prime}$-correspondences constructed from $\{P_{t}\}_{t\geq0}$
as above, and let $\{T_{t}\}_{t\geq0}$ be defined by equation (\ref{idcovrep}%
). Then $\{T_{t}\}_{t\geq0}$ is a fully coisometric, completely contractive
covariant representation of $\{E(t)\}_{t\geq0}$ such that
\begin{equation}
P_{t}(S)=\widetilde{T}_{t}(I_{E(t)}\otimes S)\widetilde{T}_{t}^{\ast}%
\text{,}\label{implement}%
\end{equation}
for all $t\geq0$ and all $S\in\mathcal{M}$. We call $\{T_{t}\}_{t\geq0}$
\emph{the identity representation of} $\{E(t)\}_{t\geq0}$.
\end{theorem and definition}

\begin{proof}
Since $T_{t}$ is given by left multiplication by an operator between Hilbert
spaces of norm at most one, viz. $\iota_{t}^{\ast}$, $T_{t}$ is completely
contractive. To check that $\{T_{t}\}_{t\geq0}$ is multiplicative, we identify
$E(t+s)$ with $E(t)\otimes E(s)$ as above and proceed to show that under this
identification, $T_{t+s}=T_{t}\otimes T_{s}$. For this purpose, let
$\mathfrak{p}_{1}$ be a partition in $\mathfrak{P}(t)$ and let $\mathfrak{p}%
_{2}$ be a partition in $\mathfrak{P}(s)$. Also, fix $X_{1}\in\mathcal{L}%
_{\mathcal{M}}(H,H_{\mathfrak{p}_{1},t})$ and $X_{2}\in\mathcal{L}%
_{\mathcal{M}}(H,H_{\mathfrak{p}_{2},t})$. Then the map that sends
$X_{1}\otimes X_{2}$ to $(I_{s}\otimes X_{1})X_{2}$ (where $I_{s}$ denotes the
identity map on $\mathcal{M}\otimes_{P_{s_{1}}}\otimes\cdots\otimes
_{P_{s-s_{j-1}}}\mathcal{M}$, with $\mathfrak{p}_{2}=\{0=s_{0}<s_{1}%
<s_{2}<\cdots<s_{j-1}<s_{j}=s\}$) carries $\mathcal{L}_{\mathcal{M}%
}(H,H_{\mathfrak{p}_{1},t})\otimes\mathcal{L}_{\mathcal{M}}(H,H_{\mathfrak{p}%
_{2},t})$ to $\mathcal{L}_{\mathcal{M}}(H,H_{\mathfrak{p},t+s})$, where
$\mathfrak{p}=\mathfrak{p}_{2}\vee\mathfrak{p}_{1}+s$. To show that
$T_{t+s}=T_{t}\otimes T_{s}$, it follows from equation (\ref{idcovrep}) and
the properties of direct limits that we need only check that $\iota
_{\mathfrak{p}}^{\ast}\circ(I_{s}\otimes X_{1})X_{2}=\iota_{\mathfrak{p}_{1}%
}^{\ast}X_{1}\circ\iota_{\mathfrak{p}_{2}}^{\ast}X_{2}$. For this purpose,
consider the element $S_{1}\otimes S_{2}\otimes\cdots\otimes S_{j}\otimes
T_{1}\otimes T_{2}\otimes\cdots\otimes T_{n}\otimes k$ in $H_{\mathfrak{p}%
,t+s}$ (where $\mathfrak{p}$ and $\mathfrak{p}_{2}$ have just been defined and
$\mathfrak{p}_{1}=\{0=t_{0}<t_{1}<t_{2}<\cdots<t_{n-1}<t_{n}=t\}$). Then
\begin{multline*}
\iota_{\mathfrak{p}}^{\ast}(S_{1}\otimes S_{2}\otimes\cdots\otimes
T_{n}\otimes k)\\
=P_{t-t_{n-1}}(P_{t_{n-1}-t_{n-2}}(\cdots(P_{s-s_{j-1}}(\cdots(P_{s_{1}}%
(S_{1})S_{2})\cdots T_{n})k\\
=\iota_{\mathfrak{p}_{1}}^{\ast}((P_{s-s_{j-1}}(\cdots(P_{s_{1}}(S_{1}%
)S_{2})\cdots(S_{j})T_{1}\otimes T_{2}\otimes\cdots\otimes T_{n}\otimes k).
\end{multline*}
Hence, for $X_{1}(h)\in H_{\mathfrak{p}_{1},t}$, and $S_{1}$, $S_{2}$, \ldots,
$S_{j}$ in $\mathcal{M}$,
\[
\iota_{\mathfrak{p}}^{\ast}(S_{1}\otimes S_{2}\otimes\cdots\otimes
S_{j}\otimes X_{1}(h))=\iota_{\mathfrak{p}_{1}}^{\ast}(P_{s-s_{j-1}}%
(\cdots(P_{s_{1}}(S_{1})S_{2})\cdots S_{j})X_{1}(h)).
\]
Since $X_{1}$ is an $\mathcal{M}$-module map, this equation can be rewritten
as
\begin{align*}
\iota_{\mathfrak{p}}^{\ast}(S_{1}\otimes S_{2}\otimes\cdots\otimes
S_{j}\otimes X_{1}(h)) &  =\iota_{\mathfrak{p}_{1}}^{\ast}X_{1}(P_{s-s_{j-1}%
}(\cdots(P_{s_{1}}(S_{1})S_{2})\cdots S_{j})h)\\
&  =\iota_{\mathfrak{p}_{1}}^{\ast}X_{1}\iota_{\mathfrak{p}_{2}}^{\ast}%
(S_{1}\otimes S_{2}\otimes\cdots\otimes S_{j}\otimes h),
\end{align*}
as was required.

To see that each $T_{t}$ is continuous with respect to the $\sigma$-topology
on $E(t)$ and the $\sigma$-weak topology on $\mathcal{B}(H)$, take
$X\in\mathcal{L}_{\mathcal{M}}(H,H_{t})\;(\simeq E(t))$ and $h,h^{\prime}\in
E(t)$ and note that $\langle T_{t}(X)h,h^{\prime}\rangle=\langle\iota
_{t}^{\ast}(X(h)),h^{\prime}\rangle=\langle X(h),\iota_{t}(h^{\prime})\rangle
$. Since $\iota_{t}(h^{\prime})\in H_{t}=\bigvee\{Y(h)\mid Y\in\mathcal{L}%
_{\mathcal{M}}(H,H_{t}),\;h\in H\}$, and since for $Y\in\mathcal{L}%
_{\mathcal{M}}(H,H_{t})$ and $k\in H$, $\langle X(h),Y(k)\rangle=\langle
h,\langle X,Y\rangle k\rangle$, the desired continuity is evident.

Finally, we must verify  equation (\ref{implement}). For $t=0$, the equation
is clear, so we always work with a fixed $t>0$. For $k\in H$, $\widetilde
{T}_{t}^{\ast}k\in\mathcal{L}_{\mathcal{M}}(H,H_{t})\otimes_{\mathcal{M}%
^{\prime}}H$ and so we may write $\widetilde{T}_{t}^{\ast}k=\sum X_{i}\otimes
h_{i},\;X_{j}\in\mathcal{L}_{\mathcal{M}}(H,H_{t})$, and $h_{j}\in H$. Then,
for $Y\in\mathcal{L}_{\mathcal{M}}(H,H_{t})$ and $h\in H$, $\langle\sum
X_{i}(h_{i}),Y(h)\rangle=\langle\sum Y^{\ast}X_{i}h_{i},h\rangle=\langle\sum
X_{i}\otimes h_{i},Y\otimes h\rangle=\langle\widetilde{T}_{t}^{\ast}k,Y\otimes
h\rangle=\langle k,\iota_{t}^{\ast}(Y(h))\rangle=\langle\iota_{t}%
(k),Y(h)\rangle$. Since $\bigvee\{Y(h)\mid Y\in\mathcal{L}_{\mathcal{M}%
}(H,H_{t}),\;h\in H\}=H_{t}$, this equation implies that $\sum X_{i}%
(h_{i})=\iota_{t}(k)$. Consequently, $\widetilde{T}_{t}(I\otimes
S)\widetilde{T}_{t}^{\ast}k=\widetilde{T}_{t}(I\otimes S)\sum X_{i}\otimes
h_{i}=\widetilde{T}_{t}(\sum X_{i}\otimes Sh_{i})=\iota_{t}^{\ast}(\sum
X_{i}(Sh_{i}))$. Since the $X_{i}$ are $\mathcal{M}$-module maps, this last
expression equals $\iota_{t}^{\ast}(S\sum X_{i}(h_{i}))=\iota_{t}^{\ast}%
S\iota_{t}(k)$. To evaluate $\iota_{t}^{\ast}S\iota_{t}(k)$, recall that
$\iota_{t}$ is the natural embedding of $H$ into $H_{t}$, when $H_{t}$ is
viewed as the inductive limit of the $H_{\mathfrak{p},t}$. Using the relations
among $\iota_{t}$, the $v_{0,\mathfrak{p,p}^{\prime}}$ and the
$v_{0,\mathfrak{p,\infty}}$ established above, it suffices choose a partition
$\mathfrak{p}\in\mathfrak{P}(t)$, and evaluate $\iota_{\mathfrak{p}}^{\ast
}S\iota_{\mathfrak{p}}(k)$. Suppose, then that $\mathfrak{p}=\{0=t_{0}%
<t_{1}<t_{2}<\cdots<t_{n-1}<t_{n}=t\}$. Then $S\iota_{\mathfrak{p}%
}(k)=S(I\otimes I\otimes\cdots\otimes I\otimes k)=S\otimes I\otimes
\cdots\otimes I\otimes k$ and $\iota_{\mathfrak{p}}^{\ast}S\iota
_{\mathfrak{p}}(k)=P_{t-t_{n-1}}(P_{t_{n-1}-t_{n-2}}(\cdots(P_{t_{1}%
}(S)I)\cdots)I)k=P_{t}(S)k$, by the semigroup property of $\{P_{t}\}_{t\geq0}%
$. Thus $\iota_{\mathfrak{p}}^{\ast}S\iota_{\mathfrak{p}}(k)=P_{t}(S)k$ for
all $\mathfrak{p}\in\mathfrak{P}(t)$, and so, in the limit, $\iota_{t}^{\ast
}S\iota_{t}(k)=P_{t}(S)k$.

Of course, setting $S=I$ in equation (\ref{implement}) shows that
$\{T_{t}\}_{t\geq0}$ is fully coisometric.
\end{proof}

This result is really not special to our given Markov semigroup $\{P_{t}%
\}_{t\geq0}$; the next result is a converse which shows that completely
contractive representations of product systems of correspondences over a von
Neumann algebra $\mathcal{N}$ always define completely positive semigroups on
$\mathcal{N}^{\prime}$.

\begin{theorem}
\label{Converse}Let $\mathcal{N}$ be a von Neumann algebra acting on a Hilbert
space $H$, let $\{E(t)\}_{t\geq0}$ be a discrete product system of
$\mathcal{N}$-correspondences, and let $\{T_{t}\}_{t\geq0}$ be a completely
contractive representation of $\{E(t)\}_{t\geq0}$ on $H$. For $S\in
\mathcal{N}^{\prime}$ and $t\geq0$, define
\[
\Theta_{t}(S)=\widetilde{T}_{t}(1_{E(t)}\otimes S)\widetilde{T}_{t}^{\ast
}\text{.}%
\]
Then $\{\Theta_{t}\}_{t\geq0}$ is a semigroup of normal, contractive,
completely positive maps on $\mathcal{N}^{\prime}$. Further, $\{\Theta
_{t}\}_{t\geq0}$ is unital if and only if $\{T_{t}\}_{t\geq0}$ is fully
coisometric, and $\{\Theta_{t}\}_{t\geq0}$ is a semigroup of $\ast
$-endomorphisms if $\{T_{t}\}_{t\geq0}$ is isometric.
\end{theorem}

\begin{proof}
Most of the result is proved in Proposition \ref{Proposition 1.15}. We simply
need to note that $T_{0}$ is a normal $\ast$-representation of $\mathcal{N}$
on $H$ and that for each $t\geq0$, $(T_{t},T_{0})$ is a completely contractive
covariant representation of $E(t)$ on $H$. All that really requires attention
is the fact that $\{\Theta_{t}\}_{t\geq0}$ is a semigroup, i.e., that
$\Theta_{t+s}=\Theta_{t}\circ\Theta_{s}$. However, the multiplicativity of
$\{T_{t}\}_{t\geq0}$ implies that for $s,t\geq0$, $\widetilde{T}%
_{t+s}=\widetilde{T}_{t}(I_{E(t)}\otimes\widetilde{T}_{s})$ and from this we
see immediately that $\Theta_{t+s}=\Theta_{t}\circ\Theta_{s}$.
\end{proof}

We note in passing that if the $\Theta_{t}$ are multiplicative, then by
Proposition \ref{Proposition 1.15}, the $E(t)$ decompose into the direct sum
$E(t)=E(t)^{\prime}\oplus E(t)^{\prime\prime}$ so that $T_{t}|E(t)^{\prime}$
is isometric, while $T_{t}|E(t)^{\prime\prime}$ is zero. The multiplicativity
of the $\Theta_{t}$ forces relations among the $q_{t}$, where $q_{t}$ is the
projection of $E(t)$ onto $E(t)^{\prime}$, but we shall not dwell on these here.

In the presence of Theorems \ref{Theorem2.6} and \ref{Lemma2.8}, we are able
to state and prove our dilation result, which is a semigroup analogue of
Theorem \ref{Theorem1.17}.

\begin{theorem}
\label{DiscDilat}Let $\mathcal{M}$ be a von Neumann algebra acting on the
Hilbert space $H$ and let $\{P_{t}\}_{t\geq0}$ be a semigroup of normal,
unital, completely positive maps on $\mathcal{M}$ such that $P_{0}$ is the
identity mapping on $\mathcal{M}$. Further, let $\{E(t)\}_{t\geq0}$ be the
product system of $\mathcal{M}^{\prime}$-correspondences associated with
$\{P_{t}\}_{t\geq0}$, let $\{T_{t}\}_{t\geq0}$ be the identity representation
of $\{E(t)\}_{t\geq0}$ on $H$ and let $(\{V_{t}\}_{t\geq0},u_{0},K)$ be the
minimal isometric dilation of $\{T_{t}\}_{t\geq0}$. We write $\rho$ for
$V_{0}$, thereby obtaining a normal $\ast$-homomorphism of $\mathcal{M}%
^{\prime}$ into $\mathcal{B}(K)$ and we set $\mathcal{R}$ equal to
$\rho(\mathcal{M}^{\prime})^{\prime}$. Then $u_{0}^{\ast}\mathcal{R}%
u_{0}=\mathcal{M}$, and if we define $\{\alpha_{t}\}_{t\geq0}$ by the formula%
\[
\alpha_{t}(S)=\widetilde{V}_{t}(I_{E(t)}\otimes S)\widetilde{V}_{t}^{\ast
}\text{,}%
\]
$S\in\mathcal{R}$, $t\geq0$, then $\{\alpha_{t}\}_{t\geq0}$ is a semigroup of
unital, normal, $\ast$-endomorphisms of $\mathcal{R}$ such that for $t\geq0$,
$S\in\mathcal{R}$, and $T\in\mathcal{M}$,%
\begin{equation}
P_{t}(u_{0}^{\ast}Su_{0})=u_{0}^{\ast}\alpha_{t}(S)u_{0} \label{Pt1}%
\end{equation}
and%
\begin{equation}
P_{t}(T)=u_{0}^{\ast}\alpha_{t}(u_{0}Tu_{0}^{\ast})u_{0}\text{.} \label{Pt2}%
\end{equation}
\end{theorem}

\begin{proof}
Since $\{V_{t}\}_{t\geq0}$ is a completely contractive covariant
representation of $\{E(t)\}_{t\geq0}$ on $K$, we know that $V_{0}=\rho$ is a
normal $\ast$-homomorphism of $\mathcal{M}^{\prime}$ on $K$. Further, from
Theorem \ref{Converse}, with $\mathcal{N}=\mathcal{M}^{\prime}$, we see that
$\{\alpha_{t}\}_{t\geq0}$ is a semigroup of normal $\ast$-endomorphisms of
$\mathcal{R}$ ($=\rho(\mathcal{M}^{\prime})^{\prime}$) that are unital because
$\{V_{t}\}_{t\geq0}$ is fully coisometric.

By the definition of $V_{0}$ in the proof of Theorem \ref{Theorem2.6} we see
that $\rho(a)u_{s}(\eta\otimes h)=V_{0}(a)u_{s}(\eta\otimes h)=u_{s}%
(a\eta\otimes h)$, for all $a\in\mathcal{M}^{\prime}$, $s\geq0$, and
$\eta\otimes h\in E(s)\otimes H$. This implies that the range of each $u_{s}$
reduces $\rho(\mathcal{M}^{\prime})$ and, in particular that the restriction
of $\rho(\mathcal{M}^{\prime})$ to $u_{0}(H)$ is unitarily equivalent to the
identity representation of $\mathcal{M}^{\prime}$. Specifically, since
$\rho(a)u_{0}(h)=u_{0}(ah)$, for all $a\in\mathcal{M}^{\prime}$, $h\in H$,
$a=u_{0}^{\ast}\rho(a)u_{0}$, $a\in\mathcal{M}^{\prime}$. Thus, $\mathcal{M}%
^{\prime}=u_{0}^{\ast}\rho(\mathcal{M}^{\prime})u_{0}$, so by the double
commutant theorem, $\mathcal{M}=u_{0}^{\ast}\rho(\mathcal{M}^{\prime}%
)^{\prime}u_{0}=u_{0}^{\ast}\mathcal{R}u_{0}$. Moreover, from Theorem
\ref{Lemma2.8} and equation (\ref{implement}), we see that for all
$S\in\mathcal{R}$,
\begin{align*}
P_{t}(u_{0}^{\ast}Su_{0}) &  =\widetilde{T}_{t}(I_{E(t)}\otimes u_{0}^{\ast
}Su_{0})\widetilde{T}_{t}^{\ast}\\
&  =u_{0}^{\ast}\widetilde{V}_{t}(I_{E(t)}\otimes u_{0})(I_{E(t)}\otimes
u_{0}^{\ast}Su_{0})(I_{E(t)}\otimes u_{0}^{\ast})\widetilde{V}_{t}^{\ast}%
u_{0}\\
&  =u_{0}^{\ast}\widetilde{V}_{t}(I_{E(t)}\otimes u_{0}u_{0}^{\ast}Su_{0}%
u_{0}^{\ast})\widetilde{V}_{t}^{\ast}u_{0}\\
&  =u_{0}^{\ast}u_{0}u_{0}^{\ast}\widetilde{V}_{t}(I_{E(t)}\otimes u_{0}%
u_{0}^{\ast})(I\otimes S)(I_{E(t)}\otimes u_{0}u_{0}^{\ast})\widetilde{V}%
_{t}^{\ast}u_{0}u_{0}^{\ast}u_{0}\\
&  =u_{0}^{\ast}u_{0}u_{0}^{\ast}\widetilde{V}_{t}(I\otimes S)\widetilde
{V}_{t}^{\ast}u_{0}u_{0}^{\ast}u_{0}\\
&  =u_{0}^{\ast}u_{0}u_{0}^{\ast}\alpha_{t}(S)u_{0}u_{0}^{\ast}u_{0}\\
&  =u_{0}^{\ast}\alpha_{t}(S)u_{0}%
\end{align*}
where the second and fifth equations are justified by Remark \ref{alternate}
and where the fourth and sixth equations are justified by the fact that the
final projection of $u_{0}$ is $u_{0}u_{0}^{\ast}$. Equation (\ref{Pt2}) can
be verified similarly, or directly from equation (\ref{Pt1}).
\end{proof}

\section{Minimality and Continuity}

Our goal in this section is to show that under the hypothesis of separability
on the Hilbert space $H$ and the hypothesis of weak continuity on
$\{P_{t}\}_{t\geq0}$ in Theorem \ref{DiscDilat}, the Hilbert space $K$ that is
produced there is separable and the semigroup $\{\alpha_{t}\}_{t\geq0}$ is
weakly continuous. That is, $\{\alpha_{t}\}_{t\geq0}$ will be an $E_{0}%
$-semigroup. Therefore, throughout this section, we make the blanket
assumption that our underlying Hilbert space $H$ is \emph{separable} and that
our semigroup of normal, unital, completely positive maps $\{P_{t}\}_{t\geq0}$
on $\mathcal{M}$ is (weakly) continuous in the sense that for all
$T\in\mathcal{M}$ and vectors $h_{1},h_{2}\in H$, the function $t\rightarrow
\langle P_{t}(T)h_{1},h_{2}\rangle$ is continuous. Note that since the
$\sigma$-weak topology coincides with the weak topology on bounded subsets of
a von Neumann algebra, our continuity assumption on $\{P_{t}\}_{t\geq0}$ is
tantamount to assuming that $(\mathcal{M},\{P_{t}\}_{t\geq0})$ is a quantum
Markov process. Our arguments will be broken into a series of (somewhat
technical) lemmas and propositions. Basically, we will distill for our use
salient features of \cite{wA97, wA97a} and \cite{dSL}. The first proposition
is a generalization of some observations due to D. SeLegue in \cite[Section
2.7]{dSL} when $\mathcal{M}=\mathcal{B}(H)$. We offer somewhat different proofs.

\begin{proposition}
\label{Proposition3.1}\ (\cite{dSL})Under our standing assumptions on
$\{P_{t}\}_{t\geq0}$, the following assertions are valid:

\begin{enumerate}
\item The map $t\rightarrow P_{t}(X)$ is strongly continuous for all
$X\in\mathcal{M}$; i.e., for all $h\in H$, $t\rightarrow P_{t}(X)h$ is
continuous from $[0,\infty)$ to $H$.

\item Given a sequence $\{X_{n}\}\subseteq\mathcal{M}$ that converges in the
weak operator topology to $X$, and given a sequence $\{t_{n}\}_{n=1}^{\infty
}\subseteq\lbrack0,\infty)$ converging to $t$, the sequence of operators
$\{P_{t_{n}}(X_{n})\}_{n=1}^{\infty}$ converges to $P_{t}(X)$ in the weak
operator topology, i.e., $\{P_{t}\}_{t\geq0}$ is jointly continuous in the
weak operator topology.
\end{enumerate}
\end{proposition}

\begin{proof}
For 1., first note that it suffices to prove the assertion when $X=U$ is
unitary and it suffices to show that $P_{t}(U)h\rightarrow Uh$ for every $h\in
H$ as $t\rightarrow0$. Then observe that for any vector $h\in H$,
$\lim_{t\rightarrow0}\left\|  P_{t}(U)h\right\|  =\left\|  h\right\|  $. For
if not, then the $\lim\inf\left\|  P_{t}(U)h\right\|  $ is strictly less than
$\left\|  h\right\|  $. Since $|\langle P_{t}(U)h,Uh\rangle|\leq\left\|
P_{t}(U)h\right\|  $, the $\lim\inf|\langle P_{t}(U)h,Uh\rangle|$ is strictly
less than $\left\|  h\right\|  $, also. However, by our hypothesis on
$\{P_{t}\}_{t\geq0}$, $\langle P_{t}(U)h,Uh\rangle\rightarrow\langle
Uh,Uh\rangle=\left\|  h\right\|  $. Thus $\lim_{t\rightarrow0}\left\|
P_{t}(U)h\right\|  $ must be $\left\|  h\right\|  $. But then we see that for
all $h\in H$, $\left\|  P_{t}(U)h-Uh\right\|  ^{2}=\langle P_{t}%
(U)h-Uh,P_{t}(U)h-Uh\rangle=\left\|  P_{t}(U)h\right\|  ^{2}%
-2\operatorname{Re}\langle Uh,P_{t}(U)h\rangle+\left\|  Uh\right\|  ^{2}$
tends to zero, as $t\rightarrow0$, as required.

For 2., observe that the normality of $P_{t}$ means that there is a unique
bounded map $\Psi_{t}$ such that $P_{t}=\Psi_{t}^{\ast}$. The uniqueness and
the fact that $\{P_{t}\}_{t\geq0}$ is a semigroup imply the same is true for
$\{\Psi_{t}\}_{t\geq0}$, i.e., $\Psi_{t+s}=\Psi_{t}\Psi_{s}$. The continuity
of $\{P_{t}\}_{t\geq0}$ in the weak operator topology and the fact that
$\{P_{t}\}_{t\geq0}$ is uniformly bounded imply that $\omega\circ P_{t}(X)$ is
continuous in $t$ for all $X\in\mathcal{M}$, and all $\omega\in\mathcal{M}%
_{\ast}$. If we write the pairing between $\mathcal{M}_{\ast}$ and
$\mathcal{M}$ by $\langle\cdot,\cdot\rangle$ as we shall, then this means that
$\langle\Psi_{t}(\omega),X\rangle$ is continuous in $t$ for all $X$; i.e.,
$\{\Psi_{t}\}_{t\geq0}$ is weakly continuous on $\mathcal{M}_{\ast}$. But
$\mathcal{M}_{\ast}$ is separable and so by \cite[Corollary 3.1.8]{BR79}, the
semigroup $\{\Psi_{t}\}_{t\geq0}$ is strongly continuous on $\mathcal{M}%
_{\ast}$, i.e., for all $\omega\in\mathcal{M}_{\ast}$, $\left\|  \Psi
_{t}(\omega)-\Psi_{s}(\omega)\right\|  \rightarrow0$, as $t\rightarrow s$.
This means, in particular, that if $\omega_{h}$ is the vector state associated
with the vector $h\in H$, then $\left\|  \omega_{h}\circ P_{t}-\omega_{h}\circ
P_{s}\right\|  \rightarrow0$ as $t\rightarrow s$. So, if $\{X_{n}%
\}_{n=1}^{\infty}$ is a sequence in $\mathcal{M}$ that converges weakly to
$X\in\mathcal{M}$, and if $t_{n}\rightarrow t$, then
\begin{multline*}
|\langle P_{t_{n}}(X_{n})h,h\rangle-\langle P_{t}(X)h,h\rangle|=|\omega
_{h}\circ P_{t_{n}}(X_{n})-\omega_{h}\circ P_{t}(X)|\\
\leq|\omega_{h}\circ P_{t_{n}}(X_{n})-\omega_{h}\circ P_{t}(X_{n}%
)|+|\omega_{h}\circ P_{t}(X_{n})-\omega_{h}\circ P_{t}(X)|\\
\leq\left\|  \omega_{h}\circ P_{t_{n}}-\omega_{h}\circ P_{t}\right\|  \left\|
X_{n}\right\|  +|\langle P_{t}(X_{n}-X)h,h\rangle|\text{.}%
\end{multline*}
Since the norms of the $X_{n}$ are uniformly bounded by the uniform
boundedness principle, this inequality shows that $P_{t_{n}}(X_{n})\rightarrow
P_{t}(X)$ in the weak operator topology, as required.
\end{proof}

\begin{proposition}
\label{Theorem3.2}Under our standing separability and continuity assumptions,
the Hilbert space $K$ in Theorem \ref{DiscDilat} is separable.
\end{proposition}

\begin{proof}
Recall from the proof of Theorem \ref{Theorem2.6} that $K$ is defined to be
the inductive limit $\underrightarrow{\lim}(E(t)\otimes_{T_{0}}H,u_{t,s})$.
Since the sequence of spaces, $\{E(n)\otimes_{T_{0}}H\}_{n\geq0}$, is cofinal
in $\{E(t)\otimes_{T_{0}}H\}_{t\geq0}$, it suffices to show each space
$E(t)\otimes_{T_{0}}H$ is separable. However, each space $E(t)$ is isomorphic
to $\mathcal{L}_{\mathcal{M}}(H,H_{t})$ by Lemma \ref{Lemma2.1}. So, if we can
show that $H_{t}$ is separable, then $E(t)$ will be separable in the $\sigma
$-topology (which is the same as the $\sigma$-weak topology by Proposition
\ref{Lemma 1.2}.) But then, of course, $E(t)\otimes_{T_{0}}H$ will be spanned
by a sequence $\{X_{n}\otimes h_{m}\}_{m,n\geq0}$, where the $X_{n}$ run
through a countable set that is dense in $E(t)$ in the $\sigma$-topology and
the $h_{m}$ run through a countable dense set of $H$, and so $E(t)\otimes
_{T_{0}}H$ will be separable. Thus we need to show $H_{t}$ is separable.

Now $H_{t}$ is, itself, an inductive limit $\underrightarrow{\lim
}(H_{\mathfrak{p},t},v_{0,\mathfrak{p,p}^{\prime}})$ where $\mathfrak{p}$ and
$\mathfrak{p}^{\prime}$ range over $\mathfrak{P}(t)$, and $\mathfrak{p}%
^{\prime}$ refines $\mathfrak{p}$. The normality of the $P_{t}$ enables one to
see that each $H_{\mathfrak{p},t}$ is separable and the weak continuity of
$\{P_{t}\}_{t\geq0}$ enables one to replace $\mathfrak{P}(t)$ with a countable
(but not, strictly speaking, cofinal) subset. From these two observations the
separability of $H_{t}$ follows. Here are the details.

To see that $H_{\mathfrak{p},t}$ is separable, first observe that
$\mathcal{M}\otimes_{P_{t}}H$ is separable for any $t$. For this, it suffices
to show that if $\{T_{n}\}_{n\geq0}$ is any sequence that is strongly dense in
the unit ball of $\mathcal{M}$ and if $\{h_{n}\}_{n\geq0}$ is any dense
sequence of vectors in $H$, then any decomposible tensor, $T\otimes h$, with
$T$ in the unit ball of $\mathcal{M}$, is in the closure of $\{T_{n}\otimes
h_{n}\}_{n\geq0}$. So, passing to subsequences, if necessary, assume that
$T_{n}\rightarrow T$ strongly and that $h_{n}\rightarrow h$. Then
$T_{n}\otimes h_{n}-T\otimes h=(T_{n}-T)\otimes h+T\otimes(h_{n}%
-h)+(T-T_{n})\otimes(h_{n}-h)$. However, $\left\|  (T_{n}-T)\otimes h\right\|
^{2}=\langle h,P_{t}((T_{n}-T)^{\ast}(T_{n}-T))h\rangle\rightarrow0$ because
$T_{n}\rightarrow T$ strongly and $P_{t}$ is normal. On the other hand,
$\left\|  T\otimes(h_{n}-h)\right\|  ^{2}=\langle h_{n}-h,P_{t}(T^{\ast
}T)(h_{n}-h)\rangle\rightarrow0$ because $h_{n}\rightarrow h$ in $H$. Finally,
since $h_{n}\rightarrow h$ in $H$ and since $T_{n}$, $T$, and their images
under $P_{t}$ are all bounded in norm by $1$, we see that $\left\|
(T-T_{n})\otimes(h-h_{n})\right\|  ^{2}=\langle h-h_{n},P_{t}((T_{n}-T)^{\ast
}(T_{n}-T))(h-h_{n})\rangle\leq4\left\|  h-h_{n}\right\|  ^{2}\rightarrow0$.
This shows that $T_{n}\otimes h_{n}\rightarrow T\otimes h$ as required. Now
the proof that each $H_{\mathfrak{p},t}$ is separable is proved by iterating
this argument.

Let $\mathfrak{P}_{0}(t)$ be the collection of those partitions $\mathfrak{p}%
\in\mathfrak{P}(t)$ whose points lie in $t\mathbb{Q}\cap\lbrack0,t]$. Observe
that $\mathfrak{P}_{0}(t)$ is countable and write $\widetilde{H}_{t}$ for the
union $\cup\{v_{0,\mathfrak{p,\infty}}(H_{\mathfrak{p},t})\mid\mathfrak{p}%
\in\mathfrak{P}_{0}(t)\}$. Then $\widetilde{H}_{t}$ is the countable union of
separable Hilbert spaces and so its closure is separable. We will show that
its closure is all of $H_{t}$. For this purpose, it suffices to show that if
$\mathfrak{p}$ is an arbitrary partition in $\mathfrak{P}(t)$, then
$v_{0,\mathfrak{p,\infty}}(H_{\mathfrak{p},t})$ is in the closure of
$\widetilde{H}_{t}$. This, in turn, will be clear if we can show that if
$\mathfrak{p}=\{0=t_{0}<t_{1}<t_{2}<\cdots<t_{n-1}<t_{n}=t\}$ and if
$\{\mathfrak{p}_{m}\}_{m\geq0}=\{\{0=s(m)_{0}<s(m)_{1}<s(m)_{2}<\cdots
<s(m)_{n-1}<s(m)_{n}=t\}\}_{m\geq0}$ is a sequence of partitions in
$\mathfrak{P}_{0}(t)$ such that $\lim_{m\rightarrow\infty}s(m)_{k}=t_{k}$ for
every $k$, then for every $n$-tuple $T_{1},T_{2},\cdots,T_{n}\in\mathcal{M}$
and every $h\in H$,
\begin{gather*}
\lim_{m\rightarrow\infty}v_{0,\mathfrak{p}_{m}\mathfrak{,\infty}}T_{1}%
\otimes_{P_{s(m)_{1}}}T_{2}\otimes_{P_{s(m)_{2}-s(m)_{1}}}\cdots
\otimes_{P_{t-s(m)_{n-1}}}h=\\
v_{0,\mathfrak{p,\infty}}T_{1}\otimes_{P_{t_{1}}}T_{2}\otimes_{P_{t_{2}-t_{1}%
}}\cdots\otimes_{P_{t-t_{n-1}}}h\text{.}%
\end{gather*}
To verify this equation, it suffices to assume that $s(m)_{k}=t_{k}$ for all
$m$ and for all $k$ but one. So, in fact, it is enough to verify the desired
limit when $\mathfrak{p}=\{0=t_{0}<t_{1}<t_{2}=t\}$ and when each
$\mathfrak{p}_{m}$ is of the form $\{0=s(m)_{0}<s(m)_{1}<s(m)_{2}=t\}$, where
$s(m)_{1}<t_{1}$. In this event, we have%
\begin{align*}
&  \left\|  v_{0,\mathfrak{p}_{m}\mathfrak{,\infty}}T_{1}\otimes_{P_{s(m)_{1}%
}}T_{2}\otimes_{P_{t-s(m)_{1}}}h-v_{0,\mathfrak{p,\infty}}T_{1}\otimes
_{P_{t_{1}}}T_{2}\otimes_{P_{t-t_{1}}}h\right\|  ^{2}\\
&  =\langle h,P_{t-s(m)_{1}}(T_{2}^{\ast}P_{s(m)_{1}}(T_{1}^{\ast}T_{1}%
)T_{2})h\rangle\\
&  -\langle h,P_{t-t_{1}}(P_{t_{1}-s(m)_{1}}(T_{2}^{\ast}P_{s(m)_{1}}%
(T_{1}^{\ast}T_{1})T_{2})h\rangle\\
&  -\langle h,P_{t-t_{1}}(T_{2}^{\ast}P_{t_{1}-s(m)_{1}}(P_{s(m)_{1}}%
(T_{1}^{\ast}T_{1})T_{2}))h\rangle\\
&  +\langle h,P_{t-t_{1}}(T_{2}^{\ast}P_{t_{1}}(T_{1}^{\ast}T_{1}%
)T_{2})h\rangle\text{.}%
\end{align*}
A moment's reflection reveals that assertion (2) in Proposition
\ref{Proposition3.1} shows that this expression tends to zero as
$m\rightarrow\infty$.
\end{proof}

Our next goal is to show that $\{\alpha_{t}\}_{t\geq0}$ is a minimal dilation
of $\{P_{t}\}_{t\geq0}$ in the sense of \cite{wA97a}. To explain this, recall
that a projection $q\in\mathcal{R}$ is \emph{increasing }relative to
$\{\alpha_{t}\}_{t\geq0}$ in case $\alpha_{t}(q)\geq q$ for all $t\geq0$,
i.e., the family $\{\alpha_{t}(q)\}_{t\geq0}$ is an increasing family of
projections. We will show that $u_{0}u_{0}^{\ast}$ is increasing. (Recall that
from the proof of Theorem \ref{DiscDilat}, $u_{0}u_{0}^{\ast}\in\mathcal{R}$.)
A projection $q\in\mathcal{R}$ is called \emph{multiplicative} in case the map
$X\rightarrow q\alpha_{t}(X)q$ is multiplicative on $\mathcal{R}$ for each
$t$, i.e., $q\alpha_{t}(\cdot)q$ is a (non-unital) endomorphism of
$\mathcal{R}$. To say, then, that $\{\alpha_{t}\}_{t\geq0}$ is \emph{minimal
}is to say that there is no multiplicative, increasing projection
$q\in\mathcal{R}$ that dominates $u_{0}u_{0}^{\ast}$, i.e., such that $q\geq
u_{0}u_{0}^{\ast}$.

We note in passing that in \cite{wA97, wA97a} Arveson assumes that
$\{\alpha_{t}\}_{t\geq0}$ is weakly continuous. However, this is not necessary
for the definition of minimality. That is, minimality makes sense without
assuming that $\{\alpha_{t}\}_{t\geq0}$ is weakly continuous. Here, minimality
will be used to show that the $\{\alpha_{t}\}_{t\geq0}$ we constructed in
Theorem \ref{DiscDilat} is weakly continuous.

\begin{lemma}
\label{Lemma3.3}With the notation of Section 2, we have, for all $t\geq0$,

\begin{enumerate}
\item $\widetilde{V}_{t}^{\ast}u_{0}=(I_{t}\otimes u_{0})\widetilde{T}%
_{t}^{\ast}$, where $I_{t}$ denotes the identity operator on $E(t)$.

\item $\bigvee\{(I_{t}\otimes X)\widetilde{T}_{t}^{\ast}h\mid X\in\mathcal{M}%
$, $h\in H\}=\mathcal{E}_{t}\otimes_{\mathcal{M}^{\prime}}H$.

\item $\bigvee\{(I_{t}\otimes Y)\widetilde{V}_{t}^{\ast}h\mid Y\in\mathcal{R}%
$, $h\in H\}=\mathcal{E}_{t}\otimes_{\mathcal{M}^{\prime}}K$.

\item $\bigvee\{\alpha_{t}(Y)u_{0}h\mid Y\in\mathcal{R}$, $h\in H\}=\bigvee
\{V_{t}(X)k\mid X\in\mathcal{E}_{t}$, $k\in K\}$
\end{enumerate}
\end{lemma}

\begin{proof}
(1) This is an easy consequence of Remark\ \ref{alternate}. Indeed, from the
first part of that remark, we know that $u_{0}^{\ast}\widetilde{V}%
_{t}(I\otimes u_{0})=\widetilde{T}_{t}$. So, $(I\otimes u_{0})^{\ast
}\widetilde{V}_{t}^{\ast}u_{0}=\widetilde{T}_{t}^{\ast}$. Therefore,
$(I\otimes u_{0}u_{0}^{\ast})\widetilde{V}_{t}^{\ast}u_{0}=(I\otimes
u_{0})\widetilde{T}_{t}^{\ast}$. On the other hand, the second part of the
remark shows that $(I\otimes u_{0}u_{0}^{\ast})\widetilde{V}_{t}^{\ast}%
u_{0}=(I\otimes u_{0}u_{0}^{\ast})\widetilde{V}_{t}^{\ast}(u_{0}u_{0}^{\ast
})u_{0}=\widetilde{V}_{t}^{\ast}(u_{0}u_{0}^{\ast})u_{0}=\widetilde{V}%
_{t}^{\ast}u_{0}$, so that $\widetilde{V}_{t}^{\ast}u_{0}=(I_{t}\otimes
u_{0})\widetilde{T}_{t}^{\ast}$.

(2) Here, $\mathcal{E}_{t}$ is embedded as a subspace of $E(t)$ and so
$\mathcal{E}_{t}\otimes H$ is contained in $E(t)\otimes H$. To show the
desired equality, first note that if $h\in H$ and if $\mathfrak{p}%
\in\mathfrak{P}(t)$, then $\iota_{\mathfrak{p}}(h)=I\otimes I\otimes
\cdots\otimes I\otimes h$. (See the discussion after Remark \ref{alternate}.)
If we let $\mathfrak{p}_{0}=\{0=t_{0}<t_{1}=t\}$, then in the notation
developed between Remark \ref{alternate} and Theorem \ref{Lemma2.8}, we have
$v_{0,\mathfrak{p}_{0},\mathfrak{p}}(I\otimes h)=\iota_{\mathfrak{p}}(h)$. So,
for $S\in\mathcal{M}$, $S\iota_{\mathfrak{p}}(h)=S\otimes I\otimes
\cdots\otimes I\otimes h=v_{0,\mathfrak{p}_{0},\mathfrak{p}}(S\otimes h)$. By
identifying $\mathcal{M}\otimes_{P_{t}}H$ with a subspace of $H_{\mathfrak{p}%
,t}$, we write $S\iota_{\mathfrak{p}}(h)\in\mathcal{M}\otimes_{P_{t}}H$,
$S\in\mathcal{M}$. Since this holds for all partitions $\mathfrak{p}%
\in\mathfrak{P}(t)$, we have
\begin{equation}
S\iota_{t}(h)\in\mathcal{M}\otimes_{P_{t}}H,\label{(i)}%
\end{equation}
for all $S\in\mathcal{M}$ and $t\geq0$. Now fix an element of $E(t)\otimes H$
that is orthogonal to $\mathcal{E}_{t}\otimes H$ and write it as $\sum
X_{i}\otimes h_{i}$, $X_{i}\in E(t)$, $h_{i}\in H$. Then for every
$X\in\mathcal{E}_{t}$ and $k\in H$, $0=\langle X\otimes k,\sum X_{i}\otimes
h_{i}\rangle=\sum\langle k,X^{\ast}X_{i}h_{i}\rangle=\langle X(k),\sum
X_{i}(h_{i})\rangle$. Since $\bigvee\{X(h)\mid X\in\mathcal{L}_{\mathcal{M}%
}(H,\mathcal{M}\otimes_{P_{t}}H)=\mathcal{E}_{t}$, $h\in H\}=\mathcal{M}%
\otimes_{P_{t}}H$, by Lemma \ref{lemma1.7}, we see that
\begin{equation}
\sum X_{i}(h_{i})\in(\mathcal{M}\otimes_{P_{t}}H)^{\perp}\text{.}\label{(ii)}%
\end{equation}
However, we have just shown above that $\mathcal{M}\iota_{t}(H)\subseteq
\mathcal{M}\otimes_{P_{t}}H$. Hence, for $S\in\mathcal{M}$ and $h\in H$,
\begin{multline*}
\langle(I\otimes S)\widetilde{T}_{t}^{\ast}h,\sum X_{i}\otimes h_{i}%
\rangle=\langle\widetilde{T}_{t}^{\ast}h,\sum X_{i}\otimes S^{\ast}%
h_{i}\rangle\\
=\langle h,\iota_{t}^{\ast}(\sum X_{i}(S^{\ast}h_{i}))\rangle=\langle
h,\iota_{t}^{\ast}(S^{\ast}\sum X_{i}(h_{i}))\rangle\\
=\langle S\iota_{t}(h),\sum X_{i}(h_{i})\rangle=0\text{,}%
\end{multline*}
where the last equation follows from (\ref{(i)}) and (\ref{(ii)}), and the
preceding one follows from the fact that the $X_{i}$ are $\mathcal{M}$-module
maps, i.e., they commute with elements of $\mathcal{M}$. This equation thus
shows that $[(I\otimes\mathcal{M})\widetilde{T}_{t}^{\ast}H]$ is contained in
$\mathcal{E}_{t}\otimes H$. For the reverse containment, fix an element $\sum
X_{i}\otimes h_{i}$, $X_{i}\in\mathcal{E}_{t}$, $h_{i}\in H$, that is in
$\mathcal{E}_{t}\otimes H$ but orthogonal to $[(I\otimes\mathcal{M}%
)\widetilde{T}_{t}^{\ast}H]$. Then for every $S\in\mathcal{M}$ and $h\in H$,
the last equation shows that $0=\langle(I\otimes S)\widetilde{T}_{t}^{\ast
}h,\sum X_{i}\otimes h\rangle=\langle S\iota_{t}(h),\sum X_{i}(h)\rangle$.
This shows that $\sum X_{i}(h)$ is orthogonal to $\mathcal{M}\otimes_{P_{t}}%
H$. Since $X_{i}\in\mathcal{E}_{t}$, $X_{i}(h_{i})\in\mathcal{M}\otimes
_{P_{t}}H$ for all $i$ and, so, $\sum X_{i}(h)=0$. However, $\langle\sum
X_{i}\otimes h_{i},\sum X_{j}\otimes h_{j}\rangle=\sum_{i,j}\langle
h_{i},X_{i}^{\ast}X_{j}h_{j}\rangle=\sum_{i,j}\langle X_{i}(h_{i}),X_{j}%
(h_{j})\rangle=\left\|  \sum X_{i}(h)\right\|  ^{2}=0$, and so $\sum
X_{i}\otimes h_{i}=0$ as was to be proved.

(3) From (1) we may write%
\[
(I\otimes\mathcal{R})\widetilde{V}_{t}^{\ast}u_{0}H=(I\otimes\mathcal{R}%
)(I\otimes u_{0})\widetilde{T}_{t}^{\ast}H=(I\otimes\mathcal{R}u_{0}%
)\widetilde{T}_{t}^{\ast}H\text{.}%
\]
As we noted in the proof of Theorem \ref{DiscDilat}, $u_{0}u_{0}^{\ast}$ lies
in $\mathcal{R}$. Consequently, $\mathcal{R}u_{0}=\mathcal{R}u_{0}u_{0}^{\ast
}\mathcal{R}u_{0}=\mathcal{R}u_{0}\mathcal{M}$. So, using (2), we conclude
that%
\begin{align*}
\lbrack(I\otimes\mathcal{R})\widetilde{V}_{t}^{\ast}u_{0}H]  &  =[(I\otimes
\mathcal{R}u_{0})(I\otimes\mathcal{M})(I\otimes u_{0})\widetilde{T}_{t}^{\ast
}H]\\
&  =[(I\otimes\mathcal{R}u_{0})(\mathcal{E}_{t}\otimes H)]=\mathcal{E}%
_{t}\otimes\lbrack\mathcal{R}u_{0}H]\text{.}%
\end{align*}
Now let $p$ be the projection of $K$ onto $[\mathcal{R}u_{0}]$. Then
$p\in\mathcal{R}^{\prime}=\rho(\mathcal{M}^{\prime})^{\prime\prime}%
=\rho(\mathcal{M}^{\prime})$; i.e. $p=\rho(p_{0})$ for some projection
$p_{0}\in\mathcal{M}^{\prime}$. However, $pK$ contains $u_{0}H$ and
$\rho(p_{0})$ acts on $u_{0}H$ by $\rho(p_{0})u_{0}h=u_{0}p_{0}h$. Hence
$p_{0}=I$ and so $p=I$. Thus $[\mathcal{R}u_{0}H]=K$ and so $[(I\otimes
\mathcal{R)}\widetilde{V}_{t}^{\ast}u_{0}H]=\mathcal{E}_{t}\otimes K$.

(4) The last assertion follows from the previous one and the definition of
$\alpha_{t}$:%
\[
\lbrack\alpha_{t}(\mathcal{R})u_{0}H]=[\widetilde{V}_{t}(I\otimes
\mathcal{R})\widetilde{V}_{t}^{\ast}u_{0}H]=[\widetilde{V}_{t}(\mathcal{E}%
_{t}\otimes K)]=[V_{t}(\mathcal{E}_{t})K]\text{.}%
\]
\end{proof}

\begin{lemma}
\label{Lemma3.4}Let $q_{t}$ be the projection from $K$ onto $[\alpha
_{t}(\mathcal{R})u_{0}H]$. Then $q_{t}$ lies in $\mathcal{R}$ and $q_{t}%
\alpha_{t}(q_{s})$ is the projection onto $[V_{t+s}(\mathcal{E}_{t}%
\otimes\mathcal{E}_{s})K]$.
\end{lemma}

\begin{proof}
The previous lemma shows that $q_{t}$ is the projection of $K$ onto
$[V_{t}(\mathcal{E}_{t})K]$. Thus $q_{t}$ lies in $\mathcal{R}=\rho
(\mathcal{M}^{\prime})^{\prime}$. Also, the range of $q_{t}$, which is
$[\alpha_{t}(\mathcal{R})u_{0}H]$ is clearly invariant under $\alpha
_{t}(\mathcal{R})$; i.e., $q_{t}\in\alpha_{t}(\mathcal{R})^{\prime}$. Thus, in
particular, $q_{t}$ commutes with $\alpha_{t}(q_{s})$; so we see that
$q_{t}\alpha_{t}(q_{s})$ is a projection. We need to show that $q_{t}%
\alpha_{t}(q_{s})$ has range $[V_{t+s}(\mathcal{E}_{t}\otimes\mathcal{E}%
_{s})K]$. For this purpose, observe that the range of $\alpha_{t}(q_{s})$ is
$\widetilde{V}_{t}(I\otimes q_{s})\widetilde{V}_{t}^{\ast}K=\widetilde{V}%
_{t}(I\otimes q_{s})(E(t)\otimes K)=\widetilde{V}_{t}(E(t)\otimes
q_{s}(K))=\widetilde{V}_{t}(E(t)\otimes\lbrack V_{s}(\mathcal{E}%
_{s})K])=[V_{t}(E(t))V_{s}(\mathcal{E}_{s})K]$. Clearly, $[V_{t}%
(\mathcal{E}_{t})V_{s}(\mathcal{E}_{s})K]\subseteq\lbrack V_{t}(E(t))V_{s}%
(\mathcal{E}_{s})K]\cap\lbrack V_{t}(\mathcal{E}_{t})K]$. We claim that, in
fact, the two subspaces coincide. To see this, let $w$ be the isometric
embedding of $\mathcal{M}\otimes_{P_{t}}H$ into $H_{t}$, view $\mathcal{M}%
\otimes_{P_{t}}H$ as a subspace of $H_{t}$ and view $\mathcal{E}_{t}$ as a
subspace of $E(t)$ (i.e. omit reference to the canonical embeddings.) Also,
identify $\mathcal{L}_{\mathcal{M}}(H,H_{t})$ with $E(t)$ and $\mathcal{L}%
_{\mathcal{M}}(H,\mathcal{M}\otimes_{P_{t}}H)$ with $\mathcal{E}_{t}$, as we
have throughout this paper. Then the map $p$ on $E(t)=\mathcal{L}%
_{\mathcal{M}}(H,H_{t})$ defined by the formula $p(X)=ww^{\ast}\circ X$, $X\in
E(t)$, is a projection in $\mathcal{L}(E(t))$ with range $\mathcal{E}_{t}$.
For elements $V_{t}(X_{i})V_{s}(Y_{i})k_{i}$, $i=1,2$, in $[V_{t}%
(E(t))V_{s}(\mathcal{E}_{s})K]$, we have%
\begin{align*}
\langle V_{t}(X_{1})V_{s}(Y_{1})k_{1},V_{t}(pX_{2})V_{s}(Y_{2})k_{2}\rangle &
=\langle\widetilde{V}_{t}(X_{1}\otimes V_{s}(Y_{1})k_{1}),\widetilde{V}%
_{t}(p(X_{2})\otimes V_{s}(Y_{2})k_{2})\rangle\\
&  =\langle X_{1}\otimes V_{s}(Y_{1})k_{1},p(X_{2})\otimes V_{s}(Y_{2}%
)k_{2}\rangle\\
&  =\langle V_{s}(Y_{1})k_{1},\rho(X_{1}^{\ast}\circ ww^{\ast}\circ
X_{2})V_{s}(Y_{2})k_{2}\rangle\\
&  =\langle V_{s}(Y_{1})k_{1},\rho(p(X_{1})^{\ast}X_{2})V_{s}(Y_{2}%
)k_{2}\rangle\\
&  =\langle p(X_{1})\otimes V_{s}(Y_{1})k_{1},X_{2}\otimes V_{s}(Y_{2}%
)k_{2}\rangle\\
&  =\langle V_{t}(p(X_{1}))V_{s}(Y_{1})k_{1},V_{t}(X_{2})V_{s}(Y_{2}%
)k_{2}\rangle\text{.}%
\end{align*}
Thus, the map $V_{t}(X)V_{s}(Y)k\rightarrow V_{t}(p(X))V_{s}(Y)k$ is
selfadjoint and, therefore, is the orthogonal projection from $[V_{t}%
(E(t))V_{s}(\mathcal{E}_{s})K]$ onto $[V_{t}(\mathcal{E}_{t})V_{s}%
(\mathcal{E}_{s})K]$. Write $q$ for this projection. A similar computation
shows that the projection from $[V_{t}(E(t))K]$ onto $[V_{t}(\mathcal{E}%
_{t})K]$ is given by the formula $V_{t}(X)k\rightarrow V_{t}(p(X))k$. However,
this projection is just the restriction of $q_{t}$ to $[V_{t}(E(t))K]$. We can
then restrict $q_{t}$ further to $[V_{t}(E(t))V_{s}(\mathcal{E}_{s})K]$ and
then the restricted image will clearly be $[V_{t}(E(t))V_{s}(\mathcal{E}%
_{s})K]\cap\lbrack V_{t}(\mathcal{E}_{t})K]$ (because $\alpha_{t}(q_{s})$
commutes with $q_{t}$). Also, by the definition of $q$, this restriction of
$q_{t}$ is just $q$ and so its image is $[V_{t}(\mathcal{E}_{t})V_{s}%
(\mathcal{E}_{s})K]$. Thus the range of $q_{t}\alpha_{t}(q_{s})$ is
$[V_{t}(\mathcal{E}_{t})V_{s}(\mathcal{E}_{s})K]=[V_{t+s}(\mathcal{E}%
_{t}\otimes\mathcal{E}_{s})K]$.
\end{proof}

Now let $\mathfrak{p}=\{0=t_{0}<t_{1}<t_{2}<\cdots<t_{n-1}<t_{n}=t\}$ be a
partition in $\mathfrak{P}(t)$, write $q_{s}$ for the projection onto
$[\alpha_{s}(\mathcal{R})u_{0}H]$, as in the last lemma, and set%
\[
q_{\mathfrak{p},t}:=q_{t-t_{n-1}}\alpha_{t-t_{n-1}}(q_{t_{n-1}-t_{n-2}}%
)\cdots\alpha_{t_{2}-t_{1}}(q_{t_{1}})\text{.}%
\]
Repeated use of the last lemma shows that $q_{\mathfrak{p},t}(K)=[V_{t}%
(\mathcal{E}_{t-t_{n-1}}\otimes\mathcal{E}_{t_{n-1}-t_{n-2}}\otimes
\cdots\otimes\mathcal{E}_{t_{1}})K]=[V_{t}(\mathcal{L}_{\mathcal{M}%
}(H,H_{\mathfrak{p},t}))K]$. Thus, it is clear that the $q_{\mathfrak{p},t}$
increase as the partitions $\mathfrak{p}$ are refined and since $E(t)=\lim
\mathcal{L}_{\mathcal{M}}(H,H_{\mathfrak{p},t})$, we see that they converge
strongly to the projection onto $[V_{t}(E(t))K]$; call it $\overline{q}_{t}$.
However, since $\widetilde{V}_{t}$ is a coisometry, we see that $K=[\widetilde
{V}_{t}(E(t)\otimes K)]=[V_{t}(E(t))K]$. Thus, $\overline{q}_{t}=I$, $t\geq0$.

Observe that $\overline{q}_{t}$ is the same projection defined by Arveson in
Section 3 of \cite{wA97}. He uses a slightly different indexing scheme for the
partitions that enter into his $q_{\mathfrak{p},t}$, but a moment's reflection
reveals that his $q_{\mathfrak{p},t}$ are the same as ours.

\begin{proposition}
\label{Proposition3.4}The semigroup of endomorphisms, $\{\alpha_{t}\}_{t\geq
0}$, of $\mathcal{R}$ is minimal.
\end{proposition}

\begin{proof}
As Arveson indicates in \S3 of \cite{wA97} (see page 575, in particular),
since we have shown that the projections $\overline{q}_{t}$ are all equal to
$I$, it remains to show that $\alpha_{t}(u_{0}u_{0}^{\ast})\rightarrow I$, as
$t\rightarrow\infty$. However, for each $t\geq0$, $\alpha_{t}(u_{0}u_{0}%
^{\ast})$ is the projection onto $[\widetilde{V}_{t}(I\otimes u_{0}u_{0}%
^{\ast})\widetilde{V}_{t}^{\ast}K]=[\widetilde{V}_{t}(I\otimes u_{0}%
u_{0}^{\ast})E(t)\otimes K]=[\widetilde{V}_{t}(E(t)\otimes u_{0}u_{0}^{\ast
}K)]=[V_{t}(E(t))u_{0}H]=u_{t}(E(t)\otimes H)$, where, recall, $u_{t}$ is the
embedding of $E(t)\otimes H$ into $K$. Since the spaces $u_{t}(E(t)\otimes H)$
are nested and have span equal to $K$, we conclude that the projections
$\alpha_{t}(u_{0}u_{0}^{\ast})$ increase to $I$.
\end{proof}

Let $p_{+}$ be the projection of $K$ onto the span
\[
\bigvee\{\alpha_{t_{1}}(u_{0}a_{1}u_{0}^{\ast})\alpha_{t_{2}}(u_{0}a_{2}%
u_{0}^{\ast})\cdots\alpha_{t_{n}}(u_{0}a_{n}u_{0}^{\ast})u_{0}h\mid a_{i}%
\in\mathcal{M},h\in H,\text{and }t_{i}\geq0\}
\]
and let $\mathcal{R}_{+}$ be the von Neumann algebra generated by
$\{\alpha_{t}(u_{0}u_{0}^{\ast}\mathcal{R}u_{0}u_{0}^{\ast})\mid t\geq0\}$.
Then, as Arveson shows in Proposition 3.14 of \cite{wA97a}, $p_{+}$ is the
largest projection in the center of $\mathcal{R}_{+}$ that dominates
$u_{0}u_{0}^{\ast}$ and, as he shows in Theorem B of \cite{wA97a}, because
$\{\alpha_{t}\}_{t\geq0}$ is minimal, $p_{+}=I$. (Note: In the proof of
\cite[Theorem B]{wA97a}, Arveson assumes that $\mathcal{R}$ is a factor.
However, this assumption is not necessary for the implications spelled out
there that we have used.) Thus we have

\begin{corollary}
\label{Corollary3.5}The von Neumann algebra $\mathcal{R}$ is generated by
$\{\alpha_{t}(u_{0}u_{0}^{\ast}\mathcal{R}u_{0}u_{0}^{\ast})\mid t\geq0\}$ and
$K$ is the span $\bigvee\{\alpha_{t_{1}}(u_{0}a_{1}u_{0}^{\ast})\alpha_{t_{2}%
}(u_{0}a_{2}u_{0}^{\ast})\cdots\alpha_{t_{n}}(u_{0}a_{n}u_{0}^{\ast}%
)u_{0}h\mid a_{i}\in\mathcal{M}$, $h\in H$, and $t_{i}\geq0\}$.
\end{corollary}

We let $\mathcal{A}$ denote the $C^{\ast}$-algebra generated by $\{\alpha
_{t}(u_{0}u_{0}^{\ast}\mathcal{R}u_{0}u_{0}^{\ast})\mid t\geq0\}$. Then
$\mathcal{A}$ is a translation invariant $C^{\ast}$-subalgebra of
$\mathcal{R}$ that generates $\mathcal{R}$ as a von Neumann algebra. To show
that $\{\alpha_{t}\}_{t\geq0}$ is weakly continuous on $\mathcal{R}$ we show
first that it is weakly continuous on $\mathcal{A}$ and then promote the weak
continuity there to all of $\mathcal{R}$. For this purpose, we begin with the
following result proved by SeLegue \cite{dSL} in the context when
$\mathcal{M}=\mathcal{B}(H)$. Our proof is somewhat different.

\begin{proposition}
\label{Proposition3.6}(\cite[Proposition 2.27]{dSL}) For every $T\in
\mathcal{M}$, $\alpha_{t}(u_{0}Tu_{0}^{\ast})\rightarrow u_{0}Tu_{0}^{\ast}$
in the strong operator topology as $t\rightarrow0+$.
\end{proposition}

\begin{proof}
Fix $T\in\mathcal{M}$ and $k\in K$ and then%
\begin{align*}
\left\|  (\alpha_{t}(u_{0}Tu_{0}^{\ast})-u_{0}Tu_{0}^{\ast})k\right\|  ^{2} &
=\langle(\alpha_{t}(u_{0}T^{\ast}u_{0}^{\ast})-u_{0}T^{\ast}u_{0}^{\ast
})(\alpha_{t}(u_{0}Tu_{0}^{\ast})-u_{0}Tu_{0}^{\ast})k,k\rangle\\
&  =\langle\alpha_{t}(u_{0}T^{\ast}Tu_{0}^{\ast})k,k\rangle-\langle\alpha
_{t}(u_{0}T^{\ast}u_{0}^{\ast})u_{0}Tu_{0}^{\ast}k,k\rangle\\
&  -\langle\alpha_{t}(u_{0}T^{\ast}u_{0}^{\ast})u_{0}Tu_{0}^{\ast}%
k,k\rangle+\langle u_{0}T^{\ast}Tu_{0}^{\ast}k,k\rangle
\end{align*}
to realize that it suffices to show that $\alpha_{t}(u_{0}Tu_{0}^{\ast
})\rightarrow u_{0}Tu_{0}^{\ast}$ in the weak operator topology as
$t\rightarrow0+$ for every $T\in\mathcal{M}$. However, since we have shown
that $\{\alpha_{t}\}_{t\geq0}$ is minimal and since $\{\alpha_{t}\}_{t\geq0}$
is uniformly bounded, we may apply Corollary \ref{Corollary3.5} to assert that
it suffices to show that
\[
\langle\alpha_{t}(u_{0}Tu_{0}^{\ast})k_{1},k_{2}\rangle\rightarrow\langle
u_{0}Tu_{0}^{\ast}k_{1},k_{2}\rangle
\]
for all $T\in\mathcal{M}$ and all vectors $k_{i}$ of the form $\alpha_{t_{1}%
}(u_{0}a_{1}u_{0}^{\ast})\alpha_{t_{2}}(u_{0}a_{2}u_{0}^{\ast})\cdots
\allowbreak\alpha_{t_{n}}(u_{0}a_{n}u_{0}^{\ast})\allowbreak u_{0}h$, $h\in
H$, $a_{i}\in\mathcal{M}$, and $t_{i}\geq0$. Note, too, that if any $t_{j}=0$,
then $\alpha_{t_{j}}(u_{0}a_{1}u_{0}^{\ast})\alpha_{t_{j+1}}(u_{0}a_{2}%
u_{0}^{\ast})\cdots\allowbreak\alpha_{t_{n}}(u_{0}a_{n}u_{0}^{\ast})u_{0}h$
lies in $H$ and so we may assume for the discussion that every $t_{i}>0$.
Also, let $t$ be the minimal number among the $t_{i}$'s. Then we may write
$\alpha_{t_{1}}(u_{0}a_{1}u_{0}^{\ast})\alpha_{t_{2}}(u_{0}a_{2}u_{0}^{\ast
})\cdots\alpha_{t_{n}}(u_{0}a_{n}u_{0}^{\ast})u_{0}h$ as $\alpha_{t}%
(\cdots)u_{0}h$. That is, $\alpha_{t_{1}}(u_{0}a_{1}u_{0}^{\ast})\alpha
_{t_{2}}(u_{0}a_{2}u_{0}^{\ast})\cdots\alpha_{t_{n}}(u_{0}a_{n}u_{0}^{\ast
})u_{0}h$ is in the cyclic subspace $[\alpha_{t}(\mathcal{R})u_{0}h]$. Thus,
we may assume that $k_{1}=\alpha_{r}(R)u_{0}h_{1}$ and that $k_{2}=\alpha
_{s}(L)u_{0}h_{2}$, where $R$ and $L$ are in $\mathcal{R}$, the $h_{i}$ are in
$H$ and $r,s>0$. We need to show that for $T\in\mathcal{M}$,
\[
\langle\alpha_{t}(u_{0}Tu_{0}^{\ast})\alpha_{r}(R)u_{0}h_{1},\alpha
_{s}(L)u_{0}h_{2}\rangle\rightarrow\langle u_{0}Tu_{0}^{\ast}\alpha
_{r}(R)u_{0}h_{1},\alpha_{s}(L)u_{0}h_{2}\rangle
\]
as $t\rightarrow0+$. For this, we may assume at the outset that the $t$'s
under consideration are all less than $r$ and $s$. Further, since $\alpha
_{t}(u_{0}u_{0}^{\ast})\geq u_{0}u_{0}^{\ast}$ as we saw in the proof of
Theorem \ref{DiscDilat}, we find that
\begin{multline*}
\langle\alpha_{t}(u_{0}Tu_{0}^{\ast})\alpha_{r}(R)u_{0}h_{1},\alpha
_{s}(L)u_{0}h_{2}\rangle=\langle\alpha_{t}(\alpha_{s-t}(L^{\ast})u_{0}%
Tu_{0}^{\ast}\alpha_{r-t}(R))u_{0}h_{1},u_{0}h_{2}\rangle\\
=\langle\alpha_{t}(\alpha_{s-t}(L^{\ast})u_{0}Tu_{0}^{\ast}\alpha
_{r-t}(R))\alpha_{t}(u_{0}u_{0}^{\ast})u_{0}h_{1},\alpha_{t}(u_{0}u_{0}^{\ast
})u_{0}h_{2}\rangle\\
=\langle\alpha_{t}(u_{0}u_{0}^{\ast}\alpha_{s-t}(L^{\ast})u_{0}Tu_{0}^{\ast
}\alpha_{r-t}(R))u_{0}u_{0}^{\ast})u_{0}h_{1},u_{0}h_{2}\rangle\\
=\langle\alpha_{t}(u_{0}P_{s-t}(u_{0}^{\ast}L^{\ast}u_{0})u_{0}^{\ast}%
u_{0}Tu_{0}^{\ast}u_{0}P_{r-t}(u_{0}^{\ast}Ru_{0})u_{0}^{\ast})u_{0}%
h_{1},u_{0}h_{2}\rangle\\
=\langle P_{t}(P_{s-t}(u_{0}^{\ast}L^{\ast}u_{0})TP_{r-t}(u_{0}^{\ast}%
Ru_{0}))h_{1},h_{2}\rangle\text{.}%
\end{multline*}
By the first assertion in Proposition \ref{Proposition3.1}, the functions
$t\rightarrow P_{s-t}(u_{0}^{\ast}L^{\ast}u_{0})$ and $t\rightarrow
P_{r-t}(u_{0}^{\ast}Ru_{0})$ are strongly continuous. Consequently, the
function $t\rightarrow P_{s-t}(u_{0}^{\ast}L^{\ast}u_{0})TP_{r-t}(u_{0}^{\ast
}Ru_{0})$ is weakly continuous. Therefore, applying the second assertion in
Proposition \ref{Proposition3.1}, we see that
\[
\langle P_{t}(P_{s-t}(u_{0}^{\ast}L^{\ast}u_{0})TP_{r-t}(u_{0}^{\ast}%
Ru_{0}))h_{1},h_{2}\rangle\rightarrow\langle u_{0}Tu_{0}^{\ast}\alpha
_{r}(R)u_{0}h_{1},\alpha_{s}(L)u_{0}h_{2}\rangle,
\]
completing the proof.
\end{proof}

Let $\mathcal{R}_{0}:=\{R\in\mathcal{R}\mid\lim_{t\rightarrow0+}\alpha
_{t}(R)=R$ in the strong operator topology$\}$. Then, since $\{\alpha
_{t}\}_{t\geq0}$ is a semigroup of $\ast$-endomorphisms of $\mathcal{R}$,
$\mathcal{R}_{0}$ is easily seen to be a $\ast$-subalgebra of $\mathcal{R}$.
In fact, since $\left\|  \alpha_{t}(R)k-Rk\right\|  \leq\left\|  \alpha
_{t}(S)k-Sk\right\|  +\left\|  \alpha_{t}(R-S)k-(R-S)k\right\|  \leq2\left\|
R-S\right\|  \left\|  k\right\|  +\left\|  \alpha_{t}(S)k-Sk\right\|  $, we
see that any $R$ in the norm closure of $\mathcal{R}_{0}$ already is in
$\mathcal{R}_{0}$. Thus, $\mathcal{R}_{0}$ is a $C^{\ast}$-algebra. This
$C^{\ast}$-algebra contains $u_{0}\mathcal{M}u_{0}^{\ast}$ by the preceding
proposition. But also, since each $\alpha_{r}$ is a normal $\ast$-endomorphism
of $\mathcal{R}$, we see that $\mathcal{R}_{0}$ contains $\alpha_{r}%
(u_{0}\mathcal{M}u_{0}^{\ast})$ for all $r\geq0$. Indeed, the proposition
shows that for each $T\in\mathcal{M}$, $\alpha_{t}(u_{0}Tu_{0}^{\ast
})\rightarrow u_{0}Tu_{0}^{\ast}$ strongly as $t\rightarrow0+$. Therefore,
since $\alpha_{r}$ is normal, $\alpha_{r}(\alpha_{t}(u_{0}Tu_{0}^{\ast
}))\rightarrow\alpha_{r}(u_{0}Tu_{0}^{\ast})$ strongly, as $t\rightarrow0+$.
Since $\alpha_{r}(\alpha_{t}(u_{0}Tu_{0}^{\ast}))=\alpha_{t}(\alpha_{r}%
(u_{0}Tu_{0}^{\ast}))$, we see that $\alpha_{t}(\alpha_{r}(u_{0}Tu_{0}^{\ast
}))\rightarrow\alpha_{r}(u_{0}Tu_{0}^{\ast})$ strongly, as $t\rightarrow0+$.
Thus, $\mathcal{R}_{0}\supseteq\mathcal{A}$ and so $\mathcal{R}_{0}$ is weakly
dense in $\mathcal{R}$ by Corollary \ref{Corollary3.5}. We are therefore well
on our way to showing that $\mathcal{R}_{0}=\mathcal{R}$. For this, we need
the following lemma.

\begin{lemma}
\label{lemma3.7}The $\alpha_{t}$, for \emph{strictly positive} $t$, are
jointly faithful on $\mathcal{R}$, i.e.,
\[
\cap\{\ker(\alpha_{t})\mid t>0\}=\{0\}.
\]
\end{lemma}

\begin{proof}
The kernel of each $\alpha_{t}$ is a $2$-sided, $\sigma$-weakly closed ideal
in $\mathcal{R}$. Thus so is $\cap\{\ker(\alpha_{t})\mid t>0\}$. Hence, we may
write $\cap\{\ker(\alpha_{t})\mid t>0\}=q\mathcal{R}$ for some central
projection $q$ in $\mathcal{R}$. Since for $R\in\mathcal{A}$ we have
$\alpha_{t}(R)\rightarrow R$ strongly as $t\rightarrow0+$, we conclude that
$\mathcal{A}\cap q\mathcal{R}=\{0\}$. Since $\mathcal{A}$ generates
$\mathcal{R}$ as a von Neumann algebra by Corollary \ref{Corollary3.5}, we
conclude that $q=0$.
\end{proof}

We have arrived at the main theorem of the paper.

\begin{theorem}
\label{theorem3.8}Let $(\mathcal{M},\{P_{t}\}_{t\geq0})$ be a quantum Markov
process and assume that $\mathcal{M}$ acts on a separable Hilbert space $H$.
Then the discrete dilation $(K,\mathcal{R},\{\alpha_{t}\}_{t\geq0},\allowbreak
u_{0})$ constructed from $\{P_{t}\}_{t\geq0}$ in Theorem \ref{DiscDilat} is an
$E_{0}$-dilation; i.e., $\{\alpha_{t}\}_{t\geq0}$ is weakly continuous.
\end{theorem}

\begin{proof}
By Proposition \ref{Theorem3.2}, $K$ is separable and so the predual of
$\mathcal{R}$, $\mathcal{R}_{\ast}$, is separable as a Banach space. We will
write the pairing between $\mathcal{R}$ and $\mathcal{R}_{\ast}$ as
$\langle\cdot,\cdot\rangle$, i.e., $\langle\omega,R\rangle=\omega(R)$,
$\omega\in\mathcal{R}_{\ast}$, $R\in\mathcal{R}$, as we did in Proposition
\ref{Proposition3.1}. However, here we write $\Psi_{t}$ for the pre-adjoint of
$\alpha_{t}$, i.e., $\Psi_{t}(\omega)=\omega\circ\alpha_{t}$ for all
$\omega\in\mathcal{R}_{\ast}$; so $\langle\Psi_{t}(\omega),R\rangle
=\langle\omega,\alpha_{t}(R)\rangle$. Since the $\sigma$-weak topology on
$\mathcal{B}(K)$ agrees with the weak operator topology on bounded subsets, we
see from the discussion following Proposition \ref{Proposition3.6} that for
all $\omega\in\mathcal{R}_{\ast}$ and all $R\in\mathcal{R}_{0}$ the function
$t\rightarrow\langle\omega,\alpha_{t}(R)\rangle$ is continuous. However, if
$R\in\mathcal{R}$ we may find a sequence $\{R_{n}\}_{n=1}^{\infty}$ in
$\mathcal{R}_{0}$ that converges weakly to $R$. Consequently, the function
$t\rightarrow\langle\omega,\alpha_{t}(R)\rangle$ is the pointwise limit of the
sequence of continuous functions $t\rightarrow\langle\omega,\alpha_{t}%
(R_{n})\rangle$. Therefore $t\rightarrow\langle\omega,\alpha_{t}(R)\rangle$ is
measurable for each $\omega\in\mathcal{R}_{\ast}$ and each $R\in\mathcal{R}$.
That is, the semigroup of linear maps $\{\Psi_{t}\}_{t\geq0}^{\infty}$ on
$\mathcal{R}_{\ast}$ is weakly measurable with respect to the duality between
$\mathcal{R}_{\ast}$ and $\mathcal{R}$ (See \cite[Definition 3.5.4]{HP74}.)
Since $\mathcal{R}_{\ast}$ is separable, Theorem 3.5.3 of \cite{HP74} implies
that $t\rightarrow\Psi_{t}(\omega)$ is strongly measurable as an
$\mathcal{R}_{\ast}$-valued function. Thus, in the terminology of
\cite[Chapter 10]{HP74}, $\{\Psi_{t}\}_{t\geq0}^{\infty}$ is a strongly
measurable semigroup of linear maps on $\mathcal{R}_{\ast}$. But then, Theorem
10.2.3 of \cite{HP74} shows that at least for $t$ strictly larger than zero,
the function $t\rightarrow\Psi_{t}$ is strongly continuous; i.e., for each
$\omega\in\mathcal{R}_{\ast}$ the $\mathcal{R}_{\ast}$-valued function on
$(0,\infty)$, $t\rightarrow\Psi_{t}(\omega)$, is continuous with respect to
the norm topology on $\mathcal{R}_{\ast}$. To extend the continuity to all of
$[0,\infty)$, let $\widetilde{\mathcal{R}}_{\ast}$ be the closed linear span
$\bigvee\{\Psi_{t}(\mathcal{R}_{\ast})\mid t>0\}$. If $\widetilde{\mathcal{R}%
}_{\ast}$ is not all of $\mathcal{R}_{\ast}$, then there is an $R\in
\mathcal{R}$ such that $\langle\omega,R\rangle=0$ for all $\omega\in
\widetilde{\mathcal{R}}_{\ast}$. This means that for all $t>0$ and all
$\omega\in\mathcal{R}_{\ast}$, $\langle\omega,\alpha_{t}(R)\rangle=\langle
\Psi_{t}(\omega),R\rangle=0$. Thus, $R$ is in the kernels of all the
$\alpha_{t}$. However, Lemma \ref{lemma3.7} implies that $R=0$. Thus,
$\widetilde{\mathcal{R}}_{\ast}$ is all of $\mathcal{R}_{\ast}$. Now we can
appeal to Theorem 10.5.5 of \cite{HP74} to conclude that $\lim_{\rightarrow
0+}\left\|  \Psi_{t}(\omega)-\omega\right\|  =0$. Consequently, for all
$\omega\in\mathcal{R}_{\ast}$ and $R\in\mathcal{R}$ we see that $\langle
\omega,\alpha_{t}(R)\rangle=\langle\Psi_{t}(\omega),R\rangle\rightarrow
\langle\omega,R\rangle$ as $t\rightarrow0+$, which is what we wanted to prove.
\end{proof}

\end{document}